\documentclass[ejs]{imsart}

\RequirePackage{amsthm,amsmath,amsfonts,amssymb}
\RequirePackage[numbers]{natbib}
\RequirePackage{bm}
\RequirePackage{bbm}
\RequirePackage{float}
\RequirePackage[linesnumbered,ruled,vlined,spanish,onelanguage]{algorithm2e}

\RequirePackage[colorlinks,citecolor=blue,urlcolor=blue]{hyperref}
\RequirePackage{graphicx}

\newcommand{\ind}{\mathbbm 1}

\arxiv{2010.00000; Accepted in 10/23/2024; To Appear in Electronic Journal of Statistics}
\startlocaldefs
\theoremstyle{plain}

\newtheorem{theorem}{Theorem}[section]
\newtheorem{lemma}[theorem]{Lemma}
\newtheorem{proposition}[theorem]{Proposition}
\newtheorem{corollary}[theorem]{Corollary}
\theoremstyle{definition}
\newtheorem{definition}[theorem]{Definition}

\theoremstyle{remark}

\endlocaldefs

\begin{document}
\begin{frontmatter}
\title{Tail-adaptive Bayesian shrinkage}
\runtitle{Tail-adaptive Bayesian shrinkage}

\begin{aug}
\author[A]{\fnms{Se Yoon}~\snm{Lee}\ead[label=e1]{seyoonlee.stat.math@gmail.com}},
\author[A]{\fnms{Peng}~\snm{Zhao}\ead[label=e2]{pzhao@tamu.edu }},
\author[A]{\fnms{Debdeep}~\snm{Pati}\ead[label=e3]{debdeep@stat.tamu.edu }},
\and
\author[A]{\fnms{Bani K.}~\snm{Mallick}\ead[label=e4]{bmallick@stat.tamu.edu}}
\address[A]{Department of Statistics, Texas A\&M University, College Station, 3143 TAMU, TX, 77843, USA \printead[presep={,\ }]{e1}}

\runauthor{Lee et al.}
\end{aug}

\begin{abstract}
Robust Bayesian methods for high-dimensional regression problems under diverse sparse regimes are studied. Traditional shrinkage priors are primarily designed to detect a handful of signals from tens of thousands of predictors in the so-called ultra-sparsity domain. However, they may not perform desirably when the degree of sparsity is moderate. In this paper, we propose a robust sparse estimation method under diverse sparsity regimes, which has a tail-adaptive shrinkage property. In this property, the tail-heaviness of the prior adjusts adaptively, becoming larger or smaller as the sparsity level increases or decreases, respectively, to accommodate more or fewer signals, a posteriori. We propose a global-local-tail (GLT) Gaussian mixture distribution that ensures this property. We examine the role of the tail-index of the prior in relation to the underlying sparsity level and demonstrate that the GLT posterior contracts at the minimax optimal rate for sparse normal mean models. We apply both the GLT prior and the Horseshoe prior to a real data problem and simulation examples. Our findings indicate that the varying tail rule based on the GLT prior offers advantages over a fixed tail rule based on the Horseshoe prior in diverse sparsity regimes.
\end{abstract}



\end{frontmatter}

\section{Introduction}
The advancement of sophisticated data acquisition techniques in gene expression microarrays, as well as other fields, has spurred the development of innovative statistical methods \citep{friedman2001elements,buhlmann2011statistics,hastie2015statistical}. One of the crucial goals of statistical analysis in these studies is to identify relevant predictors associated with a response variable from a large pool of predictors, despite having a limited number of samples. We consider the framework of sparse linear regression, where the data is given by:
\begin{align}
	\label{eq:sparse_high_dim_linear_regression}
	\textbf{y} &= \textbf{X}\bm{\beta} + \sigma \bm{\epsilon}, \quad
	\bm{\epsilon} \sim \mathcal{N}_{n}(\textbf{0},\textbf{I}_{n}).
\end{align}

Here, $\textbf{y} \in \mathbb{R}^{n}$ represents the $n$-dimensional response vector, $\textbf{X} \in \mathbb{R}^{n \times p}$ is a deterministic $n$-by-$p$ design matrix, $\bm{\beta} \in \mathbb{R}^{p}$ is the $p$-dimensional coefficient vector, and $\bm{\epsilon} \in \mathbb{R}^{n}$ denotes the $n$-dimensional Gaussian noise vector. Specifically, we are interested in the sparse setup, where $n \leq p$, and possibly $n \ll p$, with a significant proportion of the coefficients in $\bm{\beta} = (\beta_{1}, \cdots, \beta_{p})^{\top}$ are zero or close to zero. The non-zero coefficients in $\bm{\beta}$ represent signals, while the remaining coefficients represent noise.

To facilitate the discussion, we define the sparsity level as follows throughout the paper. Let $q$ denote the number of true signals among the $p$ coefficients. The sparsity level is a ratio representing the proportion of true signals out of all coefficients.
\begin{align}
	\label{eq:sparsity_level_gene_expression_data}
	s = \frac{q}{p} = \frac{\text{the number of relevant predictors}}{\text{total number of predictors}}.
\end{align}
In most real data applications, the sparsity level $s$ (\ref{eq:sparsity_level_gene_expression_data}) is an unknown quantity because the truth of the $\bm{\beta}$ is unknown.

Much has been written and published on the penalized regression techniques for estimating $\bm{\beta}$ under the assumption of sparsity \citep{hastie2015statistical}. From a Bayesian perspective, sparsity-favoring mixture priors with separate control over the signal and noise coefficients have been proposed \citep{mitchell1988bayesian,george1995stochastic,johnson2010use,yang2016computational}. Although these priors often exhibit attractive theoretical properties \citep{castillo2012needles,castillo2015bayesian}, computational challenges and the consideration that many of the $\beta_j$'s may be small but not exactly zero have led to the development of a wide variety of continuous shrinkage priors \citep{park2008bayesian,tipping2001sparse,griffin2010inference,carvalho2010horseshoe,carvalho2009handling}, which can be unified through a global-local scale mixture representation \citep{polson2010shrink}. Among the continuous shrinkage priors, the horseshoe prior \citep{carvalho2010horseshoe,carvalho2009handling}, or simply the Horseshoe, is arguably one of the most acclaimed methods. It is important to note that the posterior obtained using the Horseshoe exhibits nice finite sample performance and possesses several optimal theoretical properties  when the underlying sparsity level is very small, in the so-called ultra-sparse regime \citep{van2016many,van2017adaptive,van2017uncertainty}.

In the literature, several research works have found that the Horseshoe estimator may approach the zero vector under certain circumstances \citep{bai2018beta, van2017uncertainty}. This has garnered significant attention from practitioners utilizing the Horseshoe for many applications because once the Horseshoe estimator vanishes, it becomes difficult to infer the coefficients of the regression problem. However, little research has been published to explore or theoretically address this issue, despite its importance in the suitable application of shrinkage priors to a variety of sparse setups. We term this phenomenon the `collapsing behavior' of the Horseshoe estimator, typically caused by a sharp underestimation of the global-scale parameter. In this paper, we empirically demonstrate that the Horseshoe estimator may collapse under a moderately sparse regime, where the degree of sparsity is not so small. We show that this collapse is related to the restricted tail behavior of the Horseshoe, where the tail index is fixed. A similar phenomenon has also been observed by \cite{bai2018beta}. Importantly, in the recent discussion by \citep{van2017uncertainty}, \citep{yoo2017contributed} drew attention to the potential collapse of the marginal maximum-likelihood estimator for the global-scale parameter under an ultra-sparse regime. However, our research focuses on the collapsing behavior of the fully Bayesian Horseshoe estimator under a moderately sparse regime, addressing a different problem from the previous findings of \citep{van2017uncertainty,yoo2017contributed}.

In this paper, we develop a robust shrinkage prior that works reasonably well across diverse sparse domains. A key idea to achieve this goal is to adapt the posterior tail behavior to the sparsity level $s$ (\ref{eq:sparsity_level_gene_expression_data}). We refer to this as tail-adaptive shrinkage property, characterized by the following:

\begin{itemize}
	\item[(i)] Under an ultra sparse regime (where the sparsity level $s$ (\ref{eq:sparsity_level_gene_expression_data}) is very small), the tail-heaviness of a shrinkage prior adaptively gets thinner to accommodate a small number of signals, \emph{a posteriori}.
	\item[(ii)]  Under a moderately sparse regime (where the sparsity level $s$ (\ref{eq:sparsity_level_gene_expression_data}) is not too small), the tail-heaviness of a shrinkage prior adaptively gets thicker to accommodate a larger number of signals, \emph{a posteriori}.
\end{itemize}

In order to possess this desirable property, we first propose a new family of continuous shrinkage priors called global-local-tail shrinkage priors, and then introduce its member, `the GLT prior.' This new shrinkage formulation can be considered a generalization of the existing global-local shrinkage priors. Theoretical demonstrations show that the posterior for the GLT contracts at the (near) minimax optimal rate under sparse normal mean models, emphasizing the crucial role of the tail-heaviness parameter for robust estimation in diverse sparse regimes. We also developed a Markov Chain Monte Carlo (MCMC) sampling algorithm for GLT posterior computation. This algorithm is designed to adapt to the unknown sparsity level and is facilitated by the combination of a modern MCMC sampler and a peaks-over-threshold method from extreme value theory \citep{embrechts2013modelling,lee2018exponentiated}.



The rest of the paper is organized as follows. In Section \ref{sec:Fixed tail rule-- Global-local shrinkage priors}, we briefly illustrate global-local shrinkage priors and the Horseshoe, and explore the restricted tail-heaviness of the Horseshoe. In Section \ref{sec:Varying tail rule-- Global-local-tail shrinkage priors}, we propose global-local-tail shrinkage priors to release the tail-heaviness of existing shrinkage priors and introduce the GLT prior as a member of this new prior formulation. Theoretical properties of prior analysis and posterior convergence for the GLT prior are explored in Section \ref{sec:Properties of the GLT prior}. In Section \ref{sec:Example-- prostate cancer data}, we apply the GLT prior and the Horseshoe to a real-world data problem. Additionally, in Section \ref{sec:Example-- Simulation study with varied sparsity level}, we investigate the two priors through a simulation study. Replicated numerical studies are conducted in Section \ref{sec:Simulations} to compare the performance of the Horseshoe and the GLT prior. Section \ref{sec:Discussion} contains a summary and discussion of the findings.

\section{Fixed tail rule-- Global-local shrinkage priors}\label{sec:Fixed tail rule-- Global-local shrinkage priors}
\subsection{\textbf{Tail-heaviness of a density}}\label{subsec:Definition of tail-heaviness of a density}
Throughout the paper, the notion of tail-heaviness is adopted from extreme value theory and regular variation. This concept applies to Lebesgue measurable functions, distribution functions, prior distributions, posterior distributions, and so on. \citep{karamata1933mode,maric2000regular,mikosch1999regular,gierz2003encyclopedia,
	embrechts2013modelling,lee2018exponentiated}. 
\begin{definition}\label{definition:regularly_varying function}
	A positive, Lebesgue measurable function $\rho$ on $(0,\infty)$ is regularly varying of index $\alpha \in \mathbb{R}$ if there exists $\alpha$ such that
	$\lim_{x \rightarrow \infty}
	\rho(c x)/\rho(x) = c^{-\alpha},
	\, \text{for any } c > 0$. If $\alpha=0$, then the function $\rho$ is said to be slowly varying.
\end{definition}


We explain how Definition \ref{definition:regularly_varying function} can be adapted to extreme value theory. Consider a positive random variable $X \sim F$, where $F$ is the distribution function of $X$. In this case, we replace the measurable function $\rho$ in Definition \ref{definition:regularly_varying function} with the tail (survival) function of $X$, denoted as $\bar{F} = 1 - F$. This modification leads to the limiting equation $\lim_{x \rightarrow \infty}
\bar{F}(c x)/\bar{F}(x) = c^{-\alpha}$ for any $c > 0$. According to Karamata's characterization theorem \citep{karamata1933mode}, it holds that $\bar{F}(x) = L(x) \cdot x^{-\alpha}$, where $L$ is a slowly varying function. 

In extreme value theory, the parameter $\alpha$ represents the tail-heaviness of the random variable $X$ and is referred to as the tail-index of $X$ or the tail-index of the density $f = F'$. Its reciprocal, denoted as $\xi = 1/\alpha$, is called the shape parameter \citep{coles2001introduction,embrechts2013modelling}. A distribution $F$ with a positive $\xi>0$ is called a heavy-tailed distribution (see page 268 of \citep{mcneil2015quantitative}). As the value of $\xi$ increases, the tail-heaviness of the density $f$ also increases. It is important to note that while we illustrated the notion of the tail-index and shape parameter using a positive random variable $X$, this concept can be generalized to a real-valued random variable in a similar fashion by considering the tail-index and shape parameter at either $\infty$ or $-\infty$.

\subsection{\textbf{Global-local shrinkage priors and the Horseshoe}}\label{subsec:Global-local shrinkage priors}
Most of the continuous shrinkage priors proposed and studied in the literature can be represented as global-local scale mixtures of Gaussian distributions
\begin{align}
	\label{eq:local-global_beta}
	\beta_j | \lambda_j, \tau, \sigma^{2} &\sim \mathcal{N}_{1}(0, \lambda_j^2 \tau^2 \sigma^2), \quad \sigma^{2} \sim h(\sigma^{2}), \quad (j = 1, \cdots, p),\\
	\label{eq:local-global_lambda_tau}
	\lambda_j &\sim f(\lambda_j) ,\quad
	\tau  \sim g(\tau), \quad (j = 1, \cdots, p),
\end{align}
where $f$, $g$, and $h$ are densities supported on $(0, \infty)$. Different choices of $f$ and $g$ for the top-level scale parameters lead to different classes of priors \citep{bhadra2017lasso}. In the high-dimensional setting, the choices of $f$ and $g$ play a key role in controlling the effective sparsity and concentration of the prior and posterior distributions \citep{polson2010shrink, pati2014posterior,song2017nearly,martin2017empirical,bai2018high, zhang2019ultra}.

The Horseshoe \citep{carvalho2010horseshoe} can be obtained by choosing the unit-scaled half-Cauchy densities, $\mathcal{C}^{+}(x|0,1) = 2/\{\pi(1 + x^{2})\},\, x>0$, for the $f$ and $g$ in (\ref{eq:local-global_lambda_tau}) under the global-local form (\ref{eq:local-global_beta}) -- (\ref{eq:local-global_lambda_tau}):
\begin{align}
	\label{eq:HS_beta}
	\beta_j | \lambda_j, \tau, \sigma^{2} &\sim \mathcal{N}_{1}(0, \lambda_j^2 \tau^2 \sigma^2), \quad \sigma^{2} \sim \pi(\sigma^{2})\propto 1/\sigma^{2} , \quad (j = 1, \cdots, p),\\
	\label{eq:HS_lambda_tau}
	\lambda_j &\sim \mathcal{C}^{+}(0,1) ,\quad
	\tau  \sim \mathcal{C}^{+}(0,1), \quad (j = 1, \cdots, p).
\end{align}


Among the continuous shrinkage priors, the Horseshoe (\ref{eq:HS_beta}) -- (\ref{eq:HS_lambda_tau}) \citep{carvalho2010horseshoe} is possibly the most studied member in recent literature. It is known that the Horseshoe estimator possesses many nice theoretical properties when considering the sparsity assumption $s \to 0$ as $n$ and $p \to \infty$ (sparsity level $s$ is defined in (\ref{eq:sparsity_level_gene_expression_data})). For instance, the Horseshoe estimator is robust and achieves the minimax-optimal rate for squared error loss, up to a multiplicative constant, under certain conditions \citep{polson2010shrink,bhadra2017lasso,song2017nearly,van2017uncertainty,bai2018high}. Several highly scalable algorithms have been recently proposed for the Horseshoe \citep{bhattacharya2016fast,johndrow2017bayes}.

\subsection{\textbf{Restricted tail-heaviness of the Horseshoe}}\label{sec:Restricted tail-heaviness of the Horseshoe}
In the following, we demonstrate that the tail-heaviness of the Horseshoe \eqref{eq:HS_beta}--\eqref{eq:HS_lambda_tau} is fixed in the sense of regular variation. For simplicity, we consider a univariate form of the Horseshoe, given as $\beta|\lambda, \tau \sim \mathcal{N}_{1}(0, \lambda^{2}\tau^{2})$, and $\lambda \sim \mathcal{C}^{+}(0,1)$ with fixed $\tau>0$. Then, the marginal density of the Horseshoe is $\pi_{\text{HS}}(\beta|\tau) = \int \mathcal{N}_{1}(\beta|0, \lambda^{2}\tau^{2}) \mathcal{C}^{+}(\lambda|0,1) \, d\lambda$. (Refer to \eqref{eq:marginal_beta_HS} in the Appendix for the closed-form expression.)

The following theorem states that the tail-index of this marginal density is fixed for any value of the global-scale parameter:
\begin{proposition}\label{corollary:tail_index_of_Horseshoe_prior}
	Assume $\beta | \lambda, \tau \sim \mathcal{N}_{1}(0, \lambda^{2}\tau^{2})$, $\lambda \sim C^{+}(0,1)$, and $\tau>0$. Then the tail-index of $\pi_{\text{HS}}(\beta | \tau)$ is $\alpha = 1$ for any $\tau > 0$. 
\end{proposition}

Proposition \ref{corollary:tail_index_of_Horseshoe_prior} is proven in Subsection \ref{sec:Proof-- Restricted tail-heaviness of Horseshoe} in the Appendix. In general, it is well known that the shape parameter of the half-Cauchy density is $\xi = 1$ \citep{embrechts2013modelling}. Proposition \ref{corollary:tail_index_of_Horseshoe_prior} implies that the tail-heaviness of the marginal density $\pi_{\text{HS}}(\beta | \tau)$ inherits that of the local-scale density $f(\lambda) = \mathcal{C}^{+}(\lambda|0,1)$ \eqref{eq:HS_lambda_tau} and is fixed for any value of $\tau > 0$. This suggests that although the marginal density $\pi_{\text{HS}}(\beta | \tau)$ is heavy-tailed (due to the positive value of $\xi$), this heaviness is constant regardless of the sparsity level, as also pointed out by \cite{piironen2017sparsity}.

The absence of a tail-controlling mechanism in the Horseshoe can pose challenges when dealing with various sparsity regimes. While the Horseshoe may perform well in ultra-sparse regimes, where the tail of the density $\pi_{\text{HS}}(\beta|\tau)$ induced by $\xi = 1$ is sufficiently thick to accommodate a small number of signals, challenges may arise as the sparsity level $s$ (\ref{eq:sparsity_level_gene_expression_data}) increases. The need to allocate more mass in the tail region may also grow to accommodate a moderate number of signals as the sparsity level increases. In such cases, the tail-heaviness induced by $\xi = 1$ may not be large enough to selectively capture the moderate number of signals, which can be problematic with weakly identified parameters, such as the regression coefficients in high-dimensional regression when the signal-to-noise ratio is too small or multicollinearity is present.

We note that this limitation does not stem from the specific choice of the local-scale density but rather from the prior formulation of global-local shrinkage priors, where only local and global scale parameters are utilized as random variables to achieve optimal shrinkage. Next, we extend the existing prior formulation to incorporate the tail-index parameter as a random variable.



\section{Varying tail rule -- Global-local-tail shrinkage priors}\label{sec:Varying tail rule-- Global-local-tail shrinkage priors}
\subsection{\textbf{Global-local-tail shrinkage priors}}\label{subsec:Global-local-tail shrinkage priors}
We propose a new hierarchical formulation of continuous shrinkage priors called global-local-tail shrinkage priors for the high-dimensional regression (\ref{eq:sparse_high_dim_linear_regression}). As the name suggests, these priors can be represented as a global-local-tail Gaussian mixture distribution.
\begin{align}
	\label{eq:glt_prior_beta}
	\beta_{j} | \lambda_{j}, \sigma^2 &\sim \mathcal{N}_{1}(0, \lambda_{j}^2 \sigma^{2}), \quad \sigma^2 \sim h (\sigma^2), \quad (j = 1, \cdots, p),\\
	\label{eq:glt_prior_lambda}
	\lambda_{j} | \tau, \xi &\sim f(\lambda_{j}|\tau, \xi),
	\quad (j = 1, \cdots, p), \\
	\label{eq:glt_prior_tau_xi}
	(\tau, \xi) &\sim g(\tau,  \xi),
\end{align}
where the $f$ is a Fréchet class-density supported on $(0,\infty)$ parameterized by scale and shape parameters, $\tau$ and $\xi$ (we detail it shortly later). The density function $h$ is supported on $(0, \infty)$, while the joint density function $g$ is supported on $(0, \infty) \times (0, \infty)$. Following the literature on continuous shrinkage priors \citep{bhadra2017lasso} and extreme value theory \citep{embrechts2013modelling}, the scale parameters $\lambda_j$ ($j=1,\cdots,p$) and $\tau$ are referred to as the local-scale parameters and the global-scale parameter, respectively, and $\xi$ is called the shape parameter. Throughout the paper, we use the Jeffreys prior \citep{jeffreys1946invariant} for the measurement error $\sigma^2$, which is denoted as $h(\sigma^2) \propto 1/\sigma^2$.


The overall shrinkage performance of the global-local-shrinkage priors (\ref{eq:glt_prior_beta})--(\ref{eq:glt_prior_tau_xi}) highly depends on the choice of the local-scale density $f$ (\ref{eq:glt_prior_lambda}). Well-designed priors can control the tail-heaviness of the marginal prior for $\beta_{j}$ by changing the value of $\xi$. To ensure the statistical validity of this hierarchy, the local-scale density $f=F'$ (\ref{eq:glt_prior_lambda}) should satisfy the condition that its corresponding distribution function $F$ belongs to the Fréchet class (as defined on page 268 of \citep{embrechts2011quantitative}):
\begin{align}
	\label{eq:Frechet_class_distributions}
	\text{MDA}_{\text{Fréchet}} = \{F\,|\,\bar{F}(\lambda) =  L(\lambda) \cdot \lambda^{-1/\xi},\,\xi>0, \, L\, \text{is slowly varying} \}.
\end{align}

Table \ref{table:members_of_f} provides several members of the family $\text{MDA}_{\text{Fréchet}}$. All distributions in the table have support on $(0,\infty)$. The half-Cauchy and half-Levy distributions are derived from the half-$\alpha$-stable distribution, with fixed tail-index values of $\alpha=1$ and $\alpha=1/2$, respectively. Additional examples for $f$ can be found in \citep{hill1975simple, embrechts2013modelling}. Essentially, these distributions are regularly varying functions with the shape parameter $\xi$ (or equivalently, tail-index $\alpha=1/\xi$). It is important to note that the decay rate of the distribution $F\in \text{MDA}_{\text{Fréchet}}$ is primarily determined by the polynomial term $\lambda^{-1/\xi}$ as $\lambda$ approaches infinity.

\begin{table}[ht]
	\small
	\caption{\baselineskip=10pt 
		Unit scaled densities $f = F^{'}$ with $F \in \text{MDA}_{\text{Fréchet}}$}\label{table:members_of_f}
	\centering
	\begin{tabular}{p{2cm}cc}
		\hline
		& $f(\lambda|\tau = 1, \xi)$ & Shape parameter $\xi$ \\
		\hline
		Half-$\alpha$-stable distribution & non-closed form & $\xi$ \\
		Half-Cauchy distribution & $2\{\pi(1 + \lambda^{2})\}^{-1}$ & $1$ \\
		Half-Levy distribution & $\lambda^{-3/2} \exp\{-1/(2\lambda)\} /\sqrt{2\pi}$ & $2$ \\
		Loggamma distribution & $\{(1 + \lambda)^{-(1/\xi + 1)}\}/\xi$ & $\xi$ \\
		Generalized extreme value distribution & $\exp\ \{-(1 + \xi \lambda)^{-1/\xi}\}(1 + \xi \lambda)^{-(1/\xi+1)}$ & $\xi$ \\
		Generalized Pareto distribution & $(1 + \xi \lambda)^{-(1/\xi+ 1)}$ & $\xi$ \\
		\hline
	\end{tabular}
\end{table}

One of the practical implementational challenges in using the global-local-tail shrinkage prior formulation is the estimation of $\xi$ \citep{Armagan2010}. Specifically, the central problem here is determining which prior $\pi(\xi)$ should be used and how to estimate $\pi(\xi)$. Because $\xi$ is located within the hierarchy, farthest from the response vector $\mathbf{y}$, it is highly likely that information propagated from $\mathbf{y}$ to $\xi$ when passing through the hierarchy would be attenuated. A second problem is that even after a good prior is chosen for $\pi(\xi)$, the estimation of the parameter $\xi$ is analytically non-trivial because $\xi$ appears in the density $f$ through an exponent. Therefore, selecting a suitable prior $\pi(\xi)$ with an efficient sampling technique is crucial and requires special attention. 

\subsection{\textbf{Relationship with global-local shrinkage priors}}\label{subsec:Relationship with global-local shrinkage priors}
Overall, two formulations of shrinkage priors (global-local shrinkage priors (\ref{eq:local-global_beta})--(\ref{eq:local-global_lambda_tau}) and global-local-tail shrinkage priors (\ref{eq:glt_prior_beta})--(\ref{eq:glt_prior_tau_xi})) commonly belong to the one-group continuous shrinkage prior formulation \citep{park2008bayesian, Armagan2010, bhattacharya2015dirichlet, griffin2017hierarchical, ghosh2017asymptotic}. However, they are based on different tail rules: the \emph{fixed tail rule} versus the \emph{varying tail rule}. A prior with a varying tail rule may retain great flexibility in the shape of the marginal density for the coefficients $\bm{\beta}$. It is the shape parameter $\xi$ that provides a means to control the tail behavior of a prior, possibly adapting it to the sparsity level $s$ (\ref{eq:sparsity_level_gene_expression_data}).

It is important to note that the global-local hierarchy (\ref{eq:local-global_beta})--(\ref{eq:local-global_lambda_tau}) can be regarded as a special case of the global-local-tail hierarchy (\ref{eq:glt_prior_beta})--(\ref{eq:glt_prior_tau_xi}) with a fixed shape parameter $\xi$. To see this, first, choose a local-scale density $f(\cdot) = f(\cdot|\tau, \xi)$ in (\ref{eq:glt_prior_lambda}), then fix the shape parameter $\xi$ to a positive value, and finally, bring up the scale parameter $\tau$ under the coefficient $\beta_{j}$ by re-parameterizing $\lambda_{j}/\tau$ for each $j$. Therefore, some members of the global-local form (\ref{eq:local-global_beta})--(\ref{eq:local-global_lambda_tau}) \citep{carvalho2010horseshoe, bhadra2017lasso} can be thought of as members of the global-local-tail form (\ref{eq:glt_prior_beta})--(\ref{eq:glt_prior_tau_xi}), with the only difference being whether the tail-heaviness is fixed in advance or not.

For example, the prior formulation of the Horseshoe (\ref{eq:HS_beta})--(\ref{eq:HS_lambda_tau}) can be derived from the global-local-tail form (\ref{eq:glt_prior_beta})--(\ref{eq:glt_prior_tau_xi}) as follows: First, choose the half-$\alpha$-stable density for the local-scale density $f$ with the scale $\tau$ distributed according to $\mathcal{C}^{+}(0,1)$ and shape $\xi$ in (\ref{eq:glt_prior_lambda}). Then, fix the shape parameter to be $\xi=1$, resulting in $f$ becoming the half-Cauchy density scaled by $\tau$. Finally, bring up $\tau$ under the coefficient $\beta_{j}$ by re-parameterizing $\lambda_{j}/\tau$ for each $j$.

\subsection{\textbf{The GLT prior}}\label{subsec:GLT prior}
We introduce a member of the global-local-tail form \eqref{eq:glt_prior_beta} -- \eqref{eq:glt_prior_tau_xi}, using the generalized Pareto distribution (GPD) \citep{pickands1975statistical} as the local-scale density:
\begin{align}
	\label{eq:gpd_beta}
	\beta_{j} | \lambda_{j}, \sigma^{2} &\sim \mathcal{N}_{1}(0, \lambda_{j}^2 \sigma^{2}), \quad \sigma^{2} \sim \pi (\sigma^{2}) \propto 1/\sigma^{2}, \quad (j = 1, \cdots, p),\\
	\label{eq:gpd_lambda}
	\lambda_{j}|\tau, \xi & \sim \mathcal{GPD} (\tau, \xi),\quad (j = 1,\cdots, p),\\
	\label{eq:gpd_tau}
	\tau | \xi & \sim \mathcal{IG} (p/\xi + 1, 1),\\
	\label{eq:gpd_xi}
	\xi &\sim \log\ \mathcal{N} (\mu, \rho^2)
	\mathcal{I}_{(1/2,\infty)}, \quad \mu \in \mathbb{R}, \, \rho^{2} > 0.
\end{align}
Here, the local-scale density is a GPD: $f(\lambda_{j})=\mathcal{GPD}(\lambda_{j}|\tau, \xi) = (1/\tau) \cdot (1 + \xi \lambda_{j} /\tau)^{-(1/\xi + 1)}$ for the $p$ local-scale parameters $\lambda_{j}$ ($j=1,\cdots,p$). The joint density of the global-scale parameter and shape parameter is given as a truncated inverse-gamma-lognormal joint density:
\[
g(\tau, \xi) = \mathcal{IG}(\tau| p/\xi + 1, 1) \cdot \mathcal{I}_{(0,\infty)}(\tau) \cdot \left\{\log\ \mathcal{N} (\xi|\mu, \rho^2) \cdot \mathcal{I}_{(1/2,\infty)}(\xi)\right\}/D,
\]
where $D = D(\mu, \rho^{2}) = \int_{1/2}^{\infty} \log\ \mathcal{N} (\xi|\mu, \rho^2) \, d\xi$ is the normalizer of $g(\tau, \xi)$.

The GPD is one of the frequently used distributions in the field of extreme value theory for modeling extreme values and tail behaviors of heavy-tailed data, where the tails of the distribution decay more slowly than the tails of a Gaussian distribution \citep{lee2018exponentiated, embrechts2013modelling}. The intuition behind using GPD as the local-scale density (as noted from (\ref{eq:gpd_lambda})) is that a signal coefficient induces an extreme value of the corresponding local-scale parameter, while a noise coefficient induces a local-scale parameter that is nearly zero. The GPD is employed to describe these extreme values while controlling the scale and shape parameters.

Noting from (\ref{eq:gpd_tau}), we use the unit-scaled inverse gamma prior $\mathcal{IG}(x|p/\xi + 1, 1)$ with a shape parameter of $p/\xi + 1$, depending on the number of predictors $p$. The main motivation is to set the prior mean $\mathbb{E}[\tau|\xi] = \xi/p$, resembling the form of the sparsity level $s=q/p$ (\ref{eq:sparsity_level_gene_expression_data}), where the numerator $q$ is the number of true signals. Noting from (\ref{eq:gpd_xi}), we use the log-normal distribution $\log \mathcal{N} (x|\mu, \rho^2)$ truncated on the interval $(1/2,\infty)$ as a prior for the shape parameter $\xi$ in (\ref{eq:gpd_xi}). The main motivation is to provide robust estimation of the shape parameter $\xi$ by offering a heavy-tailed distribution of the sub-exponential density \citep{gulisashvili2016tail, lee2022use}. The truncation of the support of the log-normal prior to $(1/2,\infty)$ is intended to produce a horseshoe shape for the density of random shrinkage coefficients \citep{carvalho2010horseshoe}.

We call this specific hierarchical form (\ref{eq:gpd_beta}) -- (\ref{eq:gpd_xi}) `the GLT prior', denoted as $\bm{\beta}\sim \pi_{\text{GLT}}(\bm{\beta})$. Note that $\mu$ and $\rho^{2}$ in (\ref{eq:gpd_xi}) are the hyper-parameters. A detailed explanation of the posterior computation can be found in Section \ref{sec:Posterior computation} of Appendix. The proposed sampling algorithm utilizes a Gibbs sampler \citep{casella1992explaining,lee2021gibbs} and incorporates automatic hyper-parameter tuning facilitated by a joint technique involving the elliptical slice sampler \citep{murray2010elliptical} and Hill estimator \citep{hill1975simple}.



In Section \ref{sec:Properties of the GLT prior}, we explore the theoretical properties of the GLT prior. Notably, one of its characteristics is that the GLT prior behaves like a two-group prior \citep{george1993variable, johnson2010use, johnson2012bayesian} despite being a continuous shrinkage prior. Furthermore, we demonstrate that the GLT posterior contracts at the minimax risk \citep{donoho1994minimax} under the sparse normal mean model.

\section{Properties of the GLT prior}\label{sec:Properties of the GLT prior}
\subsection{\textbf{Marginal density of the GLT prior}}\label{subsec:Marginal density of the GLT prior}
For simplicity, we work with a univariate form of the GLT prior (\ref{eq:gpd_beta}) -- (\ref{eq:gpd_xi}), given by 
$\beta|\lambda \sim \mathcal{N}_{1}(0, \lambda^{2}),$ and $\lambda|\tau, \xi \sim \mathcal{GPD}(\tau,\xi),
$
with fixed $\tau>0$ and $\xi>1/2$. The marginal densities of the coefficient and random shrinkage coefficient \citep{carvalho2010horseshoe}, conditioned on $(\tau, \xi)$, are given as follows:
\begin{proposition}
	\label{prop:marginal of glt shrinkage}
	Suppose $\beta|\lambda \sim \mathcal{N}_{1}(0, \lambda^{2})$, $\lambda \sim \mathcal{GPD}(\tau, \xi)$, $\tau>0$ and $\xi>1/2$. Then:
	\begin{itemize}
		\item[]\emph{(a)} density of $\beta$ given $\tau$ and $\xi$ is
		\begin{align}
			\label{eq:marginal_of_beta_glt_separate_expression}
			\pi(\beta |\tau, \xi )
			&=
			\sum_{k=0}^{\infty}
			a_{k}
			\{\psi_{k}^{\text{S}}(\beta)
			+
			\psi_{k}^{\text{R}}(\beta)
			\},
		\end{align}
		where $K =
		1/(\tau 2^{3/2} \pi^{1/2} )$, $Z(\beta) =
		\beta^{2} \xi^{2}/(2 \tau^{2})$, $a_{k} = (-1)^{k}
		\cdot
		K
		\cdot
		{1/\xi + k \choose k}$, $\psi_{k}^{\text{S}}(\beta) =
		E_{k/2 + 1}\{Z(\beta)\}$, and $\psi_{k}^{\text{R}}(\beta) =
		Z(\beta)^{-\frac{1 + 1/\xi + k}{2}}
		\gamma\{(1 + 1/\xi + k)/2 $ $
		, Z(\beta)\}$. The superscripts on $\psi_{k}^{\text{S}}$ and $
		\psi_{k}^{\text{R}}$ represent (noise) shrinkage and (tail) robustness, respectively.
		\item[]\emph{(b)} 
		density of $\kappa=1/(1+\lambda^{2})$ given $\tau$ and $\xi$ is
		\begin{align}
			\label{eq:marginal_of_kappa_GLT}
			\pi(\kappa | \tau, \xi) 
			=
			\frac{\tau^{1/\xi} }{2}
			\cdot
			\frac{\kappa^{1/(2\xi) -1  }
				(1 - \kappa)^{-1/2}}{\{
				\tau \kappa^{1/2} + \xi (1-\kappa)^{1/2}
				\}^{(1 + 1/\xi)}}.
		\end{align}
	\end{itemize}
\end{proposition}

Proposition \ref{prop:marginal of glt shrinkage} is proven in Subsection \ref{subsec:Proof-- Proposition Marginal of GLT} in Appendix. It is important to note that the marginal density $\pi(\beta|\tau, \xi)$ (\ref{eq:marginal_of_beta_glt_separate_expression}) is expressed analytically as an alternating series, with its summands separated into two terms: $\{\psi_{k}^{\text{S}}(\beta)\}_{k=0}^{\infty}$ and $\{\psi_{k}^{\text{R}}(\beta)\}_{k=0}^{\infty}$. This density involves two special functions: (i) the \emph{generalized exponential-integral function} of real order \citep{milgram1985generalized, chiccoli1992concerning}, denoted as $E_{s}(x) = \int_{1}^{\infty} e^{-xt} t^{-s} dt$ $(x > 0, s \in \mathbb{R})$, and (ii) the \emph{incomplete lower gamma function} denoted as $\gamma(s, x) = \int_{0}^{x} t^{s-1} e^{-t} dt$ $(s, x \in \mathbb{R})$. These special functions, $E_{s}(x)$ and $\gamma(s, x)$, are involved in the expression of $\pi(\beta|\tau, \xi)$ through the sequences of functions $\{\psi_{k}^{\text{S}}(\beta)\}_{k=0}^{\infty}$ and $\{\psi_{k}^{\text{R}}(\beta)\}_{k=0}^{\infty}$, respectively. The generalized binomial coefficient ${1/\xi + k \choose k}$ is defined as $(1/\xi + k)(1/\xi + k-1)\cdots (1/\xi + 1)/k!$ when $k \in \{1,2, \cdots\}$, and zero when $k = 0$.

Analytically, the marginal density of the Horseshoe $\pi_{\text{HS}}(\beta|\tau)$ discussed in Subsection~\ref{sec:Restricted tail-heaviness of the Horseshoe} is influenced by a single specific function, namely the \emph{exponential integral function} $E_{1}(x)$. The Horseshoe lacks a functional component induced from the \emph{incomplete lower gamma function} $\gamma(s, x)$, which is present in the GLT prior.

\begin{figure}[h]
	\centering
	\includegraphics[width=\textwidth]{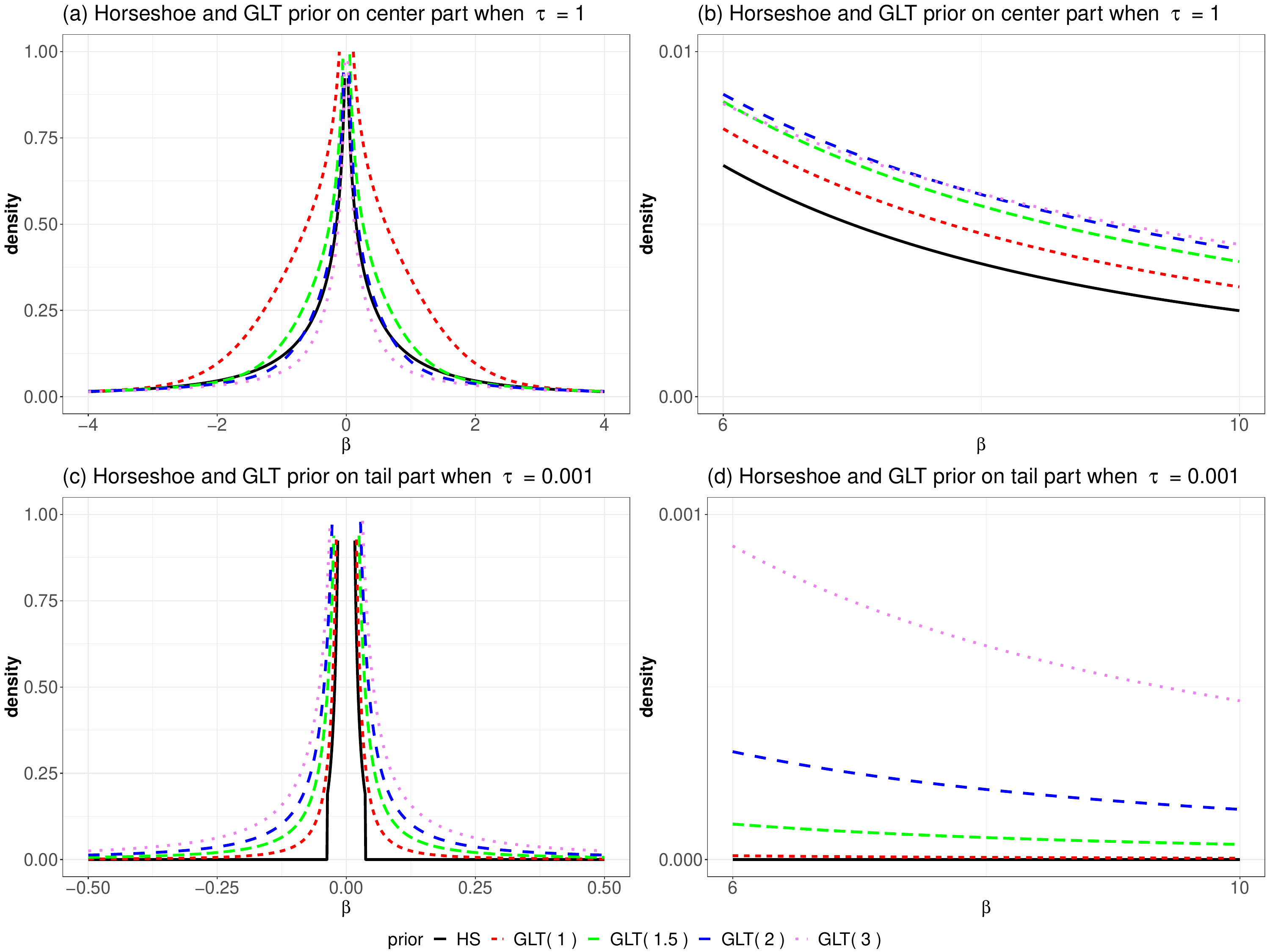}
	\caption{\baselineskip=10pt Comparison between the marginal densities of the Horseshoe and GLT prior ($\pi_{\text{HS}}(\beta|\tau)$ and $\pi(\beta|\tau, \xi)$). The global-scale parameter $\tau$ is set to $\tau=1$ (Panels (a) and (b)) and $\tau = 0.001$ (Panels (c) and (d)). The density $\pi_{\text{HS}}(\beta|\tau)$ is depicted in black, while the densities $\pi(\beta|\tau, \xi)$ are shown in red ($\xi=1$), green ($\xi=1.5$), blue ($\xi=2$), and violet ($\xi=3$), respectively.	
	}
	\label{fig:prior_comparison_beta}
\end{figure}

Figure \ref{fig:prior_comparison_beta} displays the marginal densities of the  univariate coefficient $\beta$ obtained from the Horseshoe $\pi_{\text{HS}}(\beta|\tau)$ $(\tau>0)$ (\ref{eq:marginal_beta_HS}) and the GLT prior $\pi(\beta|\tau, \xi)$ $(\tau>0, \xi>1/2)$ (\ref{eq:marginal_of_beta_glt_separate_expression}) for different values of $\tau$ and $\xi$.  The tail-heaviness of the GLT prior (varying tail rule) gets heavier as the shape $\xi$ increases, while the tail-heaviness of the Horseshoe (fixed tail rule) remains constant.

Although the GLT prior is classified as a one-group prior \citep{ghosh2017asymptotic}, Corollary \ref{corollary:separate_roles_in_beta_glt} suggests that it behaves asymptotically like a two-group prior, commonly known as the `spike-and-slab prior' \citep{george1993variable,johnson2010use,johnson2012bayesian}. This behavior is attributed to the two sequences of functions, $\{\psi_{k}^{\text{S}}(\beta)\}_{k=0}^{\infty}$ and $\{\psi_{k}^{\text{R}}(\beta)\}_{k=0}^{\infty}$.

Roughly speaking, $\{\psi_{k}^{\text{S}}(\beta)\}_{k=0}^{\infty}$ and $\{\psi_{k}^{\text{R}}(\beta)\}_{k=0}^{\infty}$ perform similar roles to the `spike' and `slab' distributions, respectively, of a two-group prior in the asymptotic cases as $|\beta|\rightarrow 0$ and $|\beta| \rightarrow \infty$.

\begin{corollary} 
	\label{corollary:separate_roles_in_beta_glt}
	Suppose $\beta|\lambda \sim N_{1}(0, \lambda^{2})$, $\lambda \sim \mathcal{GPD} (\tau, \xi)$, $\tau>0$, and $ \xi >1/2$. Let $k \in \{0\} \cup \{1,2, \cdots \}$. Then: 
	\begin{enumerate}
		\item[(a)] If $k=0$, then
		$\lim_{|\beta| \rightarrow 0}\psi_{k}^{\text{S}}(\beta) = \infty$; if $k \in \{1,2, \cdots \}$, then $\lim_{|\beta| \rightarrow 0}\pi_{k}^{\text{S}}(\beta) = 2/k < \infty$.
		\item[(b)] If $k \in \{0\} \cup  \{1,2, \cdots \}$, then 
		$\lim_{|\beta| \rightarrow \infty} \psi_{k}^{\text{S}}(\beta) = 0$ with squared exponential rate.
		\item[(c)] 
		If $k \in \{0\} \cup  \{1,2, \cdots \}$, then 
		$\lim_{|\beta| \rightarrow 0}\psi_{k}^{\text{R}}(\beta) = 
		2/(1+1/\xi + k)
		< \infty$.
		\item[(d)]
		If $k \in \{0\} \cup  \{1,2, \cdots \}$, then $\psi_{k}^{\text{R}}(\beta)$ is regularly varying with index $1+1/\xi + k$.
	\end{enumerate}
\end{corollary}
Corollary \ref{corollary:separate_roles_in_beta_glt} is proven in Subsection \ref{subsec:proof-- corollary:separate_roles_in_beta_glt} in Appendix. Interpretations of the Corollary \ref{corollary:separate_roles_in_beta_glt} are as follows. $(a)$ implies that the marginal density $\pi(\beta|\tau, \xi)$ (\ref{eq:marginal_of_beta_glt_separate_expression}) retains the infinite spike at origin for any $\tau>0, \xi > 1/2$, as seen in the Figure \ref{fig:prior_comparison_beta}, which is a common feature of the Horseshoe \citep{carvalho2010horseshoe}. Technically, this infinite spike is caused by the exponential integral function $E_{1}(x)$ ($\lim_{x\rightarrow 0^{+}} E_{1}(x) = \infty$ \citep{chiccoli1992concerning}),  allowing a very strong pulling of the $\beta$ towards zero. By $(a)$ and $(c)$ of Corollary \ref{corollary:separate_roles_in_beta_glt}, it holds $
\lim_{|\beta| \rightarrow 0}\pi_{k}^{\text{S}}(\beta) 
= 
2/k
>
\lim_{|\beta| \rightarrow 0}\psi_{k}^{\text{R}}(\beta) = 2/(1+1/\xi + k)
\,k  \in \{1,2, \cdots \},$
which implies that the contribution of $\{ \pi_{k}^{\text{S}}(\beta)\}_{k=0}^{\infty}$ is larger than that of $\{\psi_{k}^{\text{R}}(\beta)\}_{k=0}^{\infty}$ in shrinking the $\beta$ towards zero. By $(b)$, the squared exponential decay rates of the functions in $\{ \pi_{k}^{\text{S}}(\beta)\}_{k=0}^{\infty}$ as $|\beta| \rightarrow  \infty$ indicates that the contribution of $\{ \pi_{k}^{\text{S}}(\beta)\}_{k=0}^{\infty}$ in controlling the tail region of the density $\pi(\beta|\tau, \xi)$ gets negligible as $|\beta|$ goes to infinity. Finally, $(d)$ implies the density $\pi(\beta|\tau, \xi)$ possesses a systematic  mechanism to control the tail region by controlling the $\xi$ via the sequence of functions $\{\psi_{k}^{\text{R}}(\beta)\}_{k=0}^{\infty}$.

The main role of the sequence of functions $\{\psi_{k}^{\text{R}}(\beta)\}_{k=0}^{\infty}$ is to lift the tail part of the density $\pi(\beta|\tau, \xi)$ (\ref{eq:marginal_of_beta_glt_separate_expression}) by increasing the value of $\xi$. Its inclusion in the marginal prior $\pi(\beta|\tau, \xi)$ offers great flexibility in shaping the density, as demonstrated in the panels of Figure \ref{fig:prior_comparison_beta}. This flexibility proves particularly valuable in handling diverse sparsity regimes.

In contrast, the marginal density of the Horseshoe, denoted as $\pi_{\text{HS}}(\beta|\tau)$ (\ref{eq:marginal_beta_HS}), lacks a tail-controlling mechanism, as illustrated in Corollary \ref{corollary:tail_index_of_Horseshoe_prior}. This poses a significant issue, particularly when estimating a very small value for $\tau$ (e.g., $\tau = 0.001$). Panels (c) and (d) in Figure \ref{fig:prior_comparison_beta} demonstrate a discrepancy between the theoretical support $\mathbb{R}$ and the numerical support $(-\epsilon, \epsilon)$, where $\epsilon \approx 0$, of the density $\pi_{\text{HS}}(\beta|\tau = 0.001)$ (1). When $\tau$ is extremely small, such as $\tau = 10^{-10}$, the density $\pi_{\text{HS}}(\beta|\tau = 10^{-10})$ numerically converges to the Dirac delta function, potentially resulting in a collapsing behavior.

Analytic characteristics of the density for the random shrinkage coefficient $\kappa$ (\ref{eq:marginal_of_kappa_GLT}) are:
\begin{proposition}\label{prop:GLT_random_shrinkage_coefficient}
	Suppose $\lambda \sim \mathcal{GPD}(\tau, \xi)$, $\kappa = 1/(1 + \lambda^{2}) \in (0,1)$, $\tau>0$ and $\xi>1/2$. Then:
	\begin{enumerate}
		\item[(a)]
		$\lim_{\kappa \rightarrow 1^{-}} \pi(\kappa|\tau, \xi)
		=\infty
		$ and 
		$\lim_{\kappa \rightarrow 0^{+}} \pi(\kappa|\tau, \xi) = \infty$.
		\item[(b)] 
		$
		\pi(\kappa|\tau = 1, \xi = 1) 
		=
		\{\kappa^{-1/2}(1-\kappa)^{-1/2}\}/[2\cdot \{ \kappa^{1/2} + (1-\kappa)^{1/2}\}^{2}]
		.$
	\end{enumerate}
\end{proposition}
The probability mass allocated to the small regions $(1 - \epsilon, 1)$ and $(0, \epsilon)$, $\epsilon \approx 0$ under the density $\pi(\kappa|\tau, \xi)$ (\ref{eq:marginal_of_kappa_GLT}) are related to the shrinkage and the robustness, respectively \citep{carvalho2010horseshoe}. The infinite spikes of $\pi(\kappa|\tau, \xi)$ at $k=0$ and $k=1$ imply that the GLT prior has the desired shrinkage property. The density $\pi(\kappa|\tau = 1, \xi = 1)$ is not a standardly known distribution but resembles a horseshoe shape.
\begin{figure}[h]
	\centering
	\includegraphics[width=\textwidth]{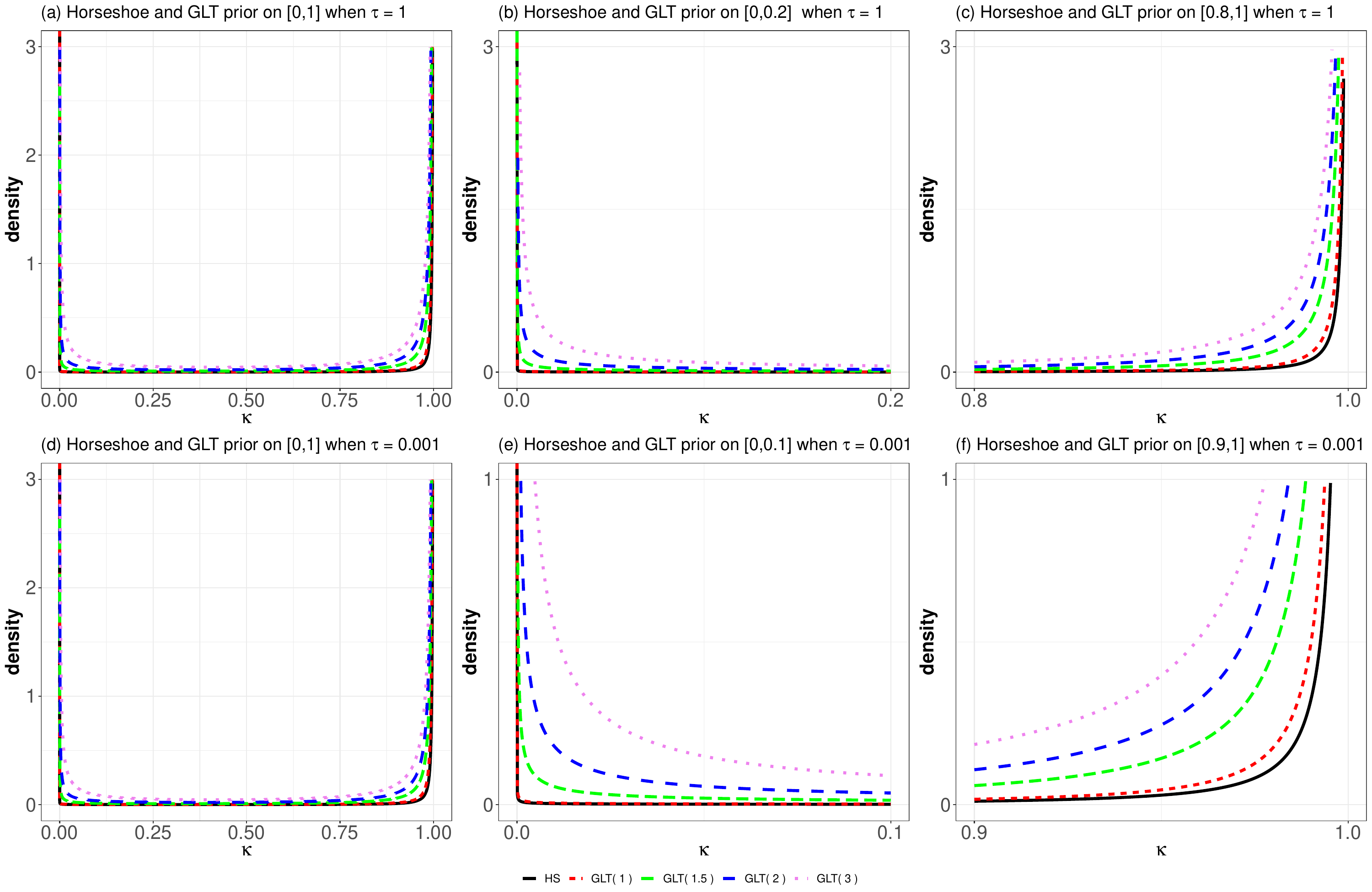}
	\caption{\baselineskip=10pt Comparison between two densities of the random shrinkage coefficient ($\pi_{\text{HS}}(\kappa|\tau)$ and $\pi(\kappa|\tau, \xi)$).  The global-scale parameter $\tau$ is set to $\tau=1$ (Panels (a), (b), and (c)) and $\tau = 0.001$ (Panels (d), (e), and (f)). The density $\pi_{\text{HS}}(\kappa|\tau)$ is depicted in black, while the densities $\pi(\kappa|\tau, \xi)$ are shown in red ($\xi=1$), green ($\xi=1.5$), blue ($\xi=2$), and violet ($\xi=3$), respectively.	
	}
	\label{fig:prior_comparison_kappa}
\end{figure}

Figure \ref{fig:prior_comparison_kappa} compares the densities of random shrinkage coefficient from the Horseshoe and GLT prior, $\pi_{\text{HS}}(\kappa|\tau)$ (refer to (\ref{eq:marginal_kappa_HS}) in Appendix) and $\pi(\kappa	|\tau, \xi)$ (\ref{eq:marginal_of_kappa_GLT}), with different values of $\tau$ and $\xi$. When $\tau=1$, the Panel (a) visualizes horseshoe shapes for both densities $\pi_{\text{HS}}(\kappa|\tau=1)$ and $\pi(\kappa|\tau=1, \xi)$. However, when $\tau=0.001$, the apparent difference is shown on Panel (e), where $\pi_{\text{HS}}(\kappa|\tau= 0.001)$ places essentially zero-mass on $(0 ,\epsilon)$, $\epsilon \approx 0$. This implies that the robustness property of the Horseshoe can be deteriorated when $\tau$ is very small. On the other hand, the GLT prior $\pi(\kappa|\tau= 0.001, \xi)$ still places a positive mass on $(0 ,\epsilon)$, $\epsilon \approx 0$, and the mass increases as the $\xi$ increases. This implies that the robustness property of the GLT prior can be maintained even when $\tau$ is very small, and is adjustable by controlling the $\xi$.

\subsection{\textbf{Propriety of the posterior of shape parameter}}\label{subsec:Tail learnability of the GLT prior}
In this subsection, we examine the propriety of the posterior distribution of the shape parameter $\xi$ under the GLT formulation (\ref{eq:gpd_beta})--(\ref{eq:gpd_xi}). This is important because the shape parameter $\xi$ is located the farthest from the response vector within the hierarchy; hence, a violation of propriety may imply that the shape parameter cannot be estimated properly. Here, we establish the proprieties of the posterior distribution of $\xi$ under two forms of the GLT prior: (a) the posterior density $\pi(\xi|y, \beta, \tau, \lambda)$ under a univariate hierarchy and (b) the full conditional posterior density $\pi(\xi|-) = \pi(\xi|\bm{\lambda}, \tau)$ used in a Gibbs sampler.


\begin{lemma}
	\label{lemma:learnability_of_xi} The following statements hold under the two formulations of the GLT prior:
	\begin{enumerate}
		\item[(a)] Assume $y|\beta \sim \mathcal{N}_{1}(\beta, 1), \beta|\lambda \sim \mathcal{N}_{1}(0, \lambda^{2}), \lambda|\tau, \xi \sim \mathcal{GPD}(\tau, \xi),$ and $\tau|\xi \sim \mathcal{IG}(1/\xi + 1, 1)$. Let $\pi(\xi)$ be any proper density of $\xi$ supported on $(1/2, \infty)$, i.e., $\int_{1/2}^{\infty} \pi(\xi) d\xi=1$. Then $\pi(\xi|y, \beta, \tau, \lambda)$ is proper on $(1/2, \infty)$.
		\item[(b)] Assume $\textbf{y} \sim \mathcal{N}_{n}(\textbf{X}\bm{\beta},\sigma^{2}\textbf{I}_{n})$ (\ref{eq:sparse_high_dim_linear_regression}) and $\bm{\beta}\sim \pi_{\text{GLT}}(\bm{\beta})$ (\ref{eq:gpd_beta}) -- (\ref{eq:gpd_xi}). Then a proportional part of the full conditional posterior for $\xi$ is represented as:
		\begin{align}
			\label{eq:general_condi_of_xi}
			\pi(\xi | - ) = 
			\pi(\xi |\bm{\lambda},\tau) 
			&\propto
			\mathcal{V}_{p}(\xi)
			\cdot
			\log\ \mathcal{N}_{1}(\xi|\mu, \rho^{2})
			\cdot
			\mathcal{I}_{(1/2, \infty)}(\xi)
			,\\
			\nonumber
			\mathcal{V}_{p}(\xi)
			&=
			\frac{\pi^{p/2} }{\Gamma(p/\xi + 1)}
			\prod_{j=1}^{p}r_{j}(\xi)
			,
		\end{align}
		where $\{r_{j}(\xi)\}_{j=1}^{p} = (\tau + \xi \lambda_{j})^{-(1/\xi + 1 )}$. The density $\pi(\xi|-)$ is proper on $(1/2, \infty)$. Here, $\mathcal{V}$ stands for volume.
	\end{enumerate}
\end{lemma}
Lemma \ref{lemma:learnability_of_xi} is proven in Subsection \ref{subsec:proof- lemma:learnability_of_xi} in Appendix. Interestingly, the likelihood part of the full conditional posterior density $\pi(\xi|-)$ (\ref{eq:general_condi_of_xi}) has a nice geometric interpretation: when $\xi=2$, the value of $\mathcal{V}_{p}(2)$ of the density corresponds to the volume of a $p$-dimensional ellipsoid with $p$ radii $\{r_{j}(2) = (\tau + 2 \lambda_{j})^{-(3/2)}\}_{j=1}^{p}$.

\subsection{\textbf{Convergence properties of the GLT posterior}}\label{subsec:Properties of the GLT posterior}
Suppose we observe $\textbf{y} = \{y_{i}\}_{i=1}^{n}$, where $y_{i} = \beta_{0i} + \epsilon_{i}$ and $\epsilon_{i} \sim \mathcal{N}_{1}(0,1)$. Here, $\bm{\beta}_0 =(\beta_{01},\cdots,\beta_{0n})^{\top}$ represents the true data generating parameter with support $S=\{k: \beta_{0k}\ne 0\}$ and the number of signals $|S|=q$. Therefore, the sparsity level (\ref{eq:sparsity_level_gene_expression_data}) is given by $s=q/n$. The assumptions of the GLT prior are as follows: $\beta_{i}|\lambda_{i} \sim \mathcal{N}_{1}(0, \lambda_{i}^{2})$ and $\lambda_{i}\sim \mathcal{GPD}(\tau, \xi)$ for $i=1,\cdots,n$, where $0< \tau<1/2$ and $1/2 < \xi$. We denote the posterior mean $\mathbb{E}[\bm{\beta}|\textbf{y}]$ as $T(\textbf{y}) = (T(y_{1}),\cdots,T(y_{n}) )^{\top} = (\mathbb{E}[\beta_1|y_1], ...,\mathbb{E}[\beta_n|y_n])^{\top} \in \mathbb{R}^{n}$, which is referred to as the `GLT estimator'. In the following, we study the mean square error (MSE) of the GLT estimator and the spread of the posterior distribution, similar to research of Horesehoe estimator done by \citep{VanDerPas2014}. Eventually, we show that the optimal thresholding value is $r_{\tau,\xi} = \sqrt{ (2/\xi) \log\left(1/\tau\right)}$. 


First, we provide an upper bound on the MSE for the GLT estimator:
\begin{theorem}[MSE for the GLT estimator]\label{thm:MSE}	
	Suppose $\textbf{y} \sim \mathcal{N}_{n}(\bm{\beta}_{0},\textbf{I})$. Then the GLT estimator $T(\textbf{y})$ satisfies
	\begin{equation}
		\label{eq:MSE_results}
		\mathbb{E}_{\bm \beta_0}\|T(\textbf{y})-\bm{\beta}_0\|_{2}^2 \lesssim \frac{q}{\xi} \log \frac{1}{\tau} +(n-q) \tau^{\frac{1}{\xi}} \sqrt{\log \frac{1}{\tau}},
	\end{equation}
	for $\tau\rightarrow 0$, as $n,q \rightarrow \infty$, and $q=o(n)$.
\end{theorem}
Note that the choice of $\tau = (q/n)^{\alpha}$ and $1/2 < \xi < \alpha$, for any constant $\alpha>1/2$, yields an upper bound (\ref{eq:MSE_results}) that is of the same order as the minimax risk, $q \log (n/q)$ \citep{donoho1992maximum}. 



Next, we describe the spread of the GLT posterior distribution: 
\begin{theorem}[Spread of the GLT posterior]\label{thm:variance}
	Suppose $\textbf{y} \sim \mathcal{N}_{n}(\bm{\beta}_{0},\textbf{I})$. Then, the following relations hold for the variance of the posterior distribution corresponding to the GLT prior.
	The total posterior variance satisfies
	\begin{equation}
		\label{eq:thm_2_results}
		\mathbb{E}_{\bm \beta_0} \sum_{i=1}^{n} \mbox{Var}[\beta_i \mid y_i] \lesssim \frac{q}{\xi} \log \frac{1}{\tau} +(n-q) \tau^{\frac{1}{\xi}} \bigg(\log \frac{1}{\tau}\bigg)^{\max\{1/2,1-1/(2\xi)\}},
	\end{equation}
	for $\tau\rightarrow 0$, as $n,q \rightarrow \infty$, and $q=o(n)$.
\end{theorem}

The following theorem provides an upper bound on the rate of contraction of the full posterior distribution around the true underlying mean vector:
\begin{theorem}[Posterior contraction]\label{thm:posterior}
	Suppose $\textbf{y} \sim \mathcal{N}_{n}(\bm{\beta}_{0},\textbf{I})$. Let $\tau = (q/n)^\alpha$ and $ 1/2< \xi <\alpha $, $\alpha>1/2$. Then, we have
	$$\mathbb{E}_{\bm \beta_0}\Pi\left(\|\bm{\beta}-\bm{\beta}_0\|_2^2 \geq M_n q\log (n/q) \mid \textbf{y} \right) \rightarrow 0,$$
	for every $M_n \rightarrow \infty$, as $n,q \rightarrow \infty$, and $q=o(n)$.
\end{theorem}

The preceding theorem suggests that several choices of $\tau$ and $\xi$, such as $\tau = (q/n)^\alpha$ and $1/2 < \xi < \alpha$ for any $\alpha > 1/2$, result in an upper bound of order $q \log(n/q)$ on both the worst-case $l_2$ risk and the rate of posterior contraction. The proofs for these theorems are provided in Section \ref{sec:Convergence properties of the GLT posterior} of Appendix.

We conduct simulation experiments to investigate the alignment between (a) the theoretically optimal values of $\xi$ based on the aforementioned theorems and (b) the posterior means $\widehat{\xi} = \mathbb{E}[\xi|\textbf{y}]$ obtained through posterior computation, under various sparse regimes. The theoretical values for $\xi$ were derived to achieve the near-minimax rate of $q \log n$ \citep{van2017adaptive}, and its analytic expression is given by $\xi = \log n / (\log n - \log q) = \log n / \log s$ (possibly up to a multiplicative constant; refer to the Appendix for the derivation).

\begin{figure}[h]
	\centering
	\includegraphics[width=\textwidth]{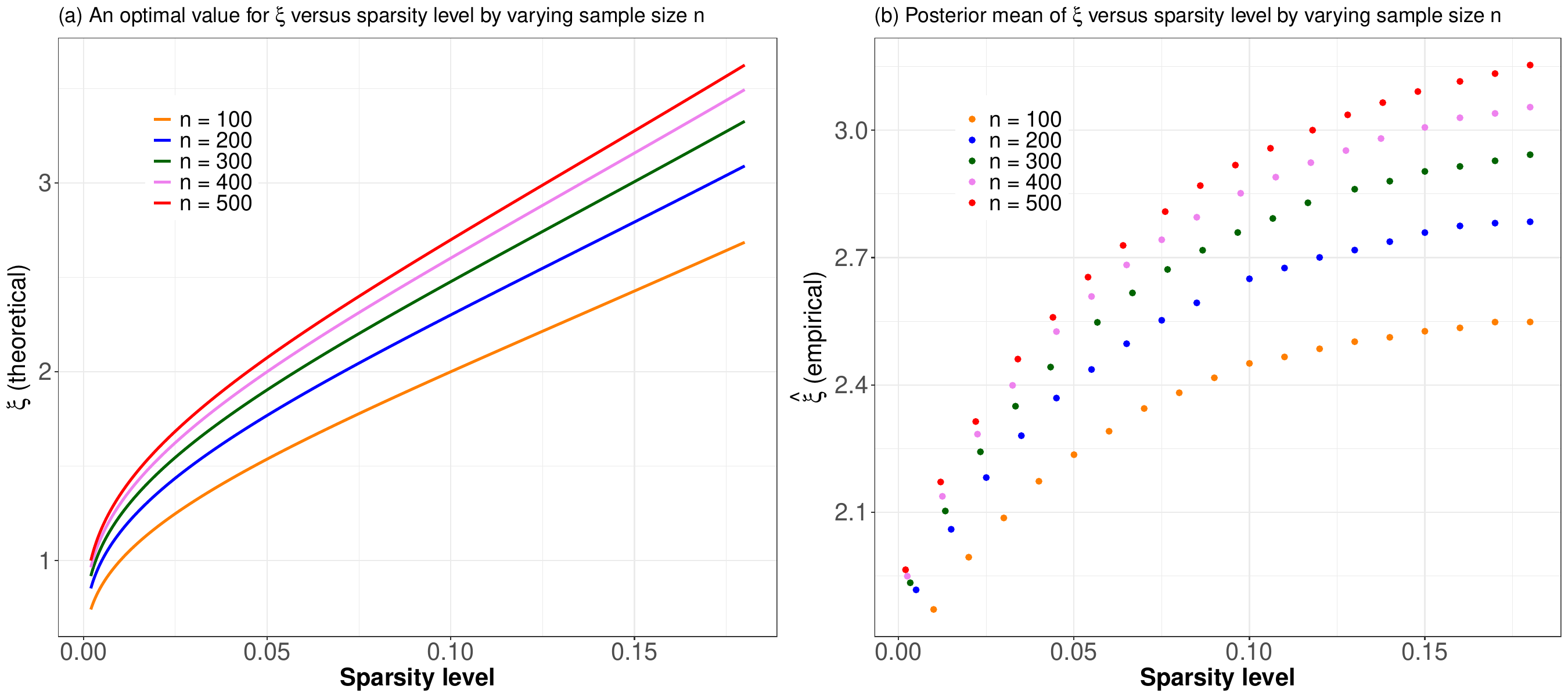}
	\caption{\baselineskip=10pt Comparisons between the optimal value $\xi$ (up to a multiplicative constant) (Panel (a)) and the posterior mean $\widehat{\xi} = \mathbb{E}[\xi|\textbf{y}]$ (Panel (b)) under sparse normal mean models. $x$-axis of the panels represents the sparsity level $s=q/n$ ranging from around 0.001 to 0.2. 
	}
	\label{fig:Theoretical_and_Empirical_Value_of_xi}
\end{figure}

Figure \ref{fig:Theoretical_and_Empirical_Value_of_xi} presents the simulation results for sample sizes $n=100, 200, 300, 400$, and $500$, while varying the sparsity level $s$ (\ref{eq:sparsity_level_gene_expression_data}) ranging from approximately 0.001 to 0.2. Panels (a) and (b) display the trends of the theoretical and empirical values of $\xi$, respectively, across the range of ultra-sparse to moderate-sparse regimes. Each point in Panel (b) represents the median of the posterior means for $\xi$ based on 50 replicated datasets. Both panels reveal a monotonically increasing trend in $\xi$ with respect to the sparsity level. These results suggest that as the sparsity level $s$ increases up to a certain constant (say, around 0.2 in this simulation), a larger value of $\xi$ leads to an optimal shrinkage effect on $\bm{\beta}$.

Next, we present motivating examples in the next two sections to appreciate the necessity of the varying tail rule compared to the fixed tail rule by contrasting the performance and posterior behaviors of the GLT prior and the Horseshoe. More simulation experiments with other shrinkage priors are contained in the Supplementary Material.
\section{Motivating example 1-- Prostate cancer data}\label{sec:Example-- prostate cancer data}
\subsection{\textbf{Prostate cancer data}}\label{subsec:Prostate cancer data}
The prostate cancer data can be downloaded from the \texttt{R} package \texttt{sda}, as detailed on page 272 of \citep{efron2016computer}. This dataset consists of a matrix $\textbf{X} \in \mathbb{R}^{102 \times 6033}$, summarizing gene expression levels measured on microarrays from two classes. The first 50 rows of $\textbf{X}$, denoted as $\textbf{X}[1:50,\cdot] \in \mathbb{R}^{50 \times 6033}$, correspond to healthy controls, while the remaining 52 rows, denoted as $\textbf{X}[51:102,\cdot] \in \mathbb{R}^{52 \times 6033}$, correspond to cancer patients. Each column vector $\textbf{X}[\cdot,j] \in \mathbb{R}^{102}$ represents the gene expression levels of the $j$-th gene, with $j$ ranging from 1 to 6,033.

The main objective of this study is to identify a subset of $q$ interesting genes out of the total 6,033 genes whose expression levels exhibit differences between the two groups \citep{efron2012large}. These selected genes are then further investigated to explore any causal links to the development of prostate cancer. The sparsity level, denoted as $s = q/p$ (\ref{eq:sparsity_level_gene_expression_data}), where $p = 6{,}033$, is unknown because it is real data.

\begin{figure}[h]
	\centering
	\includegraphics[width=\textwidth]{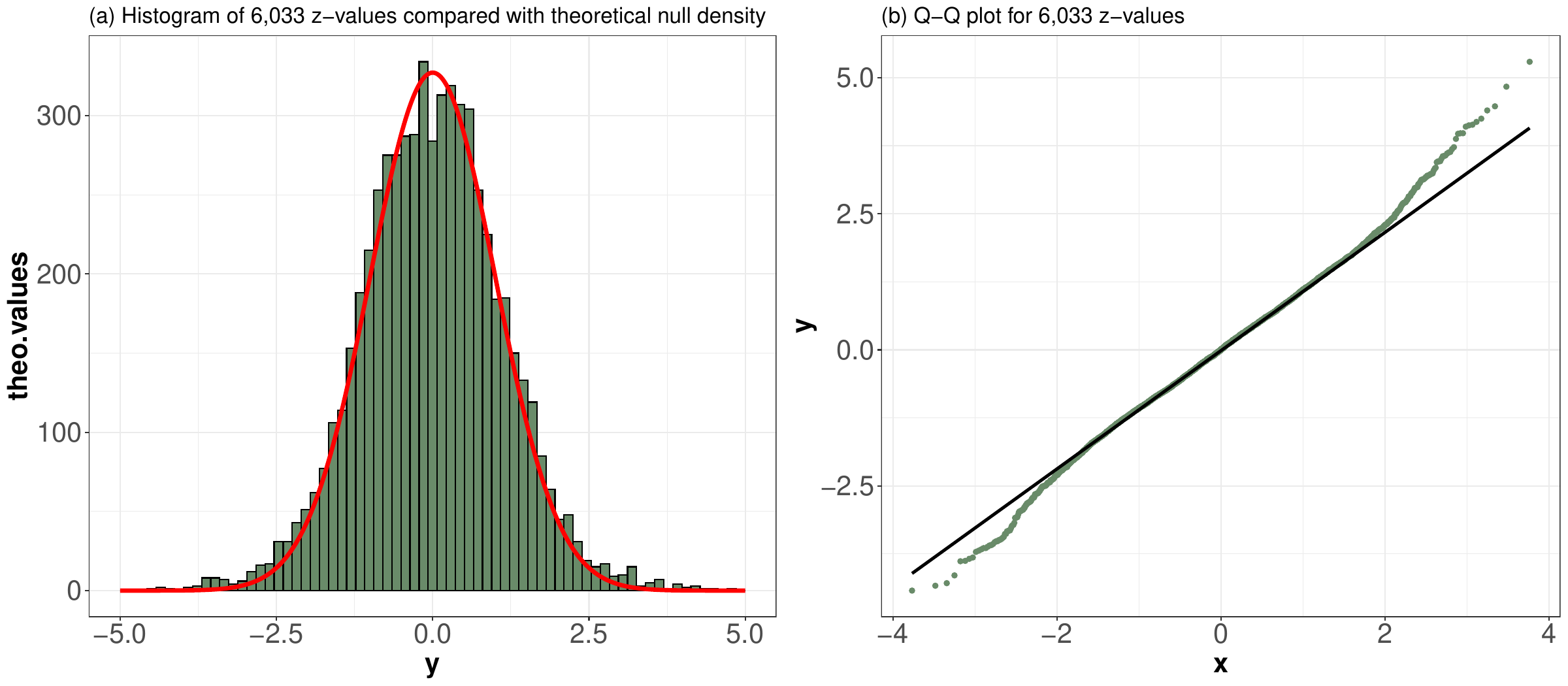}
	\caption{\baselineskip=10pt
		Histogram of $z$-values $\{y_{j} \}_{j=1}^{6033}$ obtained from prostate cancer data (Panel (a)) and the Q-Q plot (Panel (b))}
	\label{fig:Prostate_cancer_data_Efron}
\end{figure}

We adopt the data transformation method proposed by \citep{efron2010future}. The primary aim of this transformation is to indirectly address the challenge of multiple hypothesis testing \citep{benjamini1995controlling,benjamini2005false} by estimating coefficients within a sparse normal mean model \citep{bhattacharya2015dirichlet}:
\begin{align}
	\label{eq:sparse_normal_mean_model}
	y_{j}  = \beta_{j} + \sigma\epsilon_{j}, \quad \epsilon \sim \mathcal{N}_{1}(0,1), \quad j = 1,\cdots,p, \quad (p=6,033 \text{ genes}),
\end{align}
where $\sigma$ is unknown.

The responses $\{y_{j}\}_{j=1}^{p=6,033}$ are acquired as follows. First, for each $j=1,$ $\cdots, 6,033$, obtain $t$-test statistics $t_{j}$ through a two-sample $t$-test statistic with $100$ degrees of freedom based on the $j$-th vector $\textbf{X}[\cdot,j] \in \mathbb{R}^{102}$. Second, convert the acquired $t$-test statistics to $z$-test statistics using quantile transformation $y_{j} = \Phi^{-1}(F_{100 \text{ d.f.}}(t_{j}))$, where $\Phi(\cdot)$ and $F_{100 \text{ d.f.}}(\cdot)$ are distribution functions of $\mathcal{N}_{1}(0,1)$ and $t_{100}$, respectively;  refer to Section $2.1$ of \citep{efron2012large}. The histogram of $\{y_{j}\}_{j=1}^{p=6,033}$ along with the standard normal density and its Q-Q plot are displayed in Panel (a) and (b) in the Figure \ref{fig:Prostate_cancer_data_Efron}, respectively.

Denote $H_{0 j}$ as the null hypothesis that posits no difference in the gene expression levels for the $j$-th gene between the healthy controls and cancer patients. If the global null hypothesis $\cap_{j=1}^{6033} H_{0 j}$ is true, the histogram of $\{y_{j}\}_{j=1}^{p=6,033}$ should mimic a standard normal density closely. The presence of outliers is evident from the panels. Those outliers may correspond to interesting genes (cancerous genes) that reject the null hypotheses \citep{efron2010future}. 
\subsection{\textbf{Prostate cancer data analysis via the Horseshoe}}\label{subsec:Prostate cancer data analysis via the Horseshoe}
In this example, we compare the performance of the Horseshoe prior (\ref{eq:HS_beta}) -- (\ref{eq:HS_lambda_tau}) and the GLT prior (\ref{eq:gpd_beta}) -- (\ref{eq:gpd_xi}) as sparse-inducing priors for the coefficients in the sparse normal mean model (\ref{eq:sparse_normal_mean_model}), as the number of genes considered increases. For this purpose, we construct seven prostate datasets denoted as $\mathcal{P}_{1} =\{y_{j}\}_{j=1}^{p=50}$, $\mathcal{P}_{2} =\{y_{j}\}_{j=1}^{p=100}$, $\mathcal{P}_{3} =\{y_{j}\}_{j=1}^{p=200}$, $\mathcal{P}_{4} =\{y_{j}\}_{j=1}^{p=500}$,  $\mathcal{P}_{5} =\{y_{j}\}_{j=1}^{p=1000}$, $\mathcal{P}_{6} =\{y_{j}\}_{j=1}^{p=3000}$, and $\mathcal{P}_{7} =\{y_{j}\}_{j=1}^{p=6033}$. Thus, we have the subset inclusions $\mathcal{P}_{1} \subset \mathcal{P}_{2} \subset \mathcal{P}_{3} \subset \mathcal{P}_{4} \subset \mathcal{P}_{5} \subset \mathcal{P}_{6} \subset \mathcal{P}_{7}$, where $\mathcal{P}_{7}$ represents the full dataset of 6,033 genes. Here, $\mathcal{P}$ stands for prostate.

\begin{figure}[h]
	\centering
	\includegraphics[width=\textwidth]{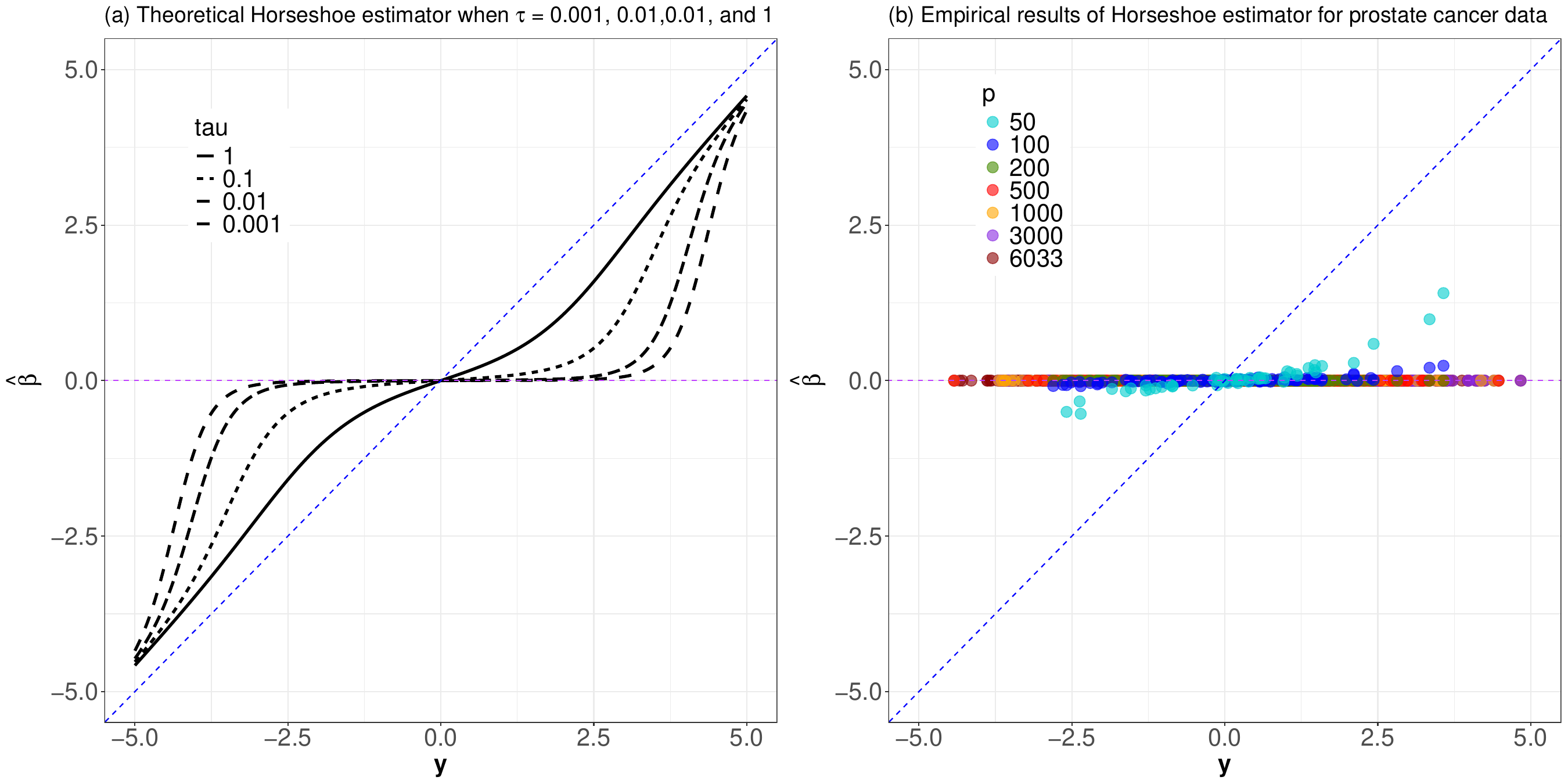}
	\caption{\baselineskip=10pt
		An idealistic reversed-$S$-shape curve (Panel (a)) is formed by pairs $\{(y_{j},\widehat{\beta}_{j})\}_{j=1}^{p}$ when the Horseshoe estimator achieves the robustness property. The posterior inference results (Panel (b)) are obtained using the Horseshoe prior on the seven datasets $\mathcal{P}_{l}$ ($l=1,\cdots,7$). The dotted line represents $y=x$. The posterior means of $\tau$ for the four datasets are $0.158$ ($\mathcal{P}_{1}$) and $0.018$ ($\mathcal{P}_{2}$), and numerically zero for $\mathcal{P}_{l}$ ($l=3,\cdots,7$).
	}
	\label{fig:idealistic_S_curve_and_HS_prostate}
\end{figure}

We present the posterior inference results by plotting ordered pairs $\{(y_{j},$ $\widehat{\beta}_{j})\}_{j=1} ^{p}$, where $\widehat{\beta}_{j}$ represents the posterior mean of $\beta_{j}$ for the $j$-th gene. If the (tail) robustness property holds, we expect to observe a reversed-$S$-shaped curve formed by the pairs $\{(y_{j},\widehat{\beta}_{j})\}_{j=1}^{p}$ in each dataset \citep{VanDerPas2014,efron2010future}. The robustness property is discussed in more detail in Section 2 of \citep{carvalho2010horseshoe}. The reversed-$S$-shaped curve, as seen in Panel (a) of Figure \ref{fig:idealistic_S_curve_and_HS_prostate}, is formed due to (i) $\widehat{\beta}_{j} \approx 0$, corresponding to a noise coefficient, (ii) $\widehat{\beta}_{j} \approx y_{j}$, corresponding to a signal coefficient, and (iii) the continuous nature of a one-group Gaussian mixture prior.

Throughout the paper, we implement the Horseshoe (\ref{eq:HS_beta}) -- (\ref{eq:HS_lambda_tau}) via the \texttt{R} function \texttt{Horseshoe} within the \texttt{R} package \texttt{Horseshoe}. More specifically, given a response vector $\textbf{y}= (y_{1},\cdots,y_{p})^{\top} \in \mathbb{R}^{p}$ formulated from one of the seven prostate cancer dataset, we use \texttt{Horseshoe(y=y, X=X,} \texttt{method.tau="halfCauchy",} \\
\texttt{method.sigma="Jeffreys",} \texttt{burn=10000, nmc=10000, thin=100)} where \texttt{y} = $\textbf{y}$ and \texttt{X} = $\textbf{I}_{p}$, to produce $100$ thinned realizations from the posterior distribution $\pi(\bm{\beta}|\textbf{y})$ via MCMC \citep{andrieu2003introduction,robert2013monte}.

The ordered pairs $\{(y_{j},\widehat{\beta}_{j})\}_{j=1}^{p}$ based on the seven datasets $\mathcal{P}_{l}$ ($l=1,\cdots,7$) are overlaid on Panel (b) of Figure \ref{fig:idealistic_S_curve_and_HS_prostate}. The results suggest that the robustness property is observed only when $p=50$, and the property disappears as $p$ increases. For $p \geq 200$, the posterior mean becomes numerically zero, i.e., $\widehat{\bm{\beta}} = (\widehat{\beta}_{1},\cdots,\widehat{\beta}_{p})^{\top} \approx \textbf{0} = (0,\cdots,0)^{\top} \in \mathbb{R}^{p}$. Therefore, the plotted dots essentially form the set $\{(y_{j}, 0)\}_{j=1}^{p}$, exhibiting the collapsing behavior.

In general, it should be noted that the ordinary estimates for the $\beta_{j}$'s in the normal mean model (\ref{eq:sparse_normal_mean_model}) are $y_{j}$, while Bayesian estimates are typically biased, as they are influenced by prior information on the $\beta_{j}$ \citep{lindley1972bayes}. However, what we observed in this study is an unusually strong shrinkage towards the prior mean. 
\subsection{\textbf{Prostate cancer data analysis via the GLT prior}}\label{subsec:Prostate cancer data analysis via the GLT prior}
\begin{figure}[h]
	\centering
	\includegraphics[width=\textwidth]{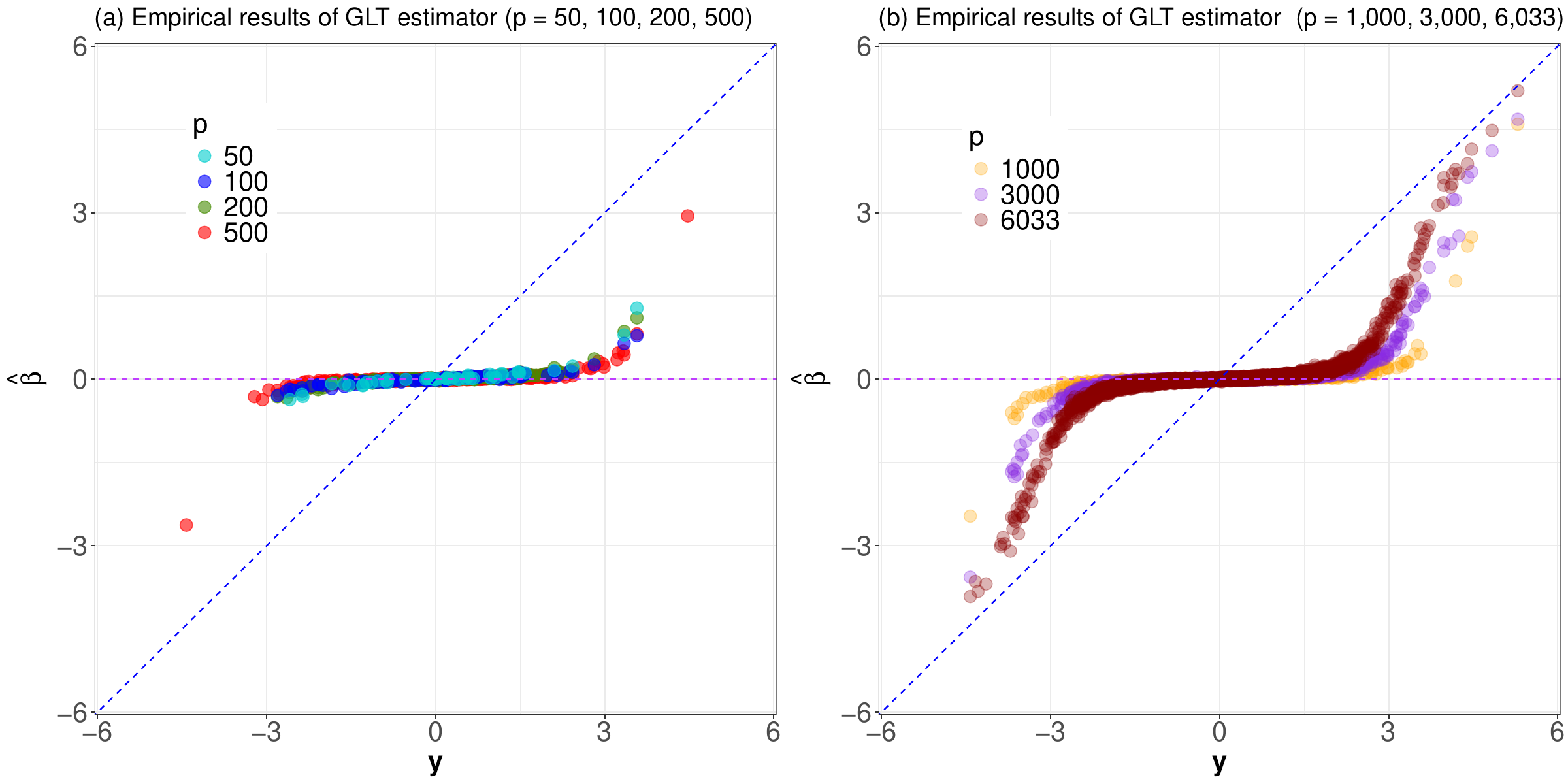}
	\caption{\baselineskip=10pt
		Posterior inference results obtained by the GLT prior applied seven prostate cancer datasets $\mathcal{P}_{l}$, $l=1,\cdots,7$. Posterior means of $(\tau,\xi)$ corresponding to the seven datasets are 
		$(0.0303, 1.620)$ ($\mathcal{P}_{1}$), 
		$(0.0154, 1.662)$ ($\mathcal{P}_{2}$), 
		$(0.0090, 1.789)$ ($\mathcal{P}_{3}$),
		$(0.0037, 1.905)$ ($\mathcal{P}_{4}$),
		$(0.0019, 1.991)$ ($\mathcal{P}_{5}$),
		$(0.0013, 2.760)$ ($\mathcal{P}_{6}$),
		and 
		$(0.0013, 3.636)$ ($\mathcal{P}_{7}$), respectively.
	}
	\label{fig:GLT_prostate_data}
\end{figure}

The posterior inference results obtained using the GLT prior (\ref{eq:gpd_beta}) -- (\ref{eq:gpd_xi}) are presented in Figure \ref{fig:GLT_prostate_data}. Panel (a) and (b) display the pairs $\{(y_{j},\widehat{\beta}_{j})\}_{j=1} ^{p}$ based on the four datasets $\mathcal{P}_{l}$ ($l=1,2,3,4$) and the three datasets $\mathcal{P}_{l}$ ($l=5,6,7$), respectively. It is observed that the desirable reversed-$S$-shape curves are manifested in all datasets, indicating that the robustness property holds regardless of the number of genes used. It should be noted that neither the Horseshoe nor the GLT prior requires any tuning procedure, ensuring a fair comparison of their performances.

Posterior means of the shape parameter $\xi$ for the seven datasets are as follows: $ 1.620$ $(\mathcal{P}_{1})$,
$1.662$ $(\mathcal{P}_{2})$,
$1.789$ $(\mathcal{P}_{3})$,
$1.905$ $(\mathcal{P}_{4})$,
$1.991$ $(\mathcal{P}_{5})$,
$2.760$ $(\mathcal{P}_{6})$, and 
$3.636$ $(\mathcal{P}_{7})$. These results indicate that the posterior tail-thickness gets heavier as the number of genes considered $p$ increases to accommodate a growing number of interesting genes, denoted as $q$, which is unknown. 

In the application using full data $(\mathcal{P}_{7})$ with 6033 genes, we expect there to be two clusters of absolute values of $\beta_{j}$'s, with one cluster concentrated closely near zero corresponding to genes that are effectively not differentially expressed (i.e., noise coefficients) and another away from zero corresponding to interesting genes for further study (i.e., signal coefficients). As a simple approach for post-processing variable selection, we use 2-mean variable selection \citep{li2017variable}, where the clusters of $|\beta_{j}|$'s at each MCMC iteration are identified using $k$-means with $k=2$. For each iteration, the number of nonzero signals is then estimated by the smaller cluster size out of the two clusters. A final estimate $(M)$ of the number of nonzero signals is obtained by taking the mode over all the MCMC iterations. The $M$ largest (in absolute magnitude) entries of the posterior median are identified as the nonzero signals. Using the above selection scheme, the GLT prior declared 184 genes as nonnull, with the most interesting gene (in absolute magnitude) indexed as 610. (This gene 610 was also identified as the most interesting gene by Efron; see Table 11.2 \citep{efron2012large}.) The Horseshoe estimator has been collapsed when applied to the full data (as observed in Panel (b) in Figure \ref{fig:idealistic_S_curve_and_HS_prostate}); hence, no gene was selected, which was similarly identified by \citep{bhattacharya2016fast} where the authors found that the Horseshoe is overly conservative, even more conservative than the classical Bonferroni correction.	
\section{Motivating example 2-- Simulation study with varied sparsity level}\label{sec:Example-- Simulation study with varied sparsity level}
\subsection{\textbf{Artificial high-dimensional data generator}}\label{subsec:Artificial high-dimensional data generator}
In this subsection, we illustrate the generating process of high-dimensional data $(\textbf{y}, \textbf{X}) \in \mathbb{R}^{n} \times \mathbb{R}^{n \times p}$ from the true data-generating distribution $p(\textbf{y}, \textbf{X})$, where the data are from the high-dimensional regression (\ref{eq:sparse_high_dim_linear_regression}) given a simulation environment $(n, p, q, \varrho, \text{SNR})$:
\begin{align}
	\label{eq:construct_high_dimensional_data}
	(\textbf{y},\textbf{X})
	\sim 
	p(\textbf{y},\textbf{X}) &=
	\mathcal{N}_{n}(\textbf{y}|\textbf{X}\bm{\beta}_{0} , \sigma_{0}^{2} \textbf{I}_{n} )
	\cdot
	\prod_{i=1}^{n}
	\mathcal{N}_{p}(\textbf{x}_{i}^{\top} | \textbf{0}, \bm{\Upsilon }(\varrho) ),\\
	\nonumber
	\bm{\Upsilon }(\varrho) &= \varrho \textbf{J}_{p} + (1 - \varrho) \textbf{I}_{p}.
\end{align}
Here, the vector $\bm{\beta}_{0} = (\beta_{0,1}, \cdots, \beta_{0,q}, \beta_{0,q+1}, \cdots, \beta_{0,p})^{\top} \in \mathbb{R}^{p}$ represents the true coefficient vector. It consists of $q$ unit signals, where $\beta_{0,1} = \cdots = \beta_{0,q} = 1$ for the first $q$ coefficients, and $p-q$ noise coefficients, where $\beta_{0,q+1} = \cdots = \beta_{0,p} = 0$, for the remaining coefficients. Sparsity level (\ref{eq:sparsity_level_gene_expression_data}) is then $s = q/p$. The matrices $\textbf{I}$ and $\textbf{J}$ in (\ref{eq:construct_high_dimensional_data}) indicate an identity matrix and a matrix whose elements are ones, respectively. The signal-to-noise ratio (SNR) is defined by $\text{SNR} = \text{var}(\textbf{X}\bm{\beta}_{0})/\text{var}(\sigma_{0}\bm{\epsilon})$. The value $\varrho$ is a number associated with column-wise correlations in the design matrix $\textbf{X}$. 

After specifying a simulation environment $(n, p, q, \varrho, \text{SNR})$, we take the following three steps to generate high-dimensional data $(\textbf{y},\textbf{X})$:
\begin{itemize}
	\item[] (i) Generate a matrix $\textbf{X}\in \mathbb{R}^{n \times p}$, where each row vector $\textbf{x}_{i} \in \mathbb{R}^{p}$ is independently sampled from $\mathcal{N}_{p}(\textbf{0}, \bm{\Upsilon}(\varrho) )$. Next, center the matrix $\textbf{X}$ column-wise so that each column vector $\textbf{X}[\cdot,j] \in \mathbb{R}^n$ ($j=1,\cdots,p$) has zero mean. After that, normalize each column vector to have a unit Euclidean $l_{2}$-norm.
	\item[] (ii) Generate an $n$-dimensional Gaussian error $\bm{\epsilon }\sim \mathcal{N}_{n}(\textbf{0}, \textbf{I}_{n})$.
	\item[] (iii) Add the mean part $\textbf{X}\bm{\beta}_{0}$ and the error part $\sigma_{0} \bm{\epsilon}$ to create the response vector $\textbf{y }= \textbf{X}\bm{\beta}_{0} + \sigma_{0} \bm{\epsilon}$, where $\sigma_{0}^{2} = \text{var}(\textbf{X}\bm{\beta}_{0})/(\text{SNR} \cdot \text{var}(\bm{\epsilon}) )$, with $\text{var}(\textbf{z}) = \sum_{i=1}^{n}(z_{i} - \bar{z})^{2}/(n-1)$ for $\textbf{z} \in \mathbb{R}^{n}$.
\end{itemize}


We use these steps in the next subsections to explore the posterior behavior of the Horseshoe and GLT prior and conduct simulation experiments in Section \ref{sec:Simulations}.
\subsection{\textbf{The Horseshoe under varied sparsity level}}\label{subsec:The Horseshoe with varied sparsity level}
To investigate the behavior of the Horseshoe (\ref{eq:HS_beta})--(\ref{eq:HS_lambda_tau}) as the sparsity level $s = q/p$ (\ref{eq:sparsity_level_gene_expression_data}) increases, we generate four high-dimensional datasets denoted as $\mathcal{A}_{l} = (\textbf{y}, \textbf{X}) \in \mathbb{R}^{n} \times \mathbb{R}^{n \times p}$, $l = 1, 2, 3, 4$. These datasets correspond to four simulation environments with $n = 100$, $p = 500$, $\varrho = 0$, and $\text{SNR} = 5$, where only the number of true signals varies, i.e., $q = 2, 5, 8, 13$. The corresponding sparsity levels of the datasets are $0.004$ $(\mathcal{A}_{1})$, $0.01$ $(\mathcal{A}_{2})$, $0.016$ $(\mathcal{A}_{3})$, and $0.026$ $(\mathcal{A}_{4})$, respectively.

\begin{figure}[h]
	\centering
	\includegraphics[width=\textwidth]{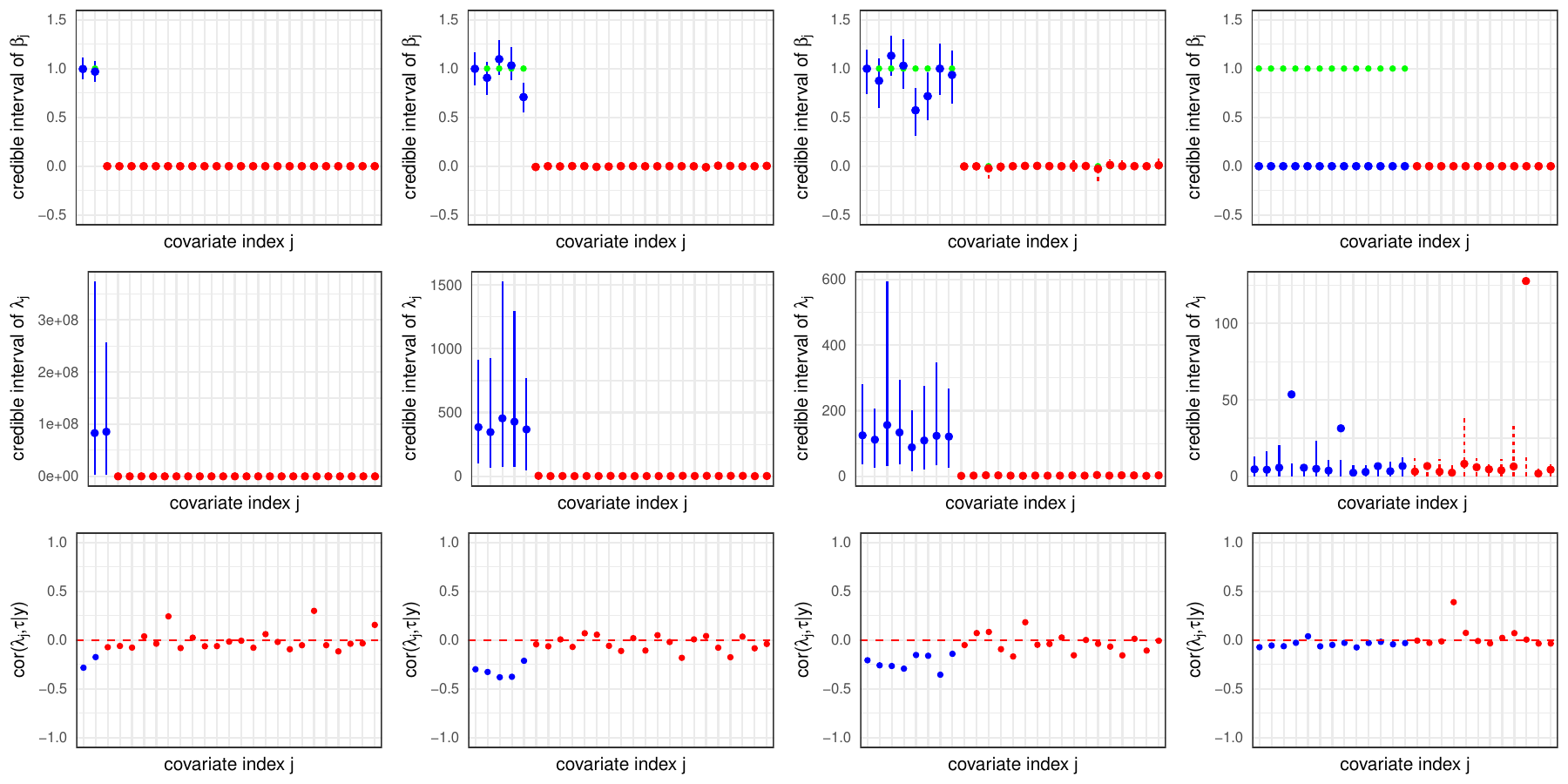}
	\caption{\baselineskip=10pt
		The results of posterior inference using the Horseshoe under varying sparsity levels: $\mathcal{A}_{1}$ ($s =0.004$, first column), $\mathcal{A}_{2}$ ($s =0.01$, second column), $\mathcal{A}_{3}$ ($s=0.016$, third column), and $\mathcal{A}_{4}$ ($s=  0.026$, fourth column). In panels, the posterior inference corresponding to signals and noises is colored in blue and red, respectively, while the ground truth $\beta_{0}$ is colored in green. The posterior means of $\tau$ for the four datasets are $1.41 \cdot 10^{-6}$ ($\mathcal{A}_{1}$), $0.05$ ($\mathcal{A}_{2}$), $0.13$ ($\mathcal{A}_{3}$), and $6.53 \cdot 10^{-15}$ ($\mathcal{A}_{4}$), respectively.
	}
	\label{fig:failure of Horseshoe}
\end{figure}

The results of posterior inference are displayed in Figure \ref{fig:failure of Horseshoe}. The panels are arranged in a way that the sparsity level increases from left to right. The first, second, and third rows of panels in Figure \ref{fig:failure of Horseshoe} show the $95\%$ credible intervals for $\{\beta_{j}\}_{j=1}^{p}$, $\{\lambda_{j}\}_{j=1}^{p}$, and the posterior correlations $\{\text{cor}(\lambda_{j}, \tau | \textbf{y}) \}_{j=1}^{p}$, respectively. For ease of visualization, only the results corresponding to the first $25$ coefficients of $\beta$ are displayed. The coefficients corresponding to signals and noises are colored blue and red, respectively, while the true coefficient vector $\beta_{0}$ is colored green.

The results show that the Horseshoe works reasonably well on the ultra-sparse datasets $\mathcal{A}_{1}$, $\mathcal{A}_{2}$, and $\mathcal{A}_{3}$, but exhibits collapsing behavior on the moderately sparse dataset $\mathcal{A}_{4}$. The posterior means of $\tau$ corresponding to the datasets are $1.41 \cdot 10^{-6}$ $(\mathcal{A}_{1})$, $0.05$ $(\mathcal{A}_{2})$, $0.13$ $(\mathcal{A}_{3})$, and $6.53 \cdot 10^{-15}$ $(\mathcal{A}_{4})$, respectively. Hence, the posterior mean of $\tau$ gradually increases as the sparsity level increases, and after reaching a certain threshold, it suddenly drops to a very small number, virtually zero.

The key to understanding how the Horseshoe detects signals lies in the relationship between the local scale parameters $\{\lambda_{j}\}_{j=1}^{p}$ and the global scale parameter $\tau$. It is important to note that $\tau$ is associated with the sparsity level \citep{polson2010shrink} and is expected to be large when there are relatively many signals present. The panels on the third row of Figure \ref{fig:failure of Horseshoe} indicate the weak negative posterior correlation between $\lambda_{j}$ and $\tau$, denoted as $\text{cor}(\lambda_{j}, \tau | \mathbf{y})$, for each $j = 1, \cdots, p$. As observed in the panels on the first and third rows of Figure \ref{fig:failure of Horseshoe}, the selected signals among the $p$ coefficients $\{\beta_{j}\}_{j=1}^{p}$, denoted as $\{\beta_{j}\}_{j \in \mathcal{Q}}$ with $\mathcal{Q} \subset \mathcal{P} = \{1, \cdots, p\}$, are those for which the corresponding posterior correlations $\{\text{cor}(\lambda_{j}, \tau | \mathbf{y})\}_{j \in \mathcal{Q}}$ exhibit even stronger negative values compared to the negative correlations $\{\text{cor}(\lambda_{j}, \tau | \mathbf{y})\}_{j \in \mathcal{P} - \mathcal{Q}}$ of the remaining coefficients $\{\beta_{j}\}_{j \in \mathcal{P} - \mathcal{Q}}$. To see this, note that, for the datasets $\mathcal{A}_{1}$ (the first column), $\mathcal{A}_{2}$ (the second column), and $\mathcal{A}_{3}$ (the third column), the signal coefficients colored as blue dots in the panels are well detected, and the corresponding correlation values are visually distinguishable from the correlations corresponding to the noise coefficients colored as red dots. On the other hand, if there are no distinguishable differences among the correlation values $\{\text{cor}(\lambda_{j}, \tau | \mathbf{y})\}_{j \in \mathcal{P}}$ as observed in the dataset $\mathcal{A}_{4}$, the signal detection mechanism of the Horseshoe may be lost, potentially leading to the collapsing behavior as observed in the panel of the fourth column.

As briefly discussed in the Introduction, the collapsing behavior of the Horseshoe estimator has been highlighted by several authors, but it has not received significant attention in the literature. One recent discussion involves the potential collapse of the empirical Bayes estimator of the Horseshoe, specifically the marginal maximum-likelihood estimator for the global-scale parameter $\tau$, when the sparsity level $s$ (\ref{eq:sparsity_level_gene_expression_data}) is very small \citep{van2017uncertainty, yoo2017contributed}. It is important to note that our research focuses on the collapsing behavior of the fully Bayesian Horseshoe estimator, where $\tau$ is distributed according to $\mathcal{C}^{+}(0,1)$, under the moderately sparse regime. Therefore, we address different problems.

\subsection{\textbf{The GLT prior under varied sparsity level}}\label{subsec:The GLT prior under varied sparsity level}

\begin{figure}[h]
	\centering
	\includegraphics[width=\textwidth]{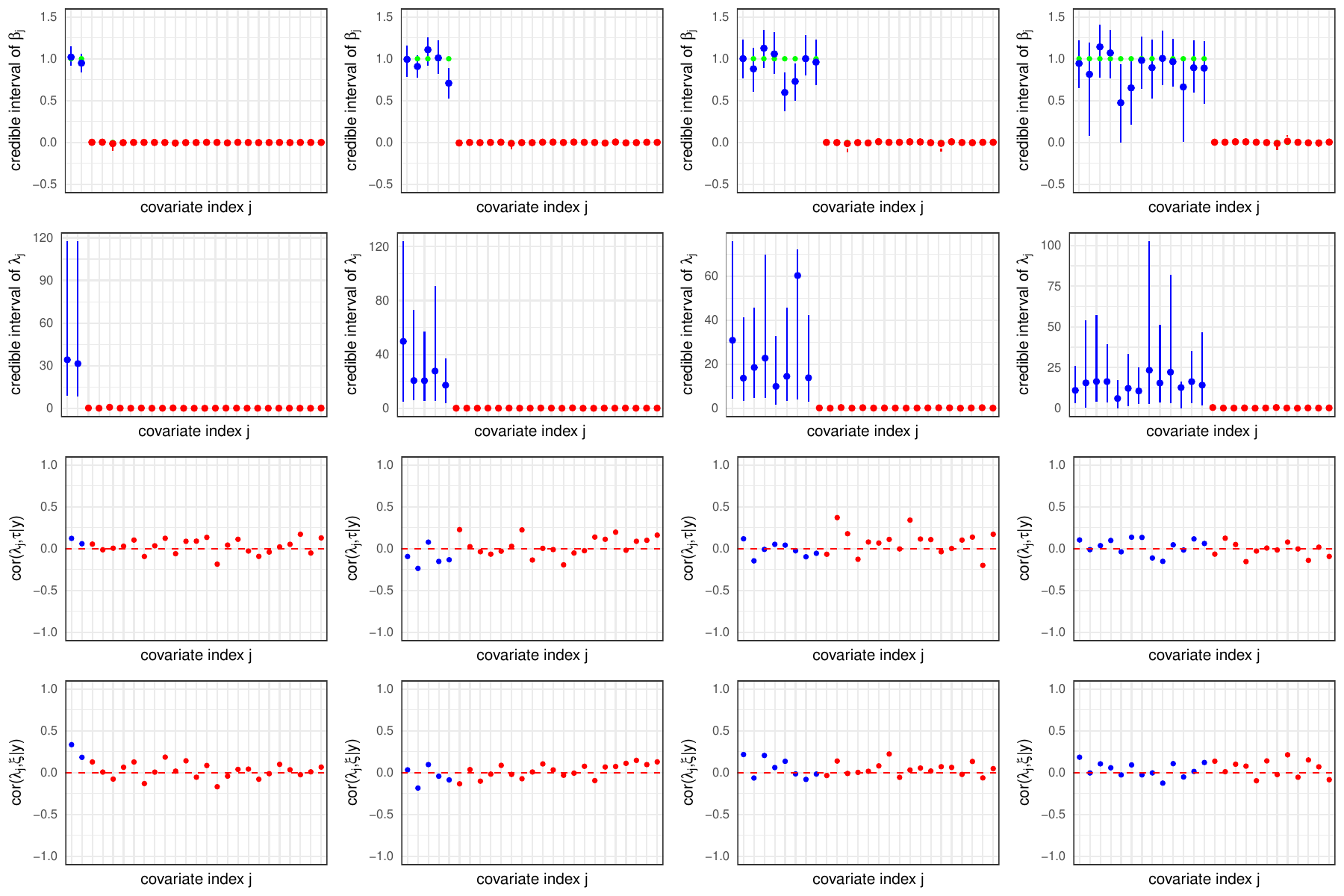}
	\caption{\baselineskip=10pt
		The results of posterior inference obtained using the GLT prior for the same datasets used in Figure \ref{fig:failure of Horseshoe}: $\mathcal{A}_{1}$ (first column), $\mathcal{A}_{2}$ (second column), $\mathcal{A}_{3}$ (third column), and $\mathcal{A}_{4}$ (fourth column). The posterior means of $(\tau,\xi)$ corresponding to the four datasets are $(0.003,2.010)$ ($\mathcal{A}_{1}$), $(0.004,2.134)$ ($\mathcal{A}_{2}$), $(0.004,2.235)$ ($\mathcal{A}_{3}$), and $(0.004,2.347)$ ($\mathcal{A}_{4}$), respectively.}
	\label{fig:signal_detection_of_GLT}
\end{figure}

We applied the GLT prior (\ref{eq:gpd_beta})--(\ref{eq:gpd_xi}) to the same four datasets $\mathcal{A}_{l}$ $(l=1,2,3,4)$ used in the previous subsection. The results of posterior inference are displayed in Figure \ref{fig:signal_detection_of_GLT}. Here, the posterior correlation between each local-scale parameter $\lambda_{j}$ and the shape parameter $\xi$ (i.e., $\{\text{cor}(\lambda_{j}, \xi | \textbf{y})\}_{j=1}^{p}$) is additionally plotted in the panels on the fourth row. It is worth noting that the GLT prior can successfully detect signals in the moderately sparse dataset $\mathcal{A}_{4}$, where the Horseshoe collapsed. Extensive numerical experiments are conducted in Section \ref{sec:Simulations}, demonstrating that the GLT prior performs well across diverse sparse regimes.

The posterior means of the shape parameter $\xi$ for the four datasets are $2.010$ $(\mathcal{A}_{1})$, $2.134$ $(\mathcal{A}_{2})$, $2.235$ $(\mathcal{A}_{3})$, and $2.347$ $(\mathcal{A}_{4})$, respectively. The monotonicity suggests that the tail of the GLT prior adapts to the sparsity level, \emph{a posteriori}, demonstrating the tail-adaptive shrinkage property.

We describe the signal detection mechanism of the GLT prior (\ref{eq:gpd_beta}) -- (\ref{eq:gpd_xi}). The GLT prior perceives the signal detection problem as the mirror image of the extreme value identification problem \citep{west1984outlier,polson2010shrink}. Under the formulation of the GLT prior, a selected (signal) coefficient $\beta_{j}$ is one whose corresponding local-scale parameter $\lambda_{j}$ is an extreme value potentially located in the tail part of the local-scale density $f(x|\tau,\xi) = \mathcal{GPD}(x|\tau,\xi)$ (\ref{eq:gpd_lambda}). Following conventional distributional theory, the global-scale parameter $\tau$ scales the local-scale density $f$, while the shape parameter $\xi$ controls the heaviness of the tail of $f$. 

The signal detection mechanism of the GLT prior is thus different from that of the Horseshoe. The former relies more on distribution theory and extreme value theory underlying the local-scale distribution, while the latter resorts to the competing nature between the local-scale parameter and the global-scale parameter, as described in Subsection \ref{subsec:The Horseshoe with varied sparsity level}. As observed from the panels on the third and fourth rows of Figure \ref{fig:signal_detection_of_GLT}, the estimates of $\tau$ and $\xi$ are nearly independent of the estimates of $\{\lambda_{j}\}_{j=1}^{p}$, \emph{a posteriori}. This contrasts with the panels on the third row of Figure \ref{fig:failure of Horseshoe}, where negative correlations are prevalent across all coefficients, and the selected signals are those with more distinguishable negative correlations.

\begin{figure}[h]
	\centering
	\includegraphics[width=\textwidth]{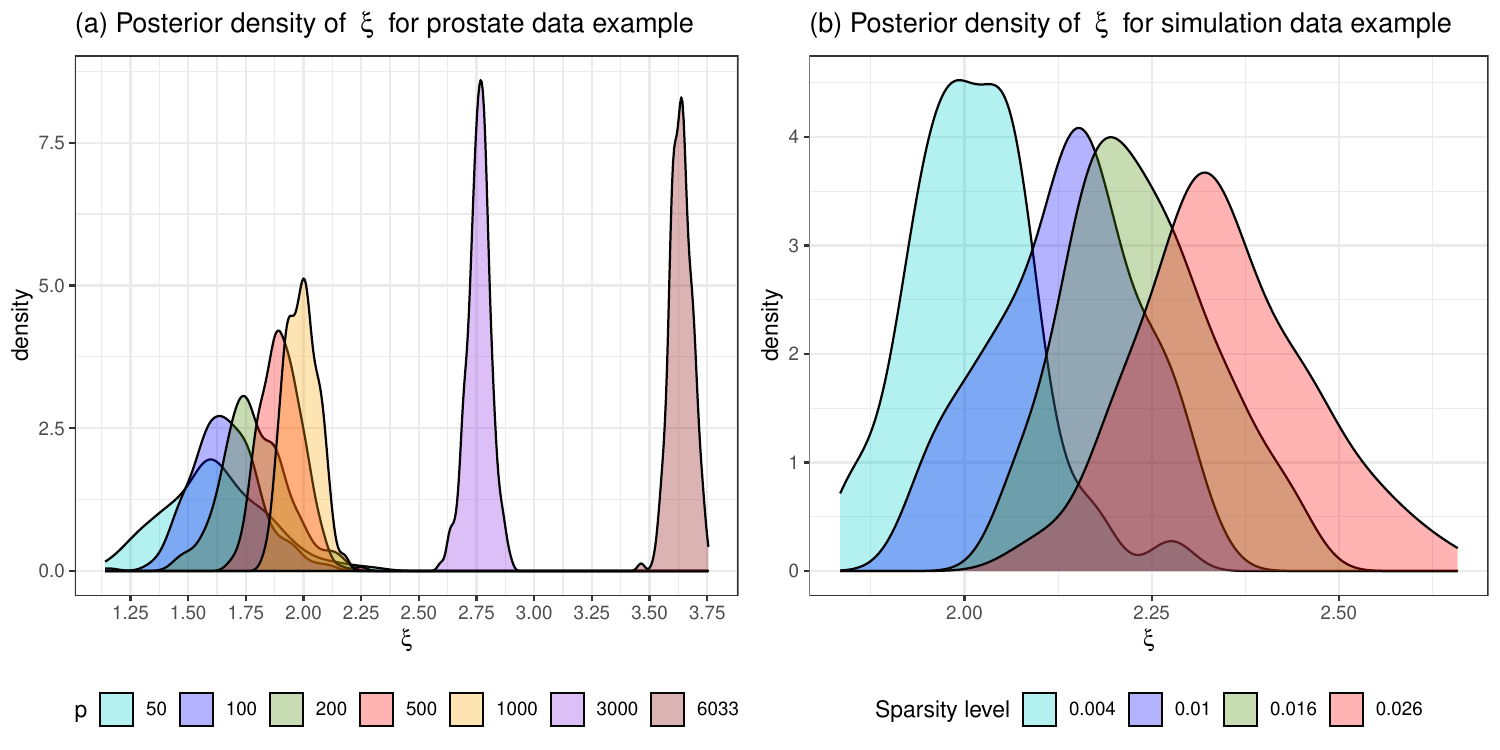}
	\caption{\baselineskip=10pt Posterior distribution of the shape parameter $\xi$ obtained using the GLT prior for the seven prostate cancer datasets $\mathcal{P}_{l}$ ($l=1, \cdots, 7$) (Panel (a)) and the four high-dimensional simulated datasets $\mathcal{A}_{l}$ ($l=1, \cdots, 4$) (Panel (b)). For the prostate cancer data example, the posterior means of $\xi$ are $1.620$ ($\mathcal{P}_{1}$) $<$ $1.662$ ($\mathcal{P}_{2}$) $<$ $1.789$ ($\mathcal{P}_{3}$) $<$ $1.905$ ($\mathcal{P}_{4}$) $<$ $1.991$ ($\mathcal{P}_{5}$) $<$ $2.760$ ($\mathcal{P}_{6}$) $<$ $3.636$ ($\mathcal{P}_{7}$). For the simulated data example, the posterior means of $\xi$ are $2.010$ ($\mathcal{A}_{1}$) $<$ $2.134$ ($\mathcal{A}_{2}$) $<$ $2.235$ ($\mathcal{A}_{3}$) $<$ $2.347$ ($\mathcal{A}_{4}$). }
	\label{fig:tail_adaptive_shrinkage_two_examples}
\end{figure}

Before proceeding to extensive simulation experiments, we explore the tail-adaptive shrinkage property of the GLT prior by examining the posterior distribution of the shape parameter \(\xi\). This examination is based on seven prostate cancer datasets \(\mathcal{P}_{l}\) (\(l=1, \cdots, 7\)) discussed in Subsection \ref{subsec:Prostate cancer data analysis via the GLT prior} and four high-dimensional simulated datasets \(\mathcal{A}_{l}\) (\(l=1,\cdots,4\)) discussed in the current subsection. Although the former uses a sparse normal mean model and the latter a high-dimensional regression model, both require sparse estimation robust against diverse sparsity regimes \citep{bhattacharya2015dirichlet}. Figure \ref{fig:tail_adaptive_shrinkage_two_examples} illustrates the posterior distribution of \(\xi\) (Panel (a) for the prostate cancer datasets and Panel (b) for the simulated datasets). Panel (b) shows that as the sparsity level increases from 0.004 to 0.026, the posterior distribution gradually shifts from left to right to accommodate a greater number of signals. In other words, the tail-heaviness of the local-scale density (\ref{eq:glt_prior_lambda}) is adaptive to the sparsity level, \emph{a posteriori}. The sparsity level of the prostate cancer datasets is unknown as it is real data, but it can be inferred that as more genes are considered, the sparsity level may increase. Panel (a) demonstrates that the posterior distribution shifts from left to right to accommodate a greater number of interesting genes, similar to what is observed in Panel (b). This behavior contrasts with the Horseshoe, where the tail behavior is pre-fixed (as proven in Proposition \ref{corollary:tail_index_of_Horseshoe_prior}).

\section{Simulations}\label{sec:Simulations}
\subsection{\textbf{Outline}}\label{subsec:Outline}
In Subsection \ref{subsec:Artificial high-dimensional data generator}, we illustrated the process of generating high-dimensional dataset $(\textbf{y}, \textbf{X}) \in \mathbb{R}^{n} \times \mathbb{R}^{n \times p}$ from high-dimensional regression (\ref{eq:sparse_high_dim_linear_regression}) under a given simulation environment $(n, p, q, \varrho, \text{SNR})$. Recall that, we assumed a true coefficient vector to be $\bm{\beta}_{0} = (\beta_{0,1}, \cdots, \beta_{0,q},$ $ \beta_{0,q+1}, \cdots, \beta_{0,p})^{\top} \in \mathbb{R}^{p}$, where $\beta_{0,1} = \cdots = \beta_{0,q} = 1$ and $\beta_{0,q+1} = \cdots = \beta_{0,p} = 0$. In this section, we conduct extensive numerical experiment to compare the performances of the Horseshoe (\ref{eq:HS_beta}) -- (\ref{eq:HS_lambda_tau}) and the GLT prior (\ref{eq:gpd_beta}) -- (\ref{eq:gpd_xi}) under three different scenarios. We set the default environmental values as $(n, p) \in \{(100, 500), (200, 1000), (300, 1500)\}$, $s = q/p = 0.02$, $\text{SNR} = 5$, and $\varrho = 0$. Then, we separately consider the following three scenarios by varying one environmental value while keeping the others fixed;
\begin{itemize}
	\item[]\textit{\textbf{Scenario 1}}: 
	varied sparsity level $s=q/p$ from $0.0001$ to $0.1$,
	\item[]\textit{\textbf{Scenario 2}}:  varied $\varrho$ from $0$ to $0.8$,
	\item[]\textit{\textbf{Scenario 3}}: varied $\text{SNR}$ from $2$  to $10$.
\end{itemize}

In general, most variable selection methods are ideally designed to work well when the sparsity level and the degree of co-linearity are small (i.e., small values for $s$ and $\varrho$), while the SNR is relatively large. The key to robust estimation is that, even in opposite situations where $s$ and $\varrho$ are relatively large and SNR is small (resulting in nearly ill-posed high-dimensional data and potential contamination by noise or co-linearity between predictors), variable selection methods should still provide reasonable estimation accuracy. 

Regarding performance metrics, we report the medians of the mean squared error (MSE), and  MSEs separately for the signal and noise coefficients, measured from 100 replicated datasets. Let $\widehat{\bm{\beta}} = (\widehat{\beta}_{1}, \cdots, \widehat{\beta}_{p} )^{\top} \in \mathbb{R}^{p}$ denote the posterior mean obtained using either the Horseshoe or the GLT prior. The MSE, MSE for signals and noises are defined as follows:

\begin{equation}
	\nonumber
	\text{MSE}= 
	\frac{1}{p}
	\sum_{j = 1}^{p}
	( \widehat{\beta}_{j} - \beta_{0,j} )^{2}\text{,}\quad
	\text{MSE}_{\text{S}} = \frac{1}{q}\sum_{j = 1}^{q} (\widehat{\beta}_{j} - 1)^{2},\quad 
	\text{MSE}_{\text{N}} = \frac{1}{p - q}\sum_{j = q + 1}^{p} (\widehat{\beta}_{j})^{2}.
\end{equation}

The MSE measures the overall accuracy of estimation for the coefficients induced by a prior, where a lower value indicates better accuracy. The MSE can be dissected into two components: (1) MSE for the signal part ($\text{MSE}_{\text{S}}$), which measures signal recovery ability, and (2) MSE for the noise parts ($\text{MSE}_{\text{N}}$), which measures noise shrinking ability. It is important to emphasize that if the estimate collapses ($\widehat{\beta}_{j} \approx 0$ for all $j$), then $\text{MSE}_{\text{S}}$ and $\text{MSE}_{\text{N}}$ will be close to $1$ and $0$, respectively.

\subsection{\bf{Scenario 1: varied sparsity level $s = q/p$}}\label{subsec:Scenario 1}
\begin{figure}[h]
	\centering
	\includegraphics[width=\textwidth]{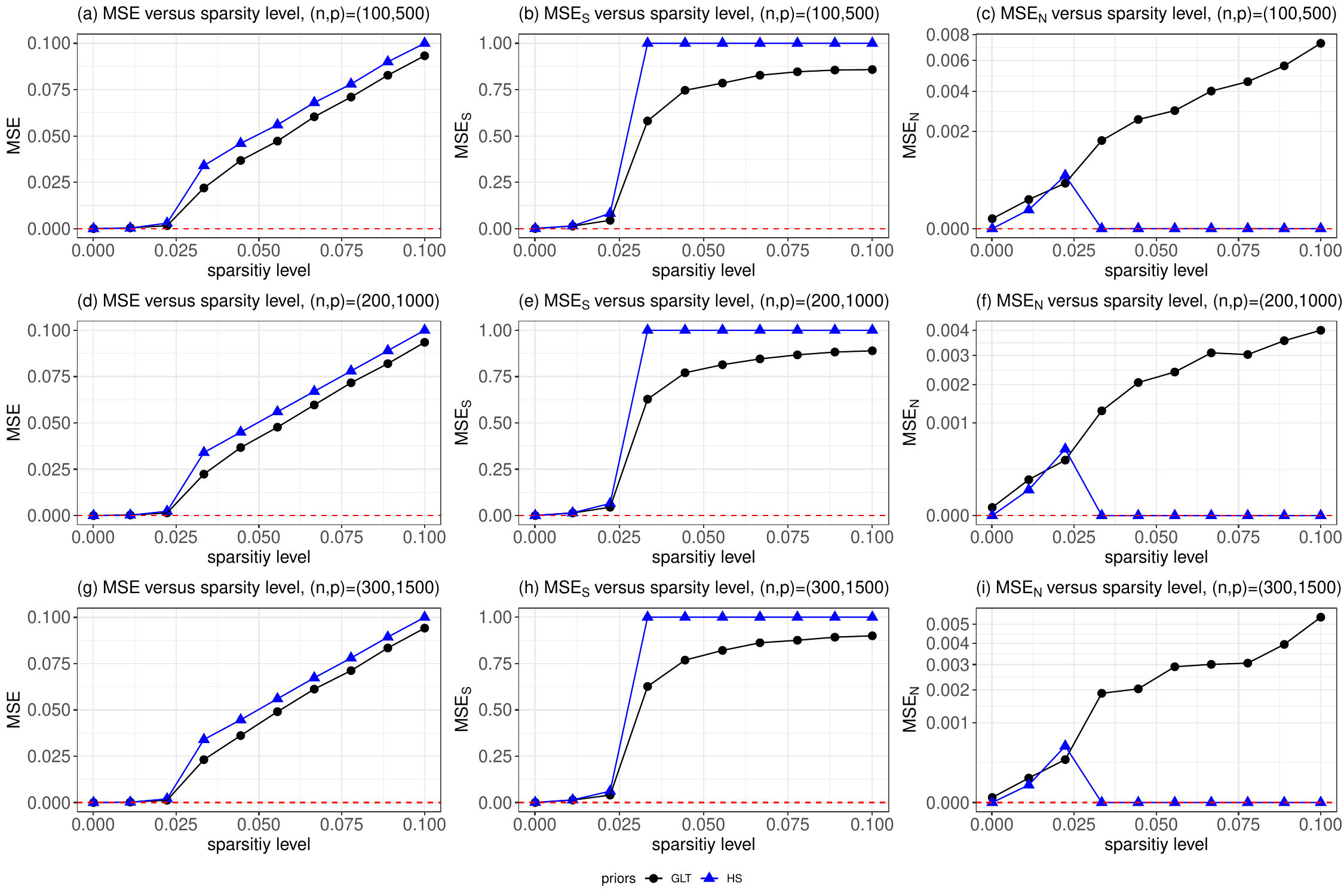}
	\caption{\baselineskip=10pt
		The medians of $\text{MSE}$, $\text{MSE}_{\text{S}}$, and $\text{MSE}_{\text{N}}$ from 100 replications across varied sparsity level $s$ from 0.0001 to 0.1: $(n,p)=(100,500)$ (top panels), $(n,p)=(200,1000)$ (middle panels), and $(n,p)=(300,1500)$ (bottom panels). Results obtained using the GLT prior and Horseshoe are indicated by a black circle dot ($\bullet$) and blue triangle dot ($\blacktriangle$). The red dotted horizontal line represents zero.
	}
	\label{fig:S1_simulation_results}
\end{figure}

Figure \ref{fig:S1_simulation_results} presents the medians  of $\text{MSE}$, $\text{MSE}_{\text{S}}$, and $\text{MSE}_{\text{N}}$ under Scenario 1. The top, middle, and bottom panels correspond to the simulation results with $(n,p) = (100,500)$, $(n,p) = (200,1000)$, and $(n,p) = (300,1500)$, respectively. Specifically, the top panels corresponds to the setting $(n=100$, $p=500$, $q$, $\varrho=0$, $\text{SNR}=5)$ with the number of signals ($q$) taking values from the set $\{1,  6, 12, 17, 23, 28, 34, 39, 45, 50\}$. Similarly, the middle panels corresponds to $(n=200$, $p=1000$, $q$, $\varrho=0$, $\text{SNR}=5)$ with $q$ ranging from $\{1,  12,  23,  34,  45,  56,  67,  78,  89, 100\}$. The bottom panels corresponds to $(n=200$, $p=1000$, $q$, $\varrho=0$, $\text{SNR}=5)$ with $q$ ranging from $\{1,  17, 34,  51,  67,  84, 101,$  $ 117, 134 ,150\}$. That way, sparsity level $s =q/p$ varies from 0.0001 to 0.1.

Observing the results from Panels (a), (d), and (g), we find that both the Horseshoe and GLT prior yield excellent performance under the ultra-sparse regime and up to a sparsity level of approximately $s = 0.022$, regardless of the sample size $n$ and number of predictors $p$. However, beyond this point, the values of MSE between the two priors suddenly diverge, and the MSE of the Horseshoe estimator stays above the MSE of the GLT estimator, indicating the superior performance of the GLT prior over the Horseshoe under a moderate sparse regime. This sudden gap is mainly due to the sharp underestimation of the global-scale parameter of the Horseshoe, causing the collapsing behavior of the Horseshoe estimator, as exemplified in Subsection \ref{subsec:The Horseshoe with varied sparsity level}. See Panels (c), (f), and (i), where the value of $\text{MSE}_{\text{N}}$ for Horseshoe beyond the sparsity level $s = 0.03$ is virtually zero. Overall, the simulation experiments demonstrate that the GLT prior provides more robust estimation across a wide range of sparsity levels than the Horseshoe, indicating the superiority of the varying tail rule over the fixed tail rule in the variable selection problem.

\subsection{\bf{Scenario 2: varied $\varrho$}}\label{subsec:Scenario 2}
\begin{figure}[h]
	\centering
	\includegraphics[width=\textwidth]{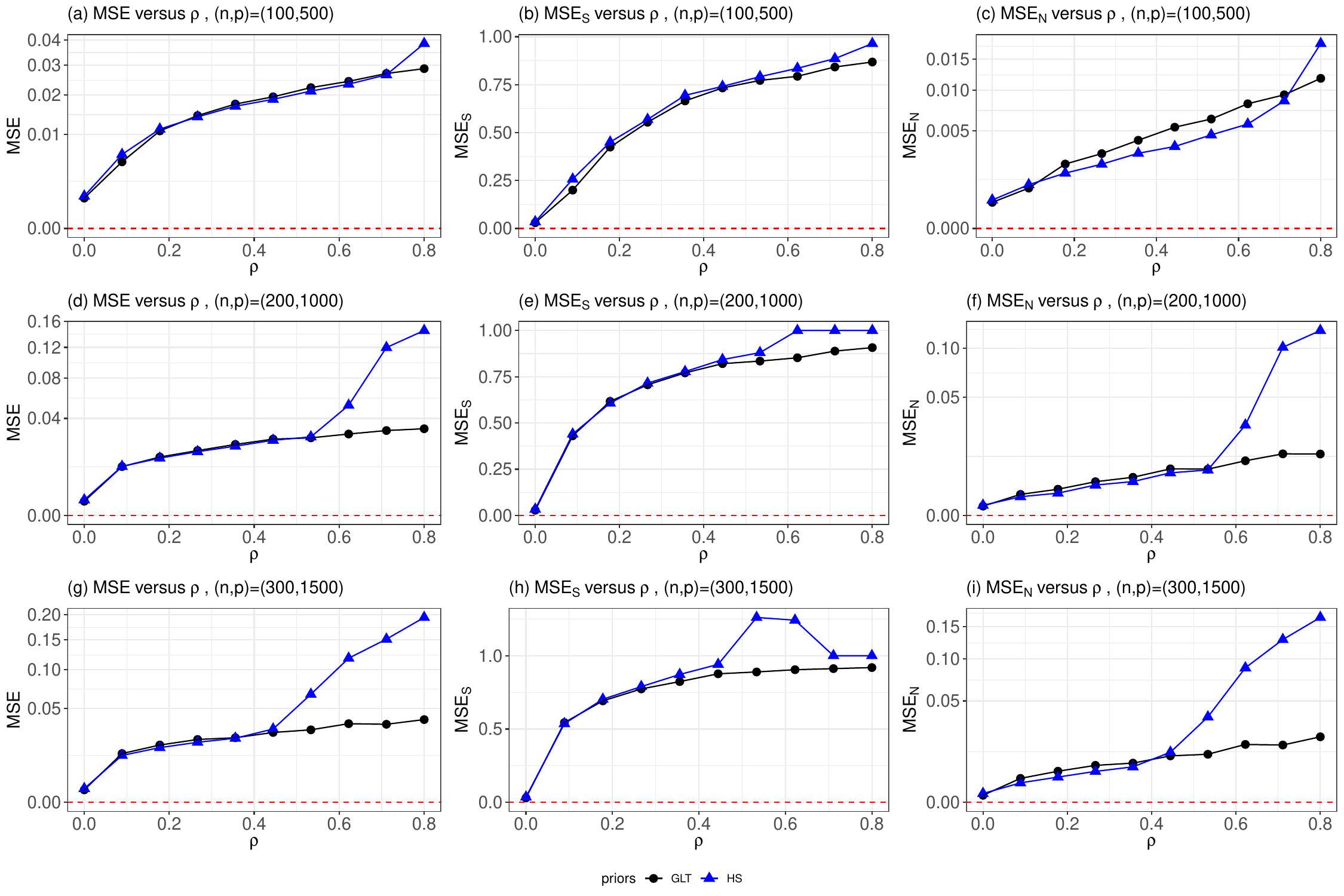}	
	\caption{\baselineskip=10pt
		The medians of $\text{MSE}$, $\text{MSE}_{\text{S}}$, and $\text{MSE}_{\text{N}}$ from 100 replications across different values of $\varrho$ from 0 to 0.8: $(n,p)=(100,500)$ (top panels), $(n,p)=(200,1000)$ (middle panels), and $(n,p)=(300,1500)$ (bottom panels).
	}
	\label{fig:S2_simulation_results}
\end{figure}

Figure \ref{fig:S2_simulation_results} displays the medians of $\text{MSE}$, $\text{MSE}_{\text{S}}$, and $\text{MSE}_{\text{N}}$ under Scenario 2. The top, middle, and bottom panels correspond to simulation results with $(n,p) = (100,500)$, $(n,p) = (200,1000)$, and $(n,p) = (300,1500)$, where the number of true predictors is $q = 10, 20, $ and $30$, respectively. 

Observing the results from Panels (a), (d), and (g), we find that both the Horseshoe and GLT prior perform well when $\varrho$ is very small. However, as $\varrho$ increases, MSE also increases. This is expected because introducing multicollinearity in linear regression leads to increased standard errors, resulting in less precise estimation of coefficients. In this simulation, it is notable that the increase in MSE of the Horseshoe estimator becomes particularly dramatic beyond a certain point—approximately $\varrho = 0.711$ for $(n,p) = (100,500)$, $\varrho = 0.533$ for $(n,p) = (200,1000)$, and $\varrho = 0.444$ for $(n,p) = (300,1500)$.  This indicates that the Horseshoe estimator becomes more vulnerable to multicollinearity as the number of predictors $p$ increases. In contrast, the increase in MSE of the GLT estimator remains relatively stable across different values of $\varrho$. Overall, the simulation experiments demonstrate that the GLT prior provides more robust estimation across a wide range of degrees of multicollinearity compared to the Horseshoe.

\subsection{\bf{Scenario 3: varied signal-to-noise ratio}}\label{subsec:Scenario 3}
\begin{figure}[h]
	\centering
	\includegraphics[width=\textwidth]{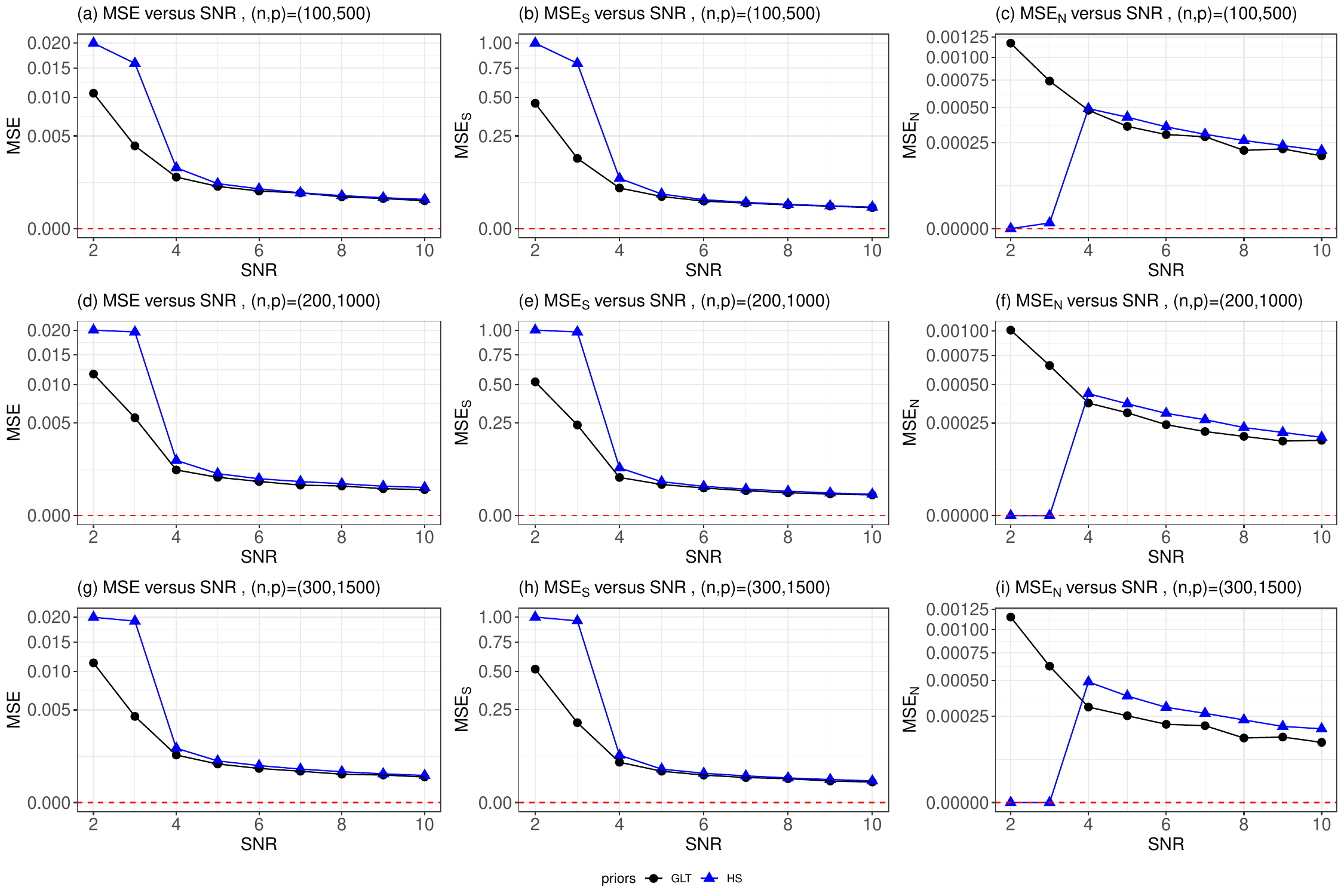}
	\caption{\baselineskip=10pt
		The medians of $\text{MSE}$, $\text{MSE}_{\text{S}}$, and $\text{MSE}_{\text{N}}$ from 100 replications across different values of $\text{SNR}$ from 2 to 10: $(n,p)=(100,500)$ (top panels), $(n,p)=(200,1000)$ (middle panels), and $(n,p)=(300,1500)$ (bottom panels).}
	\label{fig:S3_simulation_results}
\end{figure}


Figure \ref{fig:S3_simulation_results} presents the medians of $\text{MSE}$, $\text{MSE}_{\text{S}}$, and $\text{MSE}_{\text{N}}$ under Scenario 3. The top, middle, and bottom panels correspond to simulation results with $(n,p) = (100,500)$, $(n,p) = (200,1000)$, and $(n,p) = (300,1500)$, where the number of true predictors is $q = 10, 20, $ and $30$, respectively. In Panels (a), (d), and (g), it is observed that the MLE values of both the Horseshoe and GLT prior decrease as $\text{SNR}$ increases. This is expected because a higher value of $\text{SNR}$ leads to less contamination attributed to the error perturbation term (i.e., $\bm{\epsilon}$), making signal detection more trivial. Notably, when $\text{SNR}< 4$, the gap between the MSE values of the Horseshoe estimator and the GLT estimator is significant, indicating that the Horseshoe estimator is more vulnerable to a low SNR situation than the GLT estimator. Observing Panels (c), (f), and (i), this vulnerability is due to the collapsing behavior of the Horseshoe estimator when $\text{SNR}$ is extremely small. Overall, the simulation experiments demonstrate that the GLT prior outperforms the Horseshoe across different values of $\text{SNR}$ and provides more robust estimation, even when $\text{SNR}$ is quite small.

In summary, simulation experiments under Scenarios 1, 2, and 3 indicate that the tail-adaptive shrinkage property of the GLT prior provides robust estimation under moderate sparsity, higher collinearity in the design matrix, and low signal-to-noise cases, where the Horseshoe estimator may be sub-optimal.

\section{Discussion}\label{sec:Discussion}
In this paper, we proposed a new framework of shrinkage priors called global-local-tail shrinkage priors, which is an extension of global-local shrinkage priors, for high-dimensional regression problems. The main objective of using global-local-tail shrinkage priors is to provide robust Bayesian inference across diverse sparsity regimes, ranging from ultra-sparse to moderately-sparse situations. This is in contrast to the Horseshoe or other global-local shrinkage priors, which are primarily designed to yield meaningful inference results in ultra-sparse regimes. Simulation results from Section \ref{sec:Simulations} show that the advantage of using global-local-tail shrinkage priors can be further extended to problems where multicollinearity is present and the signal-to-noise ratio is small, situations where existing global-local shrinkage priors may not work well due to the fixed tail behavior. The optimal convergence properties of the GLT posterior emphasize the significant role of the shape parameter $\xi$ in sparse estimation across various sparsity regimes. This finding is further supported by real gene expression analysis and simulation studies. 


We emphasize that delicate care is required to estimate the shape parameter within the global-local-tail shrinkage framework, which we regard as one of the salient contributions of the paper. For the GLT prior, we propose an algorithm that combines the elliptical slice sampler \citep{murray2010elliptical} with the Hill estimator \citep{hill1975simple} from extreme value theory, eliminating the need for tuning any hyperparameters. This automatic tuning enables adaptive learning of the shape parameter $\xi$ according to the unknown sparsity level. The Supplementary Material includes further exploration of the performance of the GLT prior compared to other shrinkage priors, as well as its application to real gene expression data analysis and curve fitting studies.

\begin{appendix}
\section*{Appendix}
\setcounter{section}{0}
\setcounter{equation}{0}\renewcommand\theequation{A\arabic{equation}}
\renewcommand{\thesection}{A.\arabic{section}}
\renewcommand{\thesubsection}{A.\arabic{section}.\arabic{subsection}}

\section{Posterior computation}\label{sec:Posterior computation}
\subsection{\textbf{Gibbs sampler}}\label{subsec:Gibbs sampler}
Current Section provides a full description for posterior computation to implement a MCMC sampling algorithm to apply the GLT prior under the sparse high-dimensional linear regression. Consider $\textbf{y} = \textbf{X}\bm{\beta}+ \sigma\bm{\epsilon} $ , $\bm{\epsilon}\sim \mathcal{N}_{n}(\textbf{0},\textbf{I}_{n})$ (\ref{eq:sparse_high_dim_linear_regression}), $\sigma^2 \sim \pi(\sigma^2) \propto	1/\sigma^2$, and $\bm{\beta} \sim \pi_{\text{GLT}}(\bm{\beta})$ (\ref{eq:gpd_beta}) -- (\ref{eq:gpd_xi}). Let $\bm{\Omega} =(\bm{\beta}, \sigma^{2}, \bm{\lambda}, \tau, \xi)\in \mathbb{R}^{p}\times (0,\infty)\times (0,\infty)^{p} \times (0,\infty) \times(1/2, \infty)$ denote all the latent random variables. Our purpose is then to sample from the full joint posterior distribution, $\pi(\bm{\Omega}| \textbf{y} )$, which is proportional to 
\begin{small}
	\begin{align*}
		&\mathcal{N}_{n}(\textbf{y} |\textbf{X}\bm{\beta} , \sigma^{2} \textbf{I}_n)
		\mathcal{N}_{p}(\bm{\beta}|\textbf{0}, \sigma^{2} \bm{\Lambda})
		\pi(\sigma^{2})
		\bigg\{
		\prod_{j=1}^{p}
		\pi(\lambda_{j} | \tau, \xi)
		\bigg\}
		\pi(\tau, \xi)
		\\
		&\quad\propto
		\mathcal{N}_{n}(\textbf{y} |\textbf{X}\bm{\beta}, \sigma^{2} \textbf{I}_n)
		\mathcal{N}_{p}(\bm{\beta}|\textbf{0}, \sigma^{2} \bm{\Lambda})
		\pi(\sigma^{2})
		\bigg\{
		\prod_{j=1}^{p}
		\mathcal{GPD} (\lambda_{j} | \tau, \xi)
		\bigg\}\\
		&\quad \quad \mathcal{IG}(
		\tau | p/\xi + 1, 1)
		\log\ \mathcal{N} (\xi | \mu, \rho^{2})
		\mathcal{I}_{(1/2,\infty)}(\xi),\quad \bm{\Lambda} = \text{diag}(\lambda_{1}^{2}, \cdots, \lambda_{p}^{2}) \in \mathbb{R}^{p\times p}.
	\end{align*}
\end{small}
Note that the $\mu \in \mathbb{R}$ and $\rho^{2}>0$ are hyper-parameters, which typically requires an expert-tuning. In Subsection \ref{subsec:Hyperparameter specification}, we provide an automated hyper-parameter turning algorithm to obviate such tuning.

Since the full joint posterior distribution $\pi(\bm{\Omega}|y)$ is not in a closed form, we develop a Gibbs sampler \citep{casella1992explaining,lee2021gibbs} to sample from the joint density. A single cycle of the Gibbs sampling algorithm comprises the following \textbf{\emph{Step 1 - Step 5}}.
\begin{enumerate}
	\item[] \emph{\textbf{Step 1.} } Sample $\bm{\beta}$ from conditional posterior 
	\begin{align*}
		\pi(\bm{\beta} | - ) & \sim
		\mathcal{N}_{p}
		(
		\bm{\Sigma} \textbf{X}^{\top} \textbf{y}, \sigma^{2} \bm{\Sigma}
		),\quad
		\bm{\Sigma} = 
		(
		\textbf{X}^{\top}\textbf{X}+ \bm{\Lambda}^{-1}
		)^{-1} \in \mathbb{R}^{p \times p}.
	\end{align*}
	\item[]\emph{\textbf{Step 2.} }
	Sample $\sigma^{2}$ from conditional posterior
	\begin{align*}
		\pi(\sigma^{2}| - ) 
		& \sim
		\mathcal{IG}\bigg(
		\frac{n+p}{2},
		\frac{\|
			\textbf{y} - \textbf{X}\bm{\beta } 
			\|_{2}^{2} + \bm{\beta } ^{\top}\bm{\Lambda}^{-1} \bm{\beta } }{2}
		\bigg).
	\end{align*}
	\item[]\emph{\textbf{Step 3.} }
	Update $\lambda_{j}$, $j = 1,\cdots, p$, independently using slice sampler \citep{neal2003slice} within the Gibbs sampler. Proportional part of full conditional posterior is
	\begin{align}
		\label{eq:full_condi_lambda}
		\pi(\lambda_{j}| - ) & \propto
		\frac{1}{\lambda_{j}}
		\exp\bigg(- 
		\frac{\beta_{j}^{2}}{2 \sigma^{2}\lambda_{j}^{2}}
		\bigg)
		\cdot
		\bigg(1 + 
		\frac{\xi  \lambda_{j}}{\tau}
		\bigg)^{-(1/\xi +1)}.
	\end{align}
	\item[]\emph{\textbf{Step 4. }}
	Update $\tau$ using slice sampler \citep{neal2003slice} within the Gibbs sampler. Proportional part of full conditional posterior is
	\begin{align}
		\label{eq:full_condi_tau}
		\pi(\tau | - ) & \propto
		\tau^{-2}
		\exp
		(
		-1/\tau
		)
		\cdot
		\prod_{j=1}^{p}
		(\tau + \xi \lambda_{j})^{-(1/\xi + 1)}.
	\end{align}
	\item[]\emph{\textbf{Step 5. }}
	Update $\xi$ using elliptical slice sampler \citep{murray2010elliptical} after variable change $\eta = \log\ \xi$ within the Gibbs sampler. Proportional part of full conditional posterior is
	\begin{align}
		\label{eq:slice_sampler_1}
		\pi(\xi | - ) & \propto
		\mathcal{V}_{p}(\xi)
		\cdot
		\log\ \mathcal{N}_{1}(\xi|\mu, \rho^{2})
		\cdot
		\mathcal{I}_{(1/2, \infty)}(\xi),
	\end{align}
	where $\mathcal{V}_{p}(\xi)
	=
	\{\Gamma(p/\xi + 1)\}^{-1}\pi^{p/2} \prod_{j=1}^{p}r_{j}(\xi)
	$ \\ with $r_{j}(\xi) =$ $ (\tau + \xi \lambda_{j})^{-(1/\xi + 1 )}$, $j = 1, \cdots, p$. \end{enumerate}

\subsection{\textbf{Slice sampler implementation in \emph{Step 3}  and \emph{Step 4}}}
\label{subsec:Slice sampler implementation}
Slice sampler \citep{neal2003slice} is a popular technique to adapt the step-size of a MCMC algorithm and is based on the local property of the target density.  The basic idea underlying the slice sampler is parameter expansion which involves intentional introduction of auxiliary variables \citep{damlen1999gibbs}. Finding an appropriate parameter expansion depends on specific functional form of the target density. 

Let $j \in  \{1,\cdots,p\}$. To implement the slice sampler in the \textbf{\emph{Step 3}} (\ref{eq:full_condi_lambda}), first use change of variable, $\gamma_{j} = \lambda_{j}^{2}$, to get
\begin{align}
	\nonumber
	\pi(\gamma_{j}| - ) 
	& \propto
	\gamma_{j}^{-1}\exp\ (-m_{j}/\gamma_{j})
	\cdot
	(\tau + \xi \sqrt{\gamma_{j}})^{-(1/\xi +1)}\\
	\nonumber
	&= 
	\gamma_{j}^{-1}\exp\ (-m_{j}/\gamma_{j})
	\cdot
	(\sqrt{\gamma_{j}})^{-(1/\xi+1)}
	(\sqrt{\gamma_{j}})^{(1/\xi+1)}
	\cdot
	(\tau + \xi \sqrt{\gamma_{j}})^{-(1/\xi +1)}\\
	\nonumber
	&=
	\gamma_{j}^{- (1/\xi +1 )/2 -1 }\exp\ (-m_{j}/\gamma_{j})
	\cdot
	(\xi + \tau \cdot \gamma_{j}^{-1/2})^{-(1/\xi +1)}\\
	\label{eq:slice_sampler_lambda_j_1}
	&\propto
	\mathcal{IG}\{\gamma_{j}| (1/\xi +1 )/2, m_{j}\}
	\cdot
	g(\gamma_{j}),
\end{align}
where $m_{j}=\beta_{j}^{2}/(2 \sigma^{2})$ and $g(\gamma_{j}) = (\xi + \tau \cdot \gamma_{j}^{-1/2})^{-(1/\xi +1)}$. Note that the function $u_{j} = g(\gamma_{j})$ is increasing on $(0,\infty)$, and its inverse function is $\gamma_{j} = g^{-1}(u_{j}) = [\tau/\{ u_{j}^{-(\xi/(1+\xi))} - \xi \} ]^{2} $. Now, consider a density, 
$
\pi(\gamma_{j}, u_{j}| - ) 
\propto
\mathcal{IG}\{\gamma_{j}| (1/\xi +1 )/2, m_{j}\}
\cdot
\mathcal{I}_{ (0, g(\gamma_{j})) }(u_{j}).
$
Then we can show that $\int \pi(\gamma_{j}, u_{j}| - ) du_{j} = \pi(\gamma_{j}| - ) $, which means that $\pi(\gamma_{j}, u_{j}| - )$ is a valid parameter expansion of (\ref{eq:slice_sampler_lambda_j_1}). Actual sampling is executed on $\pi(\gamma_{j}, u_{j}| - ) $ using the Gibbs sampler: (i)
$u_{j}|\gamma_{j}, - \sim \pi(u_{j}|\gamma_{j}, - ) = \mathcal{U}(0, g(\gamma_{j}))$ and (ii) $
\gamma_{j}|u_{j}, - \sim \pi(\gamma_{j}| u_{j}, - ) = $ $
\mathcal{IG}\{\gamma_{j}| (1/\xi +1 )/2,$ $ m_{j}\}\cdot
\mathcal{I}_{ (g^{-1}(u_{j}), \infty) }(\gamma_{j}). 
$ After the Gibbs sampling, transform back to $\lambda_{j} = \sqrt{\gamma_{j}}$.

To implement the slice sampler in the \emph{Step 4}, note from  (\ref{eq:full_condi_tau}):
\begin{align}
	\label{eq:slice_sampler_tau_1}
	\pi(\tau| - ) 
	& \propto
	\mathcal{IG}(\tau| 1,1)
	\cdot
	\prod_{j=1}^{p}
	g_{j}(\tau),
\end{align}
where $g_{j}(\tau) = (\tau + \xi \lambda_{j})^{-(1/\xi + 1)}$, $j=1,\cdots,p$. Note that $p$-functions $v_{j} = g_{j}(\tau)$, $j=1,\cdots,p$, are decreasing on $(0,\infty)$, and their inverse functions are $\tau = g_{j}^{-1}(v_{j}) = v_{j}^{-(\xi/(1+\xi))} -\xi \lambda_{j}$, $j=1,\cdots,p$. Now, consider a density:
$\pi(\tau, v_{1},$ $\cdots, v_{p}| - ) 
$ $\propto
\mathcal{IG}(\tau| 1,1)
\cdot
\prod_{j=1}^{p}
\mathcal{I}_{(0, g_{j}(\tau))}(v_{j}).
$

Then we have $\int $ $\cdots \int \pi(\tau, v_{1},$ $\cdots, v_{p}| - )  dv_{1}\cdots dv_{p} = \pi(\tau| - ) $ and hence $\pi(\tau, v_{1},$ $\cdots, v_{p}| - )$ is a valid parameter expansion of (\ref{eq:slice_sampler_tau_1}). Actual sampling is executed on $\pi(\tau, v_{1},\cdots, v_{p}| - ) $ using the Gibbs sampler:
\begin{align}
	\label{eq:slice_sampler_tau_3}
	v_{j}|\tau, v_{-j}, - &\sim \pi(v_{j}|\tau, v_{-j}, - ) = \mathcal{U}(0, g_{j}(\tau)), \quad (j =1, \cdots,p),\\
	\nonumber
	\tau|v_{1},\cdots,v_{p}, - &\sim 
	\mathcal{IG}(\tau| 1,1)
	\cdot
	\mathcal{I}_{ ( 0, \text{min}\{g_{1}^{-1}(v_{1}), \cdots, g_{p}^{-1}(v_{p}) \} )}(\tau),
\end{align}
where in (\ref{eq:slice_sampler_tau_3}), $v_{-j}$ represents the collection of $\{ v_{j}\}_{j=1}^{p}$ except for $v_{j}$. Note also  that each full conditional posterior distribution $\pi(v_{j}|\tau, v_{-j}, - )$, $j=1,\cdots,p$, does not depend on $v_{-j}$, i.e., $\pi(v_{j}|\tau, v_{-j}, - ) = \pi(v_{j}|\tau, - )$ and hence it is possible to parallelize the update of $\{v_{j} \}_{j=1}^{p}$.
\subsection{\textbf{Summary of the Hill estimator}}
\label{subsec:Summary of the Hill estimator}
We briefly explain the Hill estimator which plays a central role in hyper-parameter specification of the $\mu$. For notational coherence, we use the Greek letter $\lambda$ to describe a random quantity. Suppose that $\bm{\lambda}=(\lambda_{1},\cdots,\lambda_{p})^{\top} \in (0,\infty)^{p}$ is $p$-dimensional random variables from a strongly stationary process whose marginal distribution is $F$ such that its tail distribution is regularly varying with the tail-index $1/\xi$ with $\xi>0$ (hence, the corresponding shape parameter is $\xi$). By the Karamata's characterization theorem \citep{karamata1933mode}, the tail distribution (survival function) is described as $\bar{F}(\lambda)  = 1  - F(\lambda)= L(\lambda) \cdot \lambda^{-1/\xi}$ for some $\xi>0$ where $L$ is a slowly varying function \citep{resnick1995consistency,drees2000make}. Denote its order statistics with $\lambda_{(1)} \geq \cdots\geq\lambda_{(p)}$.

The Hill estimator \citep{hill1975simple} is a well-known estimator of shape parameter $\xi$ principled on the peaks-over-threshold methods. Hill estimator is obtained from the $k$ upper order statistics:
\begin{align}
	\label{eq:Hill_estimator}
	\widehat{\xi}_{k}(\bm{\lambda})
	&=
	\frac{1}{k-1}
	\sum_{j=1}^{k-1}
	\log\ \bigg( \frac{\lambda_{(j)}}{\lambda_{(k)}}\bigg)
	,
	\quad \text{for } 2 \leq k\leq p.
\end{align}

It is known that the Hill estimator (\ref{eq:Hill_estimator}) is a consistent estimator for $\xi$, i.e., $\widehat{\xi}_{k}(\bm{\lambda}) \rightarrow \xi$ in probability, if $p\rightarrow \infty$, $k \rightarrow \infty$, and $k/p \rightarrow 0$ \citep{resnick1995consistency,drees2000make,embrechts2013modelling}. Empirically it is known that the Hill estimator may work effectively when $F$ is of a Pareto type \citep{drees2000make,lee2018exponentiated}. (See Fig $1$ in \citep{drees2000make}.)

Suppose we have $p$ number of observations $\bm{\lambda} = (\lambda_{1},\cdots,\lambda_{p})^{\top}$, possibly generated from the aforementioned heavy distribution $F$. In practice, the Hill estimator is used as follows. First, calculate the estimator $\widehat{\xi}_{k}(\bm{\lambda})$ at each integer $k \in \{2, \cdots, p \}$, and then plot the ordered pairs $\{(k, \widehat{\xi}_{k}(\bm{\lambda}))\}_{k=2}^{p}$: the resulting plot is called the Hill plot (See the Figure 6.4.3 of \citep{embrechts2013modelling}). Then, select value(s) from the set of Hill estimators $\{\widehat{\xi}_{k}(\bm{\lambda})\}_{k=2}^{p}$ which are stable (roughly constant) with respect to $k$: then, such stable value(s) are regarded as reasonable estimate(s) for the shape parameter $\xi$ \citep{drees2000make}. Typically, the Hill plot may display high variability when $k$ is close to $2$ or $p$. As a practical remedy, one may disregard the first or last few of the estimates: the values $\widehat{\xi}_{k}(\bm{\lambda})$ that are evaluated at integers $k \in \{k_{\text{L}},\cdots, k_{\text{U}} \}$, $2 < k_{\text{L}} < k_{\text{U}} < p$, are considered to be monitored where the integers $ k_{\text{L}}$ and $ k_{\text{U}}$ are designated by user.
\subsection{\textbf{Hyper-parameter specification of $\mu$ and $\rho^{2}$}}\label{subsec:Hyperparameter specification}
Suppose we are at the \emph{\textbf{Step 5}} of the $s$-th iteration of the Gibbs sampler described in Subsection \ref{subsec:Gibbs sampler}. At this moment, we have already acquired posterior realizations, $\bm{\lambda}^{(s+1)}=(\lambda_{1}^{(s+1)},\cdots,\lambda_{p}^{(s+1)})^{\top}$ and $\tau^{(s+1)}$, that had been sampled from the previous steps, \emph{\textbf{Step 3}} and \emph{\textbf{Step 4}}, respectively. 

By treating the indicator $\mathcal{I}_{(1/2, \infty)}(\xi)$ in (\ref{eq:slice_sampler_1}) as a part of likelihood, we consider sampling $\xi^{(s+1)}$ from the density;
\begin{align}
	\label{eq:mu_spec_1}
	\xi^{(s+1)} \sim \pi(\xi| -) &= \pi(\xi|\bm{\lambda}^{(s+1)}, \tau^{(s+1)})
	\propto \mathcal{L}(\xi) \cdot \log\ \mathcal{N}_{1}(\xi|\mu,\rho^{2}),\\
	\nonumber
	\mathcal{L}(\xi) &=
	\mathcal{V}_{p}(\xi) \mathcal{I}_{(1/2, \infty)}(\xi).
\end{align}
Henceforth, the basic idea is to strictly obey the philosophy of Gibbs sampler: as long as we are to sample $\xi^{(s+1)} \sim \pi(\xi| -)$ (\ref{eq:mu_spec_1}), every latent variables except for the target variable $\xi$ are treated as observed variables, including $\bm{\lambda}^{(s+1)}$ and $\tau^{(s+1)}$.

To start with, we choose a small value of the hyper-parameter $\rho^{2}$  so that the prior part in (\ref{eq:mu_spec_1}), that is, $\pi(\xi) =  \log\ \mathcal{N}_{1}(\xi|\mu,\rho^{2})$, is highly concentrated around its prior mean $\mathbb{E}[\xi] = \exp\ ( \mu + \rho^{2}/2)  \approx \exp\ ( \mu )$. That way, a future state $\xi^{(s+1)}$ is highly probable to be sampled around the value $\exp\ ( \mu )$, leading to an approximate relationship between the future state $\xi^{(s+1)}$ and hyper-parameter $\mu$, described by $\xi^{(s+1)} \approx \exp\ ( \mu )$, or equivalently, $\mu  \approx \log\ \xi^{(s+1)}$. This approximation will be utilized shortly later. 
Throughout this paper, we use $\rho^{2}= 0.001$ as the default value for $\rho^{2}$.

Now, we are in a position to describe how to calibrate the hyper-parameter $\mu$ via the Hill estimator (\ref{eq:Hill_estimator}). We start with ordering the realizations of the $p$ local-scale parameters $\bm{\lambda}^{(s+1)}=(\lambda_{1}^{(s+1)},\cdots,\lambda_{p}^{(s+1)})^{\top}$ to obtain $\lambda_{(1)}^{(s+1)} \geq \cdots\geq\lambda_{(p)}^{(s+1)}$. The Hill estimator based on $\bm{\lambda}^{(s+1)}$ is then
\begin{align}
	\label{eq:Hill_estimator_calibrated}
	\widehat{\xi}_{k}(\bm{\lambda}^{(s+1)}  )
	&=
	\frac{1}{k-1}
	\sum_{j=1}^{k-1}
	\log\ \Bigg( \frac{\lambda_{(j)}^{(s+1)} }{\lambda_{(k)}^{(s+1)} }\Bigg)
	,
	\quad \text{for } k_{\text{L}} \leq k\leq  k_{\text{U}},
\end{align}
where $k_{\text{L}} = \lfloor p/10 \rfloor$ and $k_{\text{U}}
= \lfloor 9p/10 \rfloor$, with $\lfloor \cdot \rfloor$ is the floor function, where $p$ is the number of covariates. In high-dimensional setting, the number of the elements of the set $\{k_{\text{L}},\cdots, k_{\text{U}} \}$ $=$ $\{\lfloor p/10 \rfloor, \cdots, \lfloor 9p/10 \rfloor \}$ $\subset$ $\{ 2, \cdots, p\}$ is still large, approximately, $\lfloor 4p/5 \rfloor $, enough to retain the consistency of the Hill estimator. Note that estimates in (\ref{eq:Hill_estimator_calibrated}) depend on $k$. To eliminate dependency on $k$, first, we average out the Hill estimators (\ref{eq:Hill_estimator_calibrated}) over $k$, and then use the approximation $\mu  \approx \log\ \xi^{(s+1)}$, to get:
\begin{align}
	\label{eq:calibrated_mu}
	\widehat{\mu}(\bm{\lambda}^{(s+1)})
	&=
	\log\
	\{
	\widehat{\xi}(\bm{\lambda}^{(s+1)}  )
	\}
	=
	\log\
	\bigg\{ 
	\frac{1}{k_{\text{U}} - k_{\text{L}} + 1}
	\sum_{k = k_{\text{L}}  }^{k_{\text{U}}}
	\widehat{\xi}_{k}(\bm{\lambda}^{(s+1)}  )
	\bigg\} .
\end{align}

Note that the value of $\widehat{\mu}(\bm{\lambda}^{(s+1)})
$ (\ref{eq:calibrated_mu}) changes at every cycle of the Gibbs sampler, and tuned by $\bm{\lambda}^{(s+1)}$ through the Hill estimator. That is, $\widehat{\mu}$ $(\bm{\lambda}^{(s+1)})
$ can be thought as a calibrated hyper-parameter adapted via the $p$ local-scale realizations $\bm{\lambda}^{(s+1)}$. By replacing $\mu$ with $
\widehat{\mu}(\bm{\lambda}^{(s+1)})$ and substituting $\rho^{2}=0.001$ in the full conditional posterior density $\pi(\xi|-)$ (\ref{eq:mu_spec_1}), the \emph{\textbf{Step 5}} within the Gibbs sampler is tuning-free.

Finally, we explain how to sample from the density $\pi(\xi|-)$ (\ref{eq:mu_spec_1}). For that, first, use a change of variable $\eta = \log\ \xi$ and sample from \begin{align}
	\label{eq:mu_spec_2}
	\eta^{(s+1)} \sim \pi(\eta| -) = \pi(\eta| \bm{\lambda}^{(s+1)}, \tau^{(s+1)})
	\propto 
	\mathcal{L}(\eta)
	\cdot
	\mathcal{N}_{1}(\eta|\widehat{\mu}(\bm{\lambda}^{(s+1)}) ,\rho^{2}=0.001),
\end{align}
where $
\mathcal{L}(\eta)
=
\mathcal{V}_{p}(e^{\eta})
\mathcal{I}_{(\log\ 1/2, \infty)}(\eta)
$
$=$
$
[
\{\Gamma(p/e^{\eta} + 1)\}^{-1}\pi^{p/2} \prod_{j=1}^{p}$ $
(\tau^{(s+1)} + e^{\eta} $ $\lambda_{j}^{(s+1)})^{-(1/e^{\eta} + 1 )}
]$ 
$\mathcal{I}_{(\log\ 1/2, \infty)}(\eta)$. Once we obtain a sample $\eta^{(s+1)}\sim \pi(\eta|-)$, then $\xi^{(s+1)} \sim \pi(\xi|-)$ is obtained via the inverse transformation through $\xi^{(s+1)} = \exp\ \eta^{(s+1)}$.

We use the elliptical slice sampler (ESS) \citep{murray2010elliptical} to sample from $\eta^{(s+1)} \sim \pi(\eta|-)$ (\ref{eq:mu_spec_2}) that exploits the Gaussian prior measure. Conceptually, ESS and Metropolis - Hastings (MH) algorithm are similar in that both comprises two steps: proposal step and criterion step. A difference between the two algorithms arises in the criterion step. If a new candidate does not pass the criterion, then MH takes the current state as the next state: whereas, ESS re-proposes a new candidate until rejection does not take place, rendering the algorithm rejection-free. Further information for ESS is referred to the original paper \citep{murray2010elliptical}. By adopting a jargon from their paper, as the calibrated $\mu$, $\widehat{\mu}(\lambda^{(s+1)})$, is positioned at the center of an ellipse \citep{murray2010elliptical,nishihara2014parallel}, hence, we refer to the following Algorithm \ref{alg:Elliptical slice sampler centered by the Hill estimator} as \emph{elliptical slice sampler centered by the Hill estimator}.
\begin{algorithm}
	\caption{Elliptical slice sampler centered by the Hill estimator}\label{alg:Elliptical slice sampler centered by the Hill estimator}
	\SetAlgoLined
	\textbf{Circumstance : } At the \emph{Step 5} of the $s$-th iteration of the Gibbs sampler in Subsection \ref{subsec:Gibbs sampler}.\\
	\textbf{Input : } Current state $\xi^{(s)}$, and posterior realizations $\bm{\lambda}^{(s+1)}$ and $\tau^{(s+1)}$ obtained from the \emph{Step 3} and \emph{Step 4}, respectively.\\
	\textbf{Output : } A new state $\xi^{(s+1)}$.\\
	1. Calibration of $\mu$: obtain $\widehat{\mu}(\bm{\lambda}^{(s+1)}) = \log\ \{ \widehat{\xi}(\bm{\lambda}^{(s+1)})\}$ (\ref{eq:calibrated_mu}) .\\
	2. Variable change ($\eta = \log \xi$): $\eta^{(s)} = \log\ \xi^{(s)}$.\\
	3. Implement elliptical slice sampler to (\ref{eq:mu_spec_2});\\
	\begin{itemize}
		\item[a. ] Choose ellipse centered by the Hill estimator: $\nu \sim \mathcal{N}_{1}( \widehat{\mu}(\bm{\lambda}^{(s+1)}) , \rho^{2} = 0.001 )$.
		\item[b. ] Define a criterion function: 
		\begin{align*}
			\alpha(\eta, \eta^{(s)}) = \text{min}\{\mathcal{L}(\eta) / \mathcal{L}(\eta^{(s)}),1\}: (\log\ 1/2, \infty) \rightarrow [0,1],
		\end{align*}
		where $\mathcal{L(\eta)}
		=
		[
		\{\Gamma(p/e^{\eta} + 1)\}^{-1}\pi^{p/2} \prod_{j=1}^{p}
		(\tau^{(s+1)} + e^{\eta} \lambda_{j}^{(s+1)})^{-(1/e^{\eta} + 1 )}
		]
		\cdot
		\mathcal{I}_{(\log\ 1/2, \infty)}(\eta)$.
		\item[c. ] Choose a threshold and fix: $u \sim \mathcal{U}[0,1]$.
		\item[d. ] Draw an initial proposal $\eta^{*}$: 
		\begin{align*}
			\theta &\sim \mathcal{U}(-\pi,\pi]\\
			\eta^{*}&=\{\eta^{(s)} - \widehat{\mu}(\bm{\lambda}^{(s+1)})\}\cos\ \theta
			+
			\{\nu - \widehat{\mu}(\bm{\lambda}^{(s+1)})\}\sin\ \theta 
			+
			\widehat{\mu}(\bm{\lambda}^{(s+1)})
		\end{align*}
		\item[e. ] \textbf{if}  \textbf{(} $u < \alpha(\eta^{*}, \eta^{(s)})$  \textbf{)} $\{$ $\eta^{(s+1)} = \eta^{*}$  $\}$ \textbf{else} $\{$ \\
		\indent$\quad$ Define a bracket : $(\theta_{\text{min}},\theta_{\text{max}} ]= (-\pi, \pi ]$.\\
		\indent$\quad$ \textbf{while} \textbf{(} $u \geq \alpha(\eta^{*}, \eta^{(s)})$ \textbf{)} $\{$\\
		\indent$\quad\quad$ Shrink the bracket and try a new point :
		\\
		\indent$\quad\quad$ \textbf{if} \textbf{(} $\theta > 0$ \textbf{)} $\theta_{\text{max}} = \theta$ \textbf{else} $\theta_{\text{min}} = \theta$ \\
		\indent$\quad\quad$ $\theta \sim \mathcal{U}(\theta_{\text{min}},\theta_{\text{max}}]$\\
		\indent$\quad\quad$
		$\eta^{*}=\{\eta^{(s)} - \widehat{\mu}(\bm{\lambda}^{(s+1)})\}\cos\ \theta
		+
		\{\nu - \widehat{\mu}(\bm{\lambda}^{(s+1)})\}\sin\ \theta 
		+
		\widehat{\mu}(\bm{\lambda}^{(s+1)})$
		\\
		\indent$\quad\quad$ $\}$\\
		\indent$\quad\quad$ $\eta^{(s+1)} = \eta^{*}$\\
		\indent$\quad$ $\}$
	\end{itemize}
	4. Variable change ($\xi = e^{\eta}$): $\xi^{(s+1)} = \exp\ \eta^{(s+1)}$.\\
\end{algorithm}

\section{Properties of the Horseshoe}\label{sec:Properties of the Horseshoe}
\subsection{\textbf{Prior analysis of  the Horseshoe}}\label{subsec:Prior analysis of  Horseshoe}
\begin{lemma}[Marginal density and random shrinkage coefficient of the Horseshoe]\label{lemma:hs_prior_prior_analysis}
	\hfill
	\begin{enumerate}
		\item[(a)] Assume $\beta | \lambda, \tau \sim \mathcal{N}_{1}(0, \tau^{2}\lambda^{2})$, $\lambda \sim \mathcal{C}^{+}(0,1)$, and $\tau>0$. Then:
		\begin{align}
			\label{eq:marginal_beta_HS}
			\pi_{\text{HS}}(\beta|\tau) 
			&= \int \mathcal{N}_{1}(\beta | 0, \tau^{2} \lambda^{2}) \pi(\lambda) d\lambda 
			=
			K_{\text{HS}} e^{Z_{\text{HS}}(\beta)}
			E_{1}\{Z_{\text{HS}}(\beta)\},
		\end{align}
		where $K_{\text{HS}} = 1/(\tau  2^{1/2}   \pi^{3/2} )$ and $Z_{\text{HS}}(\beta) = \beta^{2}/(2 \tau^{2})$. $E_{1}(x) = \int_{1}^{\infty}
		e^{-x t} $ $t^{-1} dt$, $x \in \mathbb{R}$, is the exponential integral function. 
		\item[(b)] Assume $\lambda \sim C^{+}(0,1)$, $\kappa = 1/(1 + \tau^{2} \lambda^{2}) \in (0,1)$, and $\tau>0$. Then:
		\begin{align}
			\label{eq:marginal_kappa_HS}
			\pi_{\text{HS}}(\kappa|\tau) &= 
			\frac{\tau}{\pi}
			\cdot
			\frac{\kappa^{-1/2}
				(1-\kappa)^{-1/2}}{1 - (1-\tau^{2})\kappa}. 
		\end{align}
	\end{enumerate}
\end{lemma}
\subsection{\textbf{Proof-- Proposition \ref{corollary:tail_index_of_Horseshoe_prior}}}
\label{sec:Proof-- Restricted tail-heaviness of Horseshoe}
Function $e^{x} E_{1}(x)$ satisfies tight upper and lower bounds \citep{carvalho2010horseshoe};
\begin{align}
	\label{eq:ineq_e_1_function}
	\frac{1}{2}
	\cdot 
	\log\ 
	\bigg(
	\frac{x + 2}{x}
	\bigg)
	<
	e^{x} E_{1}(x)
	< 
	\log\ 
	\bigg(
	\frac{x + 1}{x}
	\bigg)
	,\,\quad x > 0.
\end{align}
Replacing $x$ with $Z_{\text{HS}}(\beta) = \beta^{2}/(2 \tau^{2})$ and multiplying $K_{\text{HS}} = 1/(\tau  2^{1/2}   \pi^{3/2} )$ to the both sides of the inequalities (\ref{eq:ineq_e_1_function}) lead to;
\begin{align}
	\label{eq:ineq_pi(beta|tau)}
	l(\beta)
	< \pi_{\text{HS}}(\beta|\tau)
	<u(\beta),
	\,
	\quad
	\beta \in \mathbb{R}, \quad \tau>0,
\end{align}
where $l(\beta) = (K_{\text{HS}}/2)\cdot \log\ \{ (Z_{\text{HS}}(\beta) + 2)/Z_{\text{HS}}(\beta) \}
$ and $u(\beta) = K_{\text{HS}} \cdot \log\ \{ (Z_{\text{HS}}$ $(\beta) + $  $ 1)/Z_{\text{HS}}(\beta) \}$.

Now, denote the tail (survival) function of the random variable $\beta|\tau$ given $\tau>0$ by $\bar{F}_{\text{HS}}(\beta|\tau) = 1 - F_{\text{HS}}(\beta | \tau)$: then, it holds $(d/d\beta) F_{\text{HS}}(\beta | \tau) = \pi_{\text{HS}}(\beta |\tau)$. Then to show that the tail-index of $\pi_{\text{HS}}(\beta |\tau)$ is $\alpha =1$ for any $\tau>0$, we will prove that the it holds $\lim_{\beta \rightarrow \infty} \bar{F}_{\text{HS}}(c \beta|\tau)/\bar{F}_{\text{HS}}(\beta|\tau) = c^{-1}$ for any $c > 0$ and $\tau >0$. (Because $\pi_{\text{HS}}(\beta |\tau)$ is a symmetric density, showing one-directional limit $\beta \rightarrow \infty$ is sufficient.) By the L’H\^{o}pital’s Rule, it holds 
$\lim_{\beta \rightarrow \infty} \bar{F}_{\text{HS}}(c \beta|\tau)/\bar{F}_{\text{HS}}(\beta|\tau) = c \cdot \lim_{\beta \rightarrow \infty}  \pi_{\text{HS}}(c \beta|\tau)/\pi_{\text{HS}}(\beta|\tau)$, hence, our eventual goal is to prove
\begin{align*}
	\lim_{\beta \rightarrow \infty}  
	\frac{\pi_{\text{HS}}(c \beta|\tau)}{\pi_{\text{HS}}(\beta|\tau)} = c^{-2}, \quad c> 0,\quad \tau >0.
\end{align*}

Now, use inequality (\ref{eq:ineq_pi(beta|tau)}) to upper and lower bound the function $\pi_{\text{HS}}(c \beta|$ $ \tau)/\pi_{\text{HS}}(\beta|\tau)$;
\begin{align}
	\label{eq:horseshoe_fixed_tail_index_pf_1}
	\frac{l(c\beta)}{u(\beta)}
	< 
	\frac{\pi_{\text{HS}}(c \beta|\tau)}{\pi_{\text{HS}}(\beta|\tau)} 
	<
	\frac{u(c\beta)}{l(\beta)} 
	,\,\quad c >0 ,\, \beta\in \mathbb{R}, \quad \tau>0.
\end{align}
First, calculate the limit of the upper bound in the inequality (\ref{eq:horseshoe_fixed_tail_index_pf_1}) at infinity by using 
L’H\^{o}pital’s Rule;
\begin{align*}
	\lim_{\beta \rightarrow \infty} \frac{u(c\beta)}{l(\beta)} &= 
	2 \lim_{\beta \rightarrow \infty} \frac{\log\ \{(Z_{\text{HS}}(c \beta) +1)/Z_{\text{HS}}(c \beta)\} }{\log\ \{(Z_{\text{HS}}(\beta) +2)/Z_{\text{HS}}(\beta)\}} \\ &=
	2 \lim_{\beta \rightarrow \infty} 
	\frac{
		\{Z_{\text{HS}}(c \beta)/(Z_{\text{HS}}(c \beta) + 1)\}
		\cdot
		(-c^{2}/Z_{\text{HS}}(c \beta)^{2})
	}{
		\{Z_{\text{HS}}(\beta)/(Z_{\text{HS}}(\beta) + 2)\}
		\cdot
		(-2/Z_{\text{HS}}(\beta)^{2})
	}\\
	&=
	c^{2}
	\cdot
	\lim_{\beta \rightarrow \infty} 
	\frac{
		Z_{\text{HS}}(\beta) \cdot (Z_{\text{HS}}(\beta) + 2)
	}{
		Z_{\text{HS}}(c \beta) \cdot (Z_{\text{HS}}(c \beta) + 1)
	}
	= c^{-2}, \quad c> 0.
\end{align*}
By the same way, we can show $\lim_{\beta \rightarrow \infty} l(c\beta)/u(\beta) = c^{-2}$, $c >0$. Use the squeeze theorem to the inequality (\ref{eq:horseshoe_fixed_tail_index_pf_1}) to finish the proof.


\section{Prior analysis for the GLT prior}\label{sec:Proof-- Properties of the GLT prior}
\subsection{\textbf{Proof-- Proposition \ref{prop:marginal of glt shrinkage}}}\label{subsec:Proof-- Proposition Marginal of GLT}
\paragraph{$(a)$} Clearly,
\begin{align*}
	\pi(\beta|\tau, \xi) & = 
	\frac{1}{\tau \sqrt{2 \pi}}
	\int_{0}^{\infty}\
	\frac{1}{\lambda}
	\exp\bigg(
	-\frac{\beta^{2}}{2 \lambda^{2}}
	\bigg)
	\bigg(
	1 + \frac{\xi \lambda}{\tau}
	\bigg)^{-(1/\xi + 1)}
	d\lambda.
\end{align*}
Let $x = \xi \lambda/\tau$. Then
\begin{align*}
	\pi(\beta|\tau, \xi)
	&=
	\frac{1}{\tau \sqrt{2 \pi}}
	\int_{0}^{\infty}\
	\exp\bigg(
	-\frac{\beta^{2} \xi^{2}}{2 \tau^{2} x^{2}}
	\bigg)
	x^{-1}(1 + x)^{-(1/\xi + 1)}
	dx,
\end{align*}
or equivalently, for $t = 1/x^{2}$:
\begin{align}
	\label{eq:thm_1_pi(beta)}
	\pi(\beta|\tau, \xi)
	&=
	K
	\int_{0}^{\infty}
	e^{-Z t}
	(t^{1/2})^{-1 + 1/\xi}
	(1 + t^{1/2})^{-(1+1/\xi)}
	dt,
\end{align}
where $K = 1/( \tau 2 ^{3/2} \pi^{1/2})$ and $Z(\beta) = \beta^{2} \xi^{2}/(2 \tau^{2})$. Use $Z = Z(\beta)$ to avoid notation clutter. To utilize the Newton's generalized binomial theorem; 
\begin{align*}
	(x + y)^{r}
	&=
	\sum_{k = 0}^{\infty}
	{r \choose k}
	x^{r-k}
	y^{k}, \quad |x| > |y|, r \in \mathbb{C},
\end{align*}
we divide the integral in (\ref{eq:thm_1_pi(beta)}) into two parts. Then we have
\begin{align}
	\label{eq:thm_1_pi(beta)_2}
	\pi(\beta|\tau, \xi)&=K
	\bigg\{
	\int_{0}^{1}
	e^{-Z t}
	(t^{1/2})^{-1 + 1/\xi}
	(1 + t^{1/2})^{-(1+1/\xi)}
	dt \\
	\nonumber
	&+
	\int_{1}^{\infty}
	e^{-Z t}
	(t^{1/2})^{-1 + 1/\xi}
	(1 + t^{1/2})^{-(1+1/\xi)}
	dt
	\bigg\}.
\end{align}
The first integral of (\ref{eq:thm_1_pi(beta)_2}) is
\begin{align}
	\nonumber
	&\int_{0}^{1}
	e^{-Z t}
	(t^{1/2})^{-1 + 1/\xi}
	(1 + t^{1/2})^{-(1+1/\xi)}
	dt \\
	\nonumber
	&=
	\int_{0}^{1}
	e^{-Z t}
	(t^{1/2})^{-1 + 1/\xi}
	\sum_{k = 0}^{\infty}
	{ -1 - 1/\xi \choose k}
	(t^{1/2}) ^{k}
	dt \\
	\nonumber
	&=
	\sum_{k = 0}^{\infty}
	{ -1 - 1/\xi \choose k}
	\int_{0}^{1}
	e^{-Z t}
	t^{(1 + 1/\xi + k)/2 -1}
	dt\\
	\label{eq:thm1_first_integral}
	&=
	\sum_{k = 0}^{\infty}
	{ -1 - 1/\xi \choose k}
	Z^{-(1 + 1/\xi + k)/2}
	\gamma\{
	(1 + 1/\xi + k)/2
	, Z
	\},
\end{align}
where $\gamma(s, x)
=
\int_{0}^{x}
t^{s-1} e^{-t} dt \quad (s , x \in \mathbb{R})$, is the incomplete lower gamma function.
\\
The second integral of (\ref{eq:thm_1_pi(beta)_2}) is
\begin{align}
	\nonumber
	&\int_{1}^{\infty}
	e^{-Z t}
	(t^{1/2})^{-1 + 1/\xi}
	(1 + t^{1/2})^{-(1+1/\xi)}
	dt \\
	\nonumber
	&=
	\int_{1}^{\infty}
	e^{-Z t}
	(t^{1/2})^{-1 + 1/\xi}
	\sum_{k = 0}^{\infty}
	{ -1 - 1/\xi \choose k}
	(t^{1/2}) ^{-1-1/\xi -k}
	dt \\
	\nonumber
	&=
	\sum_{k = 0}^{\infty}
	{ -1 - 1/\xi \choose k}
	\int_{0}^{1}
	e^{-Z t}
	t^{-1 - k/2}
	dt \\
	\label{eq:thm1_second_integral}
	&=
	\sum_{k = 0}^{\infty}
	{ -1 - 1/\xi \choose k}
	E_{k/2 + 1} (Z),
\end{align}
where $E_{s}(x) = \int_{1}^{\infty}
e^{-x t} t^{-s} dt\quad (s , x \in \mathbb{R})$ is the generalized exponential-integral function of real order \citep{milgram1985generalized,chiccoli1992concerning}. Use ${ -1 - 1/\xi \choose k} = (-1)^{k}{ 1/\xi + k \choose k}$, (\ref{eq:thm1_first_integral}), and (\ref{eq:thm1_second_integral}) to conclude the proof.

\paragraph{$(b)$} Prove by using the change of variable;
\begin{align*}
	\pi(\kappa | \tau, \xi) 
	&= \mathcal{GPD}(\lambda|\tau, \xi)\bigg|_{\lambda = \sqrt{(1-\kappa)/\kappa}}
	\cdot
	\bigg|
	\frac{d\lambda}{dk}
	\bigg|\\
	&=
	\frac{1}{\tau}
	\bigg(
	1
	+ 
	\frac{\xi}{\tau}
	\sqrt{\frac{1-\kappa}{\kappa}}
	\bigg)^{-(1/\xi +1)}
	\cdot
	\frac{1}{2 \kappa^{2}}
	\bigg(
	\frac{1-\kappa}{\kappa}
	\bigg)^{-1/2}
	\\
	&=
	\frac{1}{2 \tau}
	(
	\tau \sqrt{\kappa}
	+ 
	\xi
	\sqrt{1-\kappa}
	)^{-(1/\xi +1)}
	(\tau \sqrt{\kappa})^{1/\xi +1}
	\cdot
	\frac{1}{\kappa^{2}}
	\bigg(
	\frac{1-\kappa}{\kappa}
	\bigg)^{-1/2}
	\\
	&=
	\frac{\tau^{1/\xi} }{2}
	\cdot
	\frac{\kappa^{1/(2\xi) -1  }
		(1 - \kappa)^{-1/2}}{\{
		\tau \kappa^{1/2} + \xi (1-\kappa)^{1/2}
		\}^{(1 + 1/\xi)}}
	.
\end{align*}
\subsection{\textbf{Proof-- Corollary \ref{corollary:separate_roles_in_beta_glt} }}\label{subsec:proof-- corollary:separate_roles_in_beta_glt}
\paragraph{$(a)$} In general, the generalized exponential-integral function has the following property; 
$\lim_{x \rightarrow 0^{+}} E_{1}(x)
= \infty$ and $\lim_{x \rightarrow 0^{+}} E_{s}(x)
= 1/(s-1) \text{ for }s > 1$ \citep{chiccoli1992concerning}.  Using this property, if $k = 0$, then $\lim_{|\beta| \rightarrow 0} \psi_{k=0}^{\text{S}}(\beta) = \lim_{|\beta| \rightarrow 0} E_{1}\{Z(\beta)\} = \infty$ because $Z(\beta) = \beta^{2}\xi^{2}/(2 \tau^{2})$. If $k \in \mathbb{N}$, then $\lim_{\beta \rightarrow 0^{+}} \psi_{k}^{\text{S}}(\beta) =$ $ \lim_{\beta \rightarrow 0^{+}} E_{k/2 + 1}$ $\{Z(\beta)\} $ $ =  2/k < \infty$. 

\paragraph{$(b)$} In general, the incomplete gamma function has the following property; 
$\lim_{x \rightarrow 0^{+}} x^{-a} \cdot \gamma(a,x) = a^{-1} \text{ for } a > 0$ \citep{jameson2016incomplete}. Using this property, $\lim_{|\beta|\rightarrow 0} \psi_{k}^{R} (\beta)
=
\lim_{|\beta|\rightarrow 0} Z(\beta)^{-(1+1/\xi + k)/2}
\cdot
\gamma\{ (1 + 1/\xi +k)/2 , Z(\beta) \} = 2/(1+1/\xi + k) < \infty$ for all $k \in \{0\} \cup \mathbb{N}$. 

\paragraph{$(c)$} In general, the generalized exponential-integral function has the following property; 
$e^{-x}/(x + s) \leq E_{s}(x) \leq e^{-x}/(x + s -1)$ for $x>0$ and $s\geq 1$ \citep{chiccoli1992concerning}. Using this property, we obtain an inequality $e^{-Z(\beta)}/\{Z(\beta) + s\} \leq 
E_{s}(Z(\beta)) \leq e^{-Z(\beta)}/\{Z(\beta) + s -1\}$ for $|\beta| >0 $ and $s \geq 1$. As $|\beta| \rightarrow \infty$, both bounds of $E_{s}(Z(\beta))$ converges to zero with squared exponential rate, and hence, $E_{s}(Z(\beta))$ also do for any $s\geq 1$.

\paragraph{$(d)$} For fixed $k \in \{0\} \cup \mathbb{N}$ and $\xi$, we have $\lim_{|\beta| \rightarrow \infty} \gamma\{ (1 + 1/\xi +k)/2 , Z(\beta) \} = \Gamma((1 + 1/\xi +k)/2)$, where $\Gamma$ is the gamma function, and hence, the function $\gamma\{ (1 + 1/\xi +k)/2 , Z(\beta) \}$ is a slowly varying function \citep{mikosch1999regular}. Using this we can re-express $\psi_{k}^{R} (\beta) = Z(\beta)^{-(1+1/\xi + k)/2}
\cdot
\gamma\{ (1 + 1/\xi +k)/2 , Z(\beta) \}$ by $\psi_{k}^{R} (\beta) = \beta^{-(1+1/\xi + k)} \cdot L(\beta)$, where $L$ is a slowly varying function. This implies that the tail-index of function $\psi_{k}^{R} (\beta)$ is $1+1/\xi + k$.
\subsection{\textbf{Proof-- Lemma \ref{lemma:learnability_of_xi}}}\label{subsec:proof- lemma:learnability_of_xi}
\paragraph{$(a)$} Start with a likelihood part:
\begin{footnotesize}
\begin{align}
	\nonumber
	f(y|\xi) &=\int_{0}^{\infty} \int_{0}^{\infty}\int_{-\infty}^{\infty} \mathcal{N}_{1}(y | \beta,1) \cdot \mathcal{N}_{1}(\beta |0,\lambda^{2})
	\cdot
	\mathcal{GPD}(\lambda|\tau, \xi)
	\cdot
	\mathcal{IG}(1/\xi +1, 1)
	d\beta d\lambda d\tau\\
	\nonumber
	&=
	\int_{0}^{\infty}
	\int_{0}^{\infty}
	\mathcal{N}_{1}(y | 0,1 + \lambda^{2}) 
	\cdot
	\mathcal{GPD}(\lambda|\tau, \xi)
	\cdot
	\mathcal{IG}(\tau|1/\xi +1, 1)
	d\lambda d\tau\\
	\nonumber
	&=
	\frac{1}{\sqrt{2\pi}}
	\int_{0}^{\infty}
	\bigg(
	\int_{0}^{\infty}
	\frac{1}{\sqrt{1+\lambda^{2}}}
	\cdot
	\exp\ \bigg\{-\frac{y^{2}}{2(1+\lambda^{2})}\bigg\}
	\cdot
	\frac{1}{\tau}
	\bigg(1 + \frac{\xi \lambda}{\tau} 
	\bigg)^{-(1/\xi +1)}
	d\lambda
	\bigg)\\
	\nonumber
	&\quad\quad \cdot	
	\mathcal{IG}(\tau|1/\xi +1, 1)
	d\tau\\
	\nonumber
	&\leq
	\frac{1}{\sqrt{2\pi}}
	\int_{0}^{\infty}
	\bigg(
	\int_{0}^{\infty}
	\frac{1}{\sqrt{1+\lambda^{2}}}
	\cdot
	\exp\ \bigg\{-\frac{y^{2}}{2(1+\lambda^{2})}\bigg\}
	\cdot
	\frac{1}{\tau + \xi \lambda}
	d\lambda
	\bigg)
	\cdot
	\mathcal{IG}(\tau|1/\xi +1, 1)
	d\tau\\
	\label{proof:finite_xi_evidence_1}
	&=
	\frac{1}{\sqrt{2\pi}}
	\int_{0}^{\infty}
	\bigg(
	\int_{0}^{1}g(y, \lambda,\tau,  \xi)d\lambda 
	+
	\int_{1}^{\infty}g(y, \lambda,\tau,  \xi)d\lambda
	\bigg)
	\cdot
	\mathcal{IG}(\tau|1/\xi +1, 1)
	d\tau,
\end{align}
\end{footnotesize}
where $g(y, \lambda,\tau,  \xi) =\{1/\sqrt{1+\lambda^{2}}\}
\cdot
\exp\ [-y^{2}/\{2(1+\lambda^{2})\}]
\cdot
\{
1/(\tau + \xi \lambda)
\}$, $y\in \mathbb{R}$ and $\lambda, \tau >0$. Because $g(y, \lambda, \tau,  \xi)$ is continuous on a closed interval $[0,1]$ as a function of $\lambda$, by mean value theorem for integral \citep{bartle2011introduction}, there exists $c \in (0,1)$ such that
\begin{footnotesize}
\begin{align}
	\nonumber
	\int_{0}^{1}g(y, \lambda, \tau, \xi)d\lambda
	&=
	g(y, c,\tau,  \xi)=
	\frac{1}{\sqrt{1+c^{2}}}
	\cdot
	\exp\ \bigg\{-\frac{y^{2}}{2(1+c^{2})}\bigg\}
	\cdot
	\frac{1}{\tau + \xi c}\\
	\label{proof:finite_xi_evidence_2}
	&\leq
	\bigg[
	\frac{1}{\sqrt{1+c^{2}}}
	\exp\ \bigg\{-\frac{y^{2}}{2(1+c^{2})}\bigg\}
	\bigg]
	\cdot
	\frac{1}{\tau}
	=
	A
	\cdot
	\frac{1}{\tau}
	\leq \frac{1}{\tau},\quad \tau \in (0,\infty),
\end{align}
\end{footnotesize}
where $A = A(y,c) = 
\{1/(\sqrt{1+c^{2}})\}\cdot
\exp\ [-y^{2}/\{2(1+c^{2})\}]
$, which is upper bounded by $1$ on $\mathbb{R}\times(0,1)$. Also, we have
\begin{align}
	\nonumber
	\int_{1}^{\infty}g(y, \lambda, \tau, \xi)d\lambda &= 
	\int_{1}^{\infty}
	\frac{1}{\sqrt{1+\lambda^{2}}}
	\cdot
	\exp\ \bigg\{-\frac{y^{2}}{2(1+\lambda^{2})}\bigg\}
	\cdot
	\frac{1}{\tau + \xi \lambda}
	d\lambda \\
	\label{proof:finite_xi_evidence_3}
	&\leq
	\int_{1}^{\infty}
	\frac{1}{\lambda}
	\cdot
	1
	\cdot
	\frac{1}{ \xi \lambda}
	d\lambda
	=
	\int_{1}^{\infty}
	\frac{1}{\lambda^{2}}
	d\lambda
	\cdot
	\frac{1}{ \xi}
	=\frac{1}{ \xi}, \quad \xi \in (1/2, \infty).
\end{align}
Using the upper bounds (\ref{proof:finite_xi_evidence_2}) and (\ref{proof:finite_xi_evidence_3}) to (\ref{proof:finite_xi_evidence_1}), then we have
\begin{footnotesize}
\begin{align*}
	f(y|\xi) &\leq
	\frac{1}{\sqrt{2\pi}}
	\int_{0}^{\infty}
	\bigg(\frac{1}{\tau} + \frac{1}{\xi}\bigg)
	\cdot
	\mathcal{IG}(\tau|1/\xi +1, 1)
	d\tau
	=
	\frac{1}{\sqrt{2\pi}}
	\bigg(
	\int_{0}^{\infty}
	\frac{1}{\tau}
	\cdot
	\mathcal{IG}(\tau|1/\xi +1, 1)
	d\tau 
	+
	\frac{1}{\xi}
	\bigg)\\
	&=
	\frac{1}{\sqrt{2\pi}}
	\bigg\{
	\bigg(
	\frac{1}{\xi} +1
	\bigg)
	+
	\frac{1}{\xi}
	\bigg\}
	=
	\frac{1}{\sqrt{2\pi}}
	\bigg(
	\frac{2}{\xi} +1
	\bigg)
	\leq
	\frac{5}{\sqrt{2\pi}}<\infty,\quad y \in \mathbb{R},\, \xi \in (1/2, \infty).
\end{align*}
\end{footnotesize}
Therefore, trivially for any proper prior $\pi(\xi)$ on $(1/2, \infty)$, we have
\begin{align*}
	m(y)&=
	\int_{1/2}^{\infty}
	f(y|\xi)
	\cdot
	\pi(\xi) d\xi
	\leq
	\frac{5}{\sqrt{2\pi}}
	\int_{1/2}^{\infty}
	\pi(\xi)
	d\xi =
	\frac{5}{\sqrt{2\pi}}
	< \infty
	, \quad y\in \mathbb{R}.
\end{align*}

\paragraph{$(b)$} Under the formulation of the GLT prior (\ref{eq:gpd_beta}) -- (\ref{eq:gpd_xi}), i.e., $\bm{\beta} \sim \pi_{\text{GLT}}(\bm{\beta})$, we have
\begin{footnotesize}
\begin{align}
	\nonumber
	\pi(\xi| - )
	&\propto
	\bigg\{
	\prod_{j=1}^{p}
	\mathcal{GPD}(\lambda_{j}|\tau,\xi)
	\bigg\}
	\cdot
	\mathcal{IG}(\tau|p/\xi + 1, 1)
	\cdot
	\log\ \mathcal{N}_{1}(\xi|\mu, \rho^{2})
	\cdot
	\mathcal{I}_{(1/2,\infty)}(\xi)
	\\
	\nonumber
	&=
	\bigg\{
	\prod_{j=1}^{p}
	\frac{1}{\tau}
	\bigg(
	1 + \frac{\xi\lambda_{j}}{\tau}
	\bigg)^{-(1/\xi + 1)}
	\bigg\}
	\cdot
	\frac{\tau^{-p/\xi - 2}
		e^{-1/\tau}}{\Gamma(p/\xi + 1)}
	\cdot
	\log\ \mathcal{N}_{1}(\xi|\mu, \rho^{2})
	\cdot
	\mathcal{I}_{(1/2,\infty)}(\xi)
	\\
	\nonumber
	&\propto
	\bigg\{
	\tau^{p/\xi}
	\cdot
	\prod_{j=1}^{p}
	(
	\tau + \xi \lambda_{j}
	)^{-(1/\xi + 1)}
	\bigg\}
	\cdot
	\frac{\tau^{-p/\xi - 2}}{\Gamma(p/\xi + 1)}
	\cdot
	\log\ \mathcal{N}_{1}(\xi|\mu, \rho^{2})
	\cdot
	\mathcal{I}_{(1/2,\infty)}(\xi)
	\\
	\nonumber
	&\propto
	\frac{\pi^{p/2}}{\Gamma(p/\xi + 1)}
	\prod_{j=1}^{p}
	(
	\tau + \xi \lambda_{j}
	)^{-(1/\xi + 1)}
	\cdot
	\log\ \mathcal{N}_{1}(\xi|\mu, \rho^{2})
	\cdot
	\mathcal{I}_{(1/2,\infty)}(\xi).
\end{align}
\end{footnotesize}

Now, our goal is to show
\begin{align*}
	m(\bm{\lambda}, \tau)=
	\int_{1/2}^{\infty}
	\frac{\pi^{p/2}}{\Gamma(p/\xi + 1)}
	\prod_{j=1}^{p}
	(
	\tau + \xi \lambda_{j}
	)^{-(1/\xi + 1)}
	\cdot
	\log\ \mathcal{N}_{1}(\xi|\mu, \rho^{2})
	d\xi < \infty,
\end{align*}
where $\bm{\lambda}\in (0,\infty)^{p}, \tau\in(0,\infty)$. 

Let $x = 1/\xi$. Then
\begin{align}
	\nonumber
	m(\bm{\lambda}, \tau)
	&=
	\int_{2}^{0}
	\frac{\pi^{p/2}}{\Gamma(p x  + 1)}
	\prod_{j=1}^{p}
	\bigg(
	\frac{x}{\lambda_{j} + \tau x}
	\bigg)^{x+1}
	\cdot
	\log\ \mathcal{N}_{1}(1/x|\mu, \rho^{2})
	\cdot
	-\frac{1}{x^{2}}
	dx
	\\
	\nonumber
	&=
	\pi^{p/2}
	\cdot
	\int_{0}^{2}
	\frac{(1/\tau)^{p(x+1)}}{\Gamma(p x  + 1)}
	\prod_{j=1}^{p}
	\bigg(
	\frac{\tau x}{\lambda_{j} + \tau x}
	\bigg)^{x+1}
	\cdot
	\log\ \mathcal{N}_{1}(1/x|\mu, \rho^{2})
	\cdot
	\frac{1}{x^{2}}
	dx
	\\
	\label{proof:evdidence_finite_5}
	&\leq
	\pi^{p/2}
	\cdot
	\int_{0}^{2}
	r(x)
	\cdot
	\log\ \mathcal{N}_{1}(1/x|\mu, \rho^{2})
	\cdot
	\frac{1}{x^{2}}
	dx,
\end{align}
where $r(x) = (1/\tau)^{p(x+1)}/\Gamma(p x  + 1)$.
Since $r(x)$ is continuous on a closed interval $[0,2]$, 
there exists $x_{0}\in[0,2]$  such that $r(x_{0}) = \sup_{x\in [0,2]}r(x) = B$. Using this bound $B$ to (\ref{proof:evdidence_finite_5}), we have
\begin{align*}
	m(\bm{\lambda}, \tau)
	&\leq
	\pi^{p/2}
	\cdot
	B
	\cdot
	\int_{0}^{2}
	\log\ \mathcal{N}_{1}(1/x|\mu, \rho^{2})
	\cdot
	\frac{1}{x^{2}}
	dx
	\\
	&\leq
	\pi^{p/2}
	\cdot
	B
	\cdot
	\int_{0}^{\infty}
	\log\ \mathcal{N}_{1}(1/x|\mu, \rho^{2})
	\cdot
	\frac{1}{x^{2}}
	dx
	=
	\pi^{p/2}
	\cdot
	B 
	<
	\infty	,
\end{align*}
where $\bm{\lambda}\in (0,\infty)^{p}, \tau\in(0,\infty)$.

\section{Convergence properties of the GLT posterior}\label{sec:Convergence properties of the GLT posterior}
\subsection{\textbf{Outline of the proofs}}
Outline of the proofs contained in the present section is as follows. First, we prove Lemma \ref{lem:2.1}, \ref{lem:2.2}, and \ref{lem:2.3} in Subsection \ref{subsec:Proofs of lemmas}. In Subsection \ref{subsec:Proofs of theorems}, these lemmas are used to demonstrate Theorem \ref{thm:MSE} and \ref{thm:variance}, followed by Theorem \ref{thm:posterior}, in the main paper. For any function of $\zeta$, $A(\zeta)$ and $B(\zeta)$, we write $A(\zeta) \lesssim B(\zeta)$ to denote that there exists a positive constant $c$ independent of $\zeta$ such that $A(\zeta) \leq c \cdot B(\zeta)$. Adapting the notation to our proofs, the $\zeta$ will be the global-scale parameter $\tau$.


Consider a univariate normal mean model with the GLT prior for a moment, underlying setup for Lemma \ref{lem:2.1}, \ref{lem:2.2}, and \ref{lem:2.3}: that is, $y\sim \mathcal{N}_{1}(\beta,1)$, $\beta|\lambda \sim \mathcal{N}_{1}(0,\lambda^{2})$, and $\lambda\sim \mathcal{GPD}(\tau, \xi)$ with fixed $0<\tau < 1/2$ and $1/2<\xi$. (The sample size $n$ under the multivariate normal mean model will be involved in Theorem \ref{thm:MSE}, \ref{thm:variance}, and \ref{thm:posterior}.) Recall that the quantity $\kappa = 1/(1+\lambda^{2}) \in (0,1)$ is referred to as the random shrinkage coefficient: refer to Proposition \ref{prop:marginal of glt shrinkage} - (b). This is the key ingredient for the derivation procedure whose prior density is given as
\begin{align}
	\label{eq:prior_density_of_RSC}
	\pi(\kappa) = \frac{\tau^{1/\xi}}{2} \frac{\kappa^{1/(2\xi)-1}(1-\kappa)^{-1/2}}{\{\tau \kappa^{1/2}+\xi(1-\kappa)^{1/2}\}^{(1+1/\xi)}}.
\end{align}
On the other hand, the conditional density of the observation $y$ given $\kappa$ is distributed according to Gaussian distribution with mean zero and precision $\kappa$
\begin{align}
	\label{eq:density_of_y_given_kappa}
	f(y\mid \kappa) = 
	\int \mathcal{N}_{1}(y|\beta,1) \mathcal{N}_{1}(\beta|0,1/\kappa-1)d\beta = \mathcal{N}_{1}(y|0,1/\kappa)
	=
	\frac{\kappa^{1/2}}{\sqrt{2\pi}} e^{-\kappa y^2/2}.
\end{align}
By multiplying (\ref{eq:prior_density_of_RSC}) and (\ref{eq:density_of_y_given_kappa}), we can obtain a proportional part of the posterior of the $\kappa$
\begin{align}
	\label{eq:posterior_of_RSC}
	\pi(\kappa \mid y) \propto f(y\mid \kappa) \pi(\kappa) \propto \frac{\tau^{1/\xi}}{2} \frac{\kappa^{1/(2\xi)-1/2}(1-\kappa)^{-1/2}}{\{\tau \kappa^{1/2}+\xi(1-\kappa)^{1/2}\}^{(1+1/\xi)}}
	e^{-\kappa y^2/2}.
\end{align}
Before moving onto the proofs, we shall comment on several reasoning and analytic techniques prevailingly used throughout the derivation procedures:
\begin{itemize}
	\item[(i)] Note that the multiplication of the term $\kappa^{1/(2\xi)-1}$ in the nominator of (\ref{eq:prior_density_of_RSC}) with the term $\kappa^{1/2}$ from (\ref{eq:density_of_y_given_kappa}) leads to the term $\kappa^{1/(2\xi)-1/2}$ in the nominator of (\ref{eq:posterior_of_RSC}), which becomes 1, having specified $\xi=1$. (This suggests that relevant derivation procedure for the Horseshoe does not have this term.) Therefore, we may need to take into account some phase transition of the shape parameter $\xi$ by separately considering two cases $1/2 < \xi\leq 1$ and $\xi>1$. This is because, in the former case, the term $\kappa^{1/(2\xi)-1/2}$ is a bounded function on the interval $[0,1]$, while in the latter case, the term $\kappa^{1/(2\xi)-1/2}$ diverges to $\infty$ as $\kappa \rightarrow 0$. 
	\item[(ii)] Let $d$ is any real-valued function of $\kappa$. We often need to derive the posterior moment of $d(\kappa)$:
	\begin{footnotesize}
	\begin{align}
		\nonumber
		\mathbb{E}[ d(\kappa)|y] &= \int_{0}^{1} d(\kappa) \pi(\kappa|y)d\kappa
		= 
		\int_{0}^{1} d(\kappa) \frac{f(y|\kappa) \pi(\kappa)}{m(y)} d\kappa
		=
		\frac{
			\int_{0}^{1} d(\kappa) f(y|\kappa) \pi(\kappa) d\kappa
		}{\int_{0}^{1} f(y|\kappa) \pi(\kappa)  d\kappa}\\
		\label{eq:denominator_RSC}
		&=
		\frac{
			\int_{0}^{1} d(\kappa) 
			\cdot
			\frac{\kappa^{1/(2\xi)-1/2}(1-\kappa)^{-1/2}}{f_{\tau,\xi}(\kappa)^{(1+1/\xi)}}
			\cdot
			e^{-\kappa y^2/2}
			d\kappa
		}{\int_{0}^{1} 
			\frac{\kappa^{1/(2\xi)-1/2}(1-\kappa)^{-1/2}}{f_{\tau,\xi}(\kappa)^{(1+1/\xi)}}
			\cdot
			e^{-\kappa y^2/2}
			d\kappa}.
	\end{align}
	\end{footnotesize}
	On the last expression (\ref{eq:denominator_RSC}), we denoted $f_{\tau,\xi}(\kappa) = \tau \kappa^{1/2}+\xi(1-\kappa)^{1/2}$ for notational simplicity. Understanding the analytic behavior of the $f_{\tau,\xi}(\kappa)$ greatly reduces calculation burden in deriving lemmas.
	\item[(iii)] We can show that the function $f_{\tau,\xi}(\kappa)$ is strictly concave on closed interval $[0,1]$, and take its global maximum value $\sqrt{\tau^{2} + \xi^{2}}$ at the maximum point $k = 1/\{ 1 + (\xi/\tau)^{2}\}$ belonging to open interval $(0,1)$. Due to the strict concavity and the existence of the maximum point in the open interval $(0,1)$, the minimum point is either $\kappa = 0$ or $1$, which leads to minimum value $\xi$ or $\tau$, respectively. To summarize, the function $f_{\tau,\xi}(\kappa)$ satisfies the following inequalities
	\begin{align}
		\label{eq:lower_and_upper_bounds_of_denom_kappa_2}
		\text{min}(\tau,\xi) \leq f_{\tau,\xi}(\kappa) \leq \sqrt{\tau^{2}+\xi^2}\quad \text{on }[0,1].
	\end{align}
	
	Note that as the $\tau$ gets closer to zero, then the maximum point ($k = 1/\{ 1 + (\xi/\tau)^{2}\}$) and lower bound ($\text{min}(\tau,\xi)$) also get closer to zeros. 
	\item[(iv)] 
	We are eventually interested in the asymptotic behaviors of the GLT estimator and, furthermore, GLT posterior as the $\tau$ converges to zero, with some fixed $\xi>1/2$, to prove Theorem \ref{thm:MSE}, \ref{thm:variance}, and \ref{thm:posterior} in the main paper. Similar proof logic has been taken in \citep{VanDerPas2014}. This logic may be also aligned with showing the minimax optimality of the Lasso estimator \citep{tibshirani1996regression} as the upper bound of the sparsity level goes to zero: refer to Equation (9) in \citep{donoho1992maximum}. Essentially, in Lemma \ref{lem:2.1}, \ref{lem:2.2}, and \ref{lem:2.3}, we will see that such asymptotic behaviors are closely related with the quantification of asymptotic behaviors of relevant posterior moments $\mathbb{E}[ d(\kappa)|y]$ (\ref{eq:denominator_RSC}) as the $\tau$ converges to zero, where $d(\kappa)=\kappa$, $1-\kappa$, or $\kappa^{2}$. Therefore, without loss of generality, we can assume that the $\tau$ is bounded above by any positive real number smaller than 1 which is independent of $\tau$, saying $1/2$, hence, the $\tau$ is confined with $0<\tau<1/2$. 
	
	Three technical uses of this upper bound (that is, `$1/2$' in $0<\tau<1/2$) in proving lemmas are as follows. As the first application, let $\nu$ be any real number in the interval $(0,1)$. Then, the function $f_{\tau,\xi}(\kappa)$ satisfies the following inequality
	\begin{align}
		\label{eq:lower_and_upper_bounds_of_denom_kappa_0nu}
		0 < \xi \sqrt{1-\nu} \leq  f_{\tau,\xi}(\kappa) \leq \sqrt{\tau^{2}+\xi^2} \leq \sqrt{(1/2)^{2}+\xi^2}
		\quad \text{on }[0,\nu].
	\end{align}
	
	Note that the lower bound ($\xi \sqrt{1-\nu}$) and upper bound ($\sqrt{(1/2)^{2}+\xi^2}$) are positive and independent of $\tau$. 
	
	
	As the second application, by using the upper bound of the inequalities (\ref{eq:lower_and_upper_bounds_of_denom_kappa_2}) and $\tau<1/2$, the following inequality holds
	\begin{align}
		\label{eq:lower_and_upper_bounds_of_denom_kappa_01}
		f_{\tau,\xi}(\kappa) \leq \sqrt{\tau^{2}+\xi^2} \leq \sqrt{(1/2)^{2}+\xi^2}\quad \text{on }[0,1].
	\end{align}
	Note that the upper bound ($\sqrt{(1/2)^{2}+\xi^2}$) is positive and independent of $\tau$. 
	
	The following is the third application handling the term $\kappa^{1/(2\xi)-1/2}$ mentioned in (i). Due to the upper bound $\tau<1/2$, it holds $1-\tau^2 > 3/4$. Thus, for any $\xi>1/2$, we have the following lower and upper bounds for the term $\kappa^{1/(2\xi)-1/2}$ on the interval $[1-\tau^{2},1]$:
	
	\begin{equation}\label{eq:lower_bound_of_tricky_term}
		\kappa^{1/(2\xi)-1/2} \geq 
		\begin{cases}
			(1 - \tau^{2})^{1/(2\xi) - 1/2} \geq (3/4)^{1/(2\xi) - 1/2} ,\quad\quad \text{if }1/2 < \xi \leq 1,\\
			1,\quad \quad \quad \quad \quad \quad \quad \quad \quad \quad \quad \quad \quad \quad \quad \quad  \text{if }\xi>1,\\
		\end{cases}
	\end{equation}
	and
	\begin{equation}\label{eq:lower_bound_of_tricky_term_2}
		\kappa^{1/(2\xi)-1/2} \leq
		\begin{cases}
			1,\quad \quad \quad \quad \quad \quad \quad \quad \quad \quad \quad \quad \quad \quad \quad \quad  \text{if }1/2 < \xi \leq 1,\\
			(1 - \tau^{2})^{1/(2\xi) - 1/2} \leq (3/4)^{1/(2\xi) - 1/2} \quad  \quad \quad  \text{if }\xi>1.\\
		\end{cases}
	\end{equation}
	Note that the lower (\ref{eq:lower_bound_of_tricky_term}) and upper (\ref{eq:lower_bound_of_tricky_term_2}) bounds of the term $\kappa^{1/(2\xi)-1/2}$ are independent of $\tau$ for any $\xi>1/2$.

	\item[(v)] Let $\nu$ be any real number in the interval $(0,1)$. Then, by using the kernel expression of the density $\pi(\kappa)$ (\ref{eq:prior_density_of_RSC}), it holds
	\begin{footnotesize}
	\begin{equation}
		\label{eq:kernel_integral}
		\int_{\nu}^{1} \frac{\kappa^{1/(2\xi)-1}(1-\kappa)^{-1/2}}{\{\tau \kappa^{1/2}+\xi(1-\kappa)^{1/2}\}^{(1+1/\xi)}} d\kappa  \leq 
		\int_{0}^{1} \frac{\kappa^{1/(2\xi)-1}(1-\kappa)^{-1/2}}{\{\tau \kappa^{1/2}+\xi(1-\kappa)^{1/2}\}^{(1+1/\xi)}} d\kappa = \frac{2}{\tau^{1/\xi}}.
	\end{equation}
\end{footnotesize}
\end{itemize}
\subsection{\textbf{Proofs of Lemma \ref{lem:2.1}, \ref{lem:2.2}, and \ref{lem:2.3}}\label{subsec:Proofs of lemmas}}
\begin{lemma}\label{lem:2.1}
	Suppose $y\sim \mathcal{N}_{1}(\beta,1)$, $\beta|\lambda \sim \mathcal{N}_{1}(0,\lambda^{2})$, and $\lambda\sim \mathcal{GPD}(\tau, \xi)$ with $0< \tau <1/2 $ and  constant $\xi>1/2 $. Let $T(y)$ be the posterior mean given the GLT prior. Define the thresholding value as $r_{\tau,\xi} = \sqrt{ (2/\xi) \log\left(1/\tau\right)}$. Then $|T(y)-y|$ can be bounded above by a real-valued function $h(y)$ such that for any $\rho>1$,  $h(\cdot)$ satisfies 
	\begin{equation*}
		\lim_{\tau \rightarrow 0} \sup_{|y|> \rho r_{\tau,\xi}} h(y) =0.
	\end{equation*}
	
\end{lemma}
\begin{proof}
	Given any $\rho>1$, choose arbitrary values $\eta,\delta\in (0,1)$ satisfying $\rho^2 >1/\{\eta(1-\delta)\}$. Then, it holds
	\begin{equation}
		\label{eq:pf_eq_1_lem_1}
		|T(y)-y| = |y \mathbb{E}[\kappa \mid y]| \leq |y\mathbb{E}[\kappa \ind_{\{\kappa<\eta\}} \mid y]| +|y\mathbb{E}[\kappa \ind_{ \{\kappa>\eta\}} \mid y]|,\quad \text{for all }y\in \mathbb{R}.		
	\end{equation}
	
	The first equality in (\ref{eq:pf_eq_1_lem_1}) holds since $T(y) = \mathbb{E}[\beta|y] = (1 - \mathbb{E}[\kappa|y])y$ under sparse normal mean model: see Section 1.4 in \citep{carvalho2010horseshoe}.
	
	We focus on the first term on the right-hand side of the inequality (\ref{eq:pf_eq_1_lem_1}), and find some upper bounding function of observation $y$:
	\begin{align}
		\nonumber
		|y\mathbb{E}[\kappa \ind_{\{\kappa<\eta\}} \mid y]| &= \Bigg|
		y \cdot \frac{ \int_{0}^\eta \kappa \frac{\kappa^{1/(2\xi)-{1/2}}(1-\kappa)^{-1/2}}{f_{\tau,\xi}(\kappa)^{(1+1/\xi)}} e^{-\kappa y^2/2}d \kappa}{\int_{0}^1  \frac{\kappa^{1/(2\xi)-{1/2}}(1-\kappa)^{-1/2}}{f_{\tau,\xi}(\kappa)^{(1+1/\xi)}} e^{-\kappa y^2/2}d \kappa} \Bigg|		
		\\
		\label{eq:lemma_pf_1}
		& \lesssim |y|  \frac{ \int_{0}^\eta \kappa^{1/(2\xi)+{1/2}}(1-\kappa)^{-1/2} e^{-\kappa y^2/2}d\kappa}{\int_0^1 \kappa^{1/(2\xi)-{1/2}}(1-\kappa)^{-1/2} e^{-\kappa y^2/2}d\kappa} 
		\\
		\nonumber
		& \leq |y|  \frac{ \int_{0}^\eta \kappa^{1/(2\xi)+{1/2}} e^{-\kappa y^2/2}d\kappa}{\int_0^1 \kappa^{1/(2\xi)-{1/2}}e^{-\kappa y^2/2}d\kappa} \\
		\label{eq:lemma_pf_2}
		&=  |y| \frac{2}{y^2}\frac{ \int_{0}^{\eta y^2/2} t^{1/(2\xi)+{1/2}} e^{-t}dt}{\int_0^{y^2/2} t^{1/(2\xi)-{1/2}}e^{-t}dt}   \\
		\label{eq:lemma_pf_3}
		&  \leq  |y|^{-1}  \frac{ \int_{0}^{\infty} t^{1/(2\xi)+{1/2}} e^{-t}dt}{\int_0^{y^2/2} t^{1/(2\xi)-{1/2}}e^{-t}dt}  \\
		\label{eq:lemma_pf_h1}
		& \lesssim |y|^{-1} \left({\int_0^{y^2/2} t^{1/(2\xi)-{1/2}}e^{-t}dt} \right)^{-1} = h_1(y).
	\end{align}
	
	Here, we used the inequalities (\ref{eq:lower_and_upper_bounds_of_denom_kappa_0nu}) and (\ref{eq:lower_and_upper_bounds_of_denom_kappa_01}) to derive the inequality (\ref{eq:lemma_pf_1}). In (\ref{eq:lemma_pf_2}), we used the change the variable $t=\kappa y^2/2$. The nominator in (\ref{eq:lemma_pf_3}) is the value $\Gamma(1/(2\xi) + 3/2)$ independent of $\tau$. 
	
	We now move to the second term on the right-hand side of the inequality (\ref{eq:pf_eq_1_lem_1}). We can obtain
	\begin{align} 
		\nonumber
		\mathbb{P}(\kappa>\eta \mid y) & = \mathbb{E}[\ind_{\{\kappa>\eta\}}\mid  y]
		\\ 
		\nonumber
		&=\frac{ \int_{\eta}^1 \ \frac{\kappa^{1/(2\xi)-1/2}(1-\kappa)^{-1/2}}{f_{\tau,\xi}(\kappa)^{(1+1/\xi)}} e^{-\kappa y^2/2}d \kappa}{\int_{0}^1  \frac{\kappa^{1/(2\xi)-1/2}(1-\kappa)^{-1/2}}{f_{\tau,\xi}(\kappa)^{(1+1/\xi)}} e^{-\kappa y^2/2}d \kappa} \\
		\nonumber
		&\leq \frac{ \int_{\eta}^1 \ \frac{\kappa^{1/(2\xi)-1/2}(1-\kappa)^{-1/2}}{f_{\tau,\xi}(\kappa)^{(1+1/\xi)}} e^{-\kappa y^2/2}d \kappa}{\int_{0}^{\eta \delta}  \frac{\kappa^{1/(2\xi)-1/2}(1-\kappa)^{-1/2}}{f_{\tau,\xi}(\kappa)^{(1+1/\xi)}} e^{-\kappa y^2/2}d \kappa} \\
		\nonumber		& \leq \frac{ e^{-\eta y^2/2}\int_{\eta}^1 \ \frac{\kappa^{1/(2\xi)-1/2}(1-\kappa)^{-1/2}}{f_{\tau,\xi}(\kappa)^{(1+1/\xi)}} d \kappa}{e^{-\eta \delta y^2/2} \int_0^{\eta \delta}  \ \frac{\kappa^{1/(2\xi)-1/2}(1-\kappa)^{-1/2}}{f_{\tau,\xi}(\kappa)^{(1+1/\xi)}} d\kappa}   \\
		\label{eq:lemma_pf_4}
		&\lesssim
		\exp \left(-\frac{\eta(1-\delta)}{2} y^2\right)
		\cdot
		\int_{\eta}^1 \ \frac{\kappa^{1/(2\xi)-1/2}(1-\kappa)^{-1/2}}{f_{\tau,\xi}(\kappa)^{(1+1/\xi)}} d \kappa
		\\
		\label{eq:lemma_pf_5}&\lesssim \frac{1}{\tau^{1/\xi} }\exp \left(-\frac{\eta(1-\delta)}{2} y^2\right),
	\end{align}
	where we applied the inequalities (\ref{eq:lower_and_upper_bounds_of_denom_kappa_0nu}) and (\ref{eq:kernel_integral}) to derive the inequalities (\ref{eq:lemma_pf_4}) and (\ref{eq:lemma_pf_5}), respectively. Therefore, we have
	\begin{align}
		\label{eq:lemma_pf_h2}
		|y\mathbb{E}[\kappa \ind_{\{\kappa>\eta\}} \mid y]| &\leq |y \mathbb P(\kappa>\eta \mid y)| \lesssim |y|  \tau^{-1/\xi} \exp \left(-\frac{\eta(1-\delta)}{2} y^2\right)= h_2(y).
	\end{align}
	
	Note that the $h_1(y)$ (\ref{eq:lemma_pf_h1}) is an even function, and monotonically decreases on interval $[0,\infty)$, hence, for the $y$ such that $|y|> \rho r_{\tau,\xi}$, we have $ h_1(y) \leq h_1(\rho r_{\tau,\xi})$, followed by $\lim_{\tau \rightarrow 0} \sup_{|y|> \rho r_{\tau,\xi}} h_1(y) = 0$. On the other hand, the $h_2(y)$ (\ref{eq:lemma_pf_h2}) is an even function, and monotonically decreases on the interval $[1/\sqrt{\eta(1-\delta)},\infty)$. Recall that the $\eta$ and $\delta$ are fixed during the calculation and chosen to satisfy the inequality $\rho^2 >1/\{\eta(1-\delta)\}$. As we are interested in asymptotic results with letting $\tau \rightarrow 0$ (hence, the thresholding value $r_{\tau,\xi}= \sqrt{ (2/\xi) \log\left(1/\tau\right)}$ grows to infinity), for sufficiently large $y$ such that $|y|> \rho r_{\tau,\xi}$, it holds inequality $h_{2}(y) \leq h_2(\rho r_{\tau,\xi})$ where 
	\begin{align*}
		h_2(\rho r_{\tau,\xi}) &= \rho \sqrt{\frac{2}{\xi} \log \frac{1}{\tau}} \cdot  \exp \left[
		\left\{ 1 -\eta(1-\alpha)\rho^{2}
		\right\}\cdot
		\frac{1}{\xi} \log\left(\frac{1}{\tau}\right)\right].
	\end{align*}
	
	Because the exponential term eventually dominates the decreasing rate of the $h_2(\rho r_{\tau,\xi})$ as $\tau \rightarrow 0$, it holds $\lim_{\tau \rightarrow 0} \sup_{|y|> \rho r_{\tau,\xi}} h_2(y)=0$. 
	
	To finalize the proof, we define the summation of the two bounding functions as $h(y)=h_1(y)+h_2(y)$, thereby, leading to the inequality $|T(y)-y|\lesssim h(y)$ for all $y\in \mathbb{R}$. Because we showed that it holds $\lim_{\tau \rightarrow 0} \sup_{|y|> \rho r_{\tau,\xi}} h_1(y)=0$ and 
	$\lim_{\tau \rightarrow 0} \sup_{|y|> \rho r_{\tau,\xi}} h_2(y)=0$, it follows $\lim_{\tau \rightarrow 0} \sup_{|y|> \rho r_{\tau,\xi}} h(y) =0$.
\end{proof}

On the next lemma, we provide bounds on posterior moment $\mathbb{E}[1-\kappa \mid y]$ for  $ \xi>1/2 $. 
\begin{lemma} \label{lem:2.2}
	Suppose the conditions in Lemma~\ref{lem:2.1}. Then, for any constant $\eta$ with $0<\eta<3/4 <1-\tau^2$, the following relations hold
	\begin{itemize}
		\item [(a)]	$\mathbb{E}[1-\kappa \mid y]\lesssim e^{y^2/2} \tau^{1/\xi}\left\{ 1/y^{1/\xi +1} + e^{-\eta y^2/2}+e^{-(1-\tau^2) y^2/2} \tau^{2-1/\xi}\right\};	
		$
		\item [(b)] $\mathbb{E}[1-\kappa \mid y]\lesssim e^{y^2/2} \tau^{1/\xi}\left\{ 1+ e^{-\eta y^2/2}+e^{-(1-\tau^2) y^2/2} \tau^{2-1/\xi}\right\}.$
	\end{itemize}

\end{lemma}
\begin{proof}
	By (\ref{eq:denominator_RSC}), we have
	\begin{align}
		\label{eq:lemma2_pf1}
		\mathbb{E}[1-\kappa \mid y] = \frac{ \int_{0}^1  \frac{\kappa^{1/(2\xi)-{1/2}}(1-\kappa)^{1/2}}{f_{\tau,\xi}(\kappa)^{(1+1/\xi)}} e^{-\kappa y^2/2}d \kappa}{\int_{0}^1  \frac{\kappa^{1/(2\xi)-{1/2}}(1-\kappa)^{-1/2}}{f_{\tau,\xi}(\kappa)^{(1+1/\xi)}} e^{-\kappa y^2/2}d \kappa} .
	\end{align}
	
	First, we focus on finding a lower bound of the denominator of $\mathbb{E}[1-\kappa \mid y]$ (\ref{eq:lemma2_pf1})
	\begin{footnotesize}
	\begin{align}
		\nonumber
		\int_{0}^1  \frac{\kappa^{1/(2\xi)-1/2}(1-\kappa)^{-1/2}}{f_{\tau,\xi}(\kappa)^{(1+1/\xi)}} e^{-\kappa y^2/2}d \kappa &\geq \int_{1-\tau^2}^1  \frac{\kappa^{1/(2\xi)-1/2}(1-\kappa)^{-1/2}}{f_{\tau,\xi}(\kappa)^{(1+1/\xi)}} e^{-\kappa y^2/2}d \kappa  \\
		\label{eq:lemma2_pf2}
		&\gtrsim \tau^{-1-1/\xi} \int_{1-\tau^2}^1  {\kappa^{1/(2\xi)-1/2}(1-\kappa)^{-1/2}} e^{-\kappa y^2/2}d \kappa \\
		\label{eq:lemma2_pf3}
		&\gtrsim    \tau^{-2-1/\xi}   \int_{1-\tau^2}^1   e^{-\kappa y^2/2}d \kappa \\
		\nonumber
		& \geq  e^{- y^2/2} \tau^{-1/\xi}.
	\end{align}
	\end{footnotesize}
	Here, to derive the inequality (\ref{eq:lemma2_pf2}), we used the inequality $f_{\tau,\xi}(\kappa) = \tau \kappa^{1/2}+\xi(1-\kappa)^{1/2}\leq (1+\xi)\tau$ on $[1-\tau^{2},1]$. To derive the inequality (\ref{eq:lemma2_pf3}), we used the inequality $(1-\kappa)^{-1/2}\geq \tau^{-1}$ on $[1-\tau^{2},1]$ and the inequality (\ref{eq:lower_bound_of_tricky_term}) on the term $\kappa^{1/(2\xi)-1/2}$.
	
	Next, we concern the nominator of the $\mathbb{E}[1-\kappa \mid y]$ (\ref{eq:lemma2_pf1}), and find its upper bound. For notational simplicity, we express the integrand of the the nominator with $g(\kappa)$: that is, 
	\begin{align*}
		g(\kappa)  = \frac{\kappa^{1/(2\xi)-{1/2}}(1-\kappa)^{1/2}}{f_{\tau,\xi}(\kappa)^{(1+1/\xi)}} e^{-\kappa y^2/2}
	\end{align*}		
	
	We separate the integral $\int_{0}^1  	g(\kappa)d \kappa$ with the following three parts
	\begin{align}
		\label{eq:sum}
		\int_{0}^1  	g(\kappa)d \kappa =\int_{0}^{\eta}g(\kappa)d \kappa + \int_{\eta}^{1-\tau^2} g(\kappa)d \kappa+ \int_{1-\tau^2}^1  	g(\kappa)d \kappa .
	\end{align}
	
	The first integral on the right hand side of (\ref{eq:sum}) has two upper bounds, which leads to the statements (a) and (b), respectively. To derive the statement (a), we use the following upper bound
	\begin{align}
		\label{eq:first_integral_pf1}	
		\int_{0}^{\eta}  	g(\kappa) d\kappa &\lesssim  \int_{0}^{\eta}  \kappa^{1/(2\xi)-1/2}e^{-\kappa y^2/2} d\kappa \\
		\label{gamma}	 &= \int_{0}^{\eta y^2/2} \frac{2}{y^2}\left(\frac{2t}{y^2}\right)^{1/(2\xi)-1/2}e^{-t} dt \\
		\nonumber	 &\lesssim \frac{1}{y^{1/\xi+1}} \int_{0}^{\eta y^2/2} t^{1/(2\xi)-1/2} e^{-t} dt \\
		\nonumber	 &\leq \frac{1}{y^{1/\xi+1}} \Gamma(1/(2\xi)+1/2) \\
		\nonumber	 & \lesssim \frac{1}{y^{1/\xi+1}},
	\end{align}
	where the inequality (\ref{eq:first_integral_pf1}) holds due to inequalities (\ref{eq:lower_and_upper_bounds_of_denom_kappa_0nu}) and $(1-\kappa)^{1/2} \leq 1$ on $[0,\eta]$. We used the change of variable $t= \kappa y^2/2$ in the equation ~\eqref{gamma}.
	
	On the other hand, the statement (b) is based on the following upper bound:
	\begin{align}
		\nonumber	\int_{0}^{\eta}  	g(\kappa) d\kappa &\lesssim  \int_{0}^{\eta}  \kappa^{1/(2\xi)-1/2}e^{-\kappa y^2/2} d\kappa \leq \int_{0}^{\eta}  \kappa^{1/(2\xi)-1/2} d\kappa \lesssim 1.
	\end{align}
	
	An upper bound for the second term on the right hand side of the equality (\ref{eq:sum}) can be obtained as follows:
	\begin{align}
		\label{eq:lem2_pf1}
		\int_{\eta}^{1-\tau^2} 	g(\kappa) d\kappa &\lesssim  		
		\int_{\eta}^{1-\tau^2} \kappa^{1/(2\xi)-{1/2}} (1-\kappa)^{-1/(2\xi)} e^{-\kappa y^2/2} d\kappa \\
		\label{eq:lem2_pf2}
		& 	\lesssim  		
		\int_{\eta}^{1-\tau^2} (1-\kappa)^{-1/(2\xi)} e^{-\kappa y^2/2} d\kappa \\
		\nonumber
		& 
		\leq  e^{-\eta y^2/2}  \int_{\eta}^{1-\tau^2} (1-\kappa)^{-1/(2\xi)} d\kappa   \\
		\nonumber
		&=  e^{-\eta y^2/2}  \cdot 1/\{1 - 1/(2\xi)\} \cdot ( (1-\eta)^{1-1/(2\xi)} - \tau^{1-1/(2\xi)} )
		\\
		\nonumber
		&\lesssim e^{-\eta y^2/2}.
	\end{align}
	
	Here, we used the inequality $f_{\tau,\xi}(\kappa) \geq \xi (1-\kappa)^{1/2}$ to derive the inequality (\ref{eq:lem2_pf1}). In the integrand of (\ref{eq:lem2_pf1}), we can show that term $\kappa^{1/(2\xi)-{1/2}}$ on the interval $[\eta,1 - \tau^{2}]$ is bounded above by (i) $1$ when $1/2 < \xi \leq 1$ and by (ii) $\eta^{1/(2\xi) - 1/2}$ when $\xi>1$: therefore, the inequality (\ref{eq:lem2_pf2}) holds.
	
	Lastly, an upper bound for the third term on the right hand side of the equality (\ref{eq:sum}) can be obtained as follows:
	\begin{align}
		\label{eq:lem2_thrid_int_pf1}
		\int_{1-\tau^2}^1 	g(\kappa) d\kappa &\lesssim 
		\tau^{-1-1/\xi} \cdot \int_{1-\tau^2}^1  \kappa^{1/(2\xi)-{1/2}}  (1-\kappa)^{1/2} e^{-\kappa y^2/2} d\kappa \\
		\label{eq:lem2_thrid_int_pf2}
		&\lesssim \tau^{-1-1/\xi} \cdot e^{-(1-\tau^2) y^2/2} \int_{1-\tau^2}^1 (1-\kappa)^{1/2}d\kappa  \\
		\nonumber
		& = \tau^{-1-1/\xi} \cdot e^{-(1-\tau^2) y^2/2} \cdot (2/3) \tau^{3} \\
		\nonumber
		&\lesssim \tau^{2-1/\xi}  \cdot e^{-(1-\tau^2) y^2/2}.
	\end{align}
	
	The inequality (\ref{eq:lem2_thrid_int_pf1}) holds since it holds $f_{\tau,\xi}(\kappa) \geq \tau \cdot (3/4)^{1/2}$ on the interval $[1-\tau^{2},1]$. To derive the inequality (\ref{eq:lem2_thrid_int_pf2}), we used the upper bounds of the term $\kappa^{1/(2\xi)-1/2}$ (\ref{eq:lower_bound_of_tricky_term_2}).
	
	Finish the proofs by aggregating the derived results for the lower bound of the denominator and the upper bound of the nominator of the moment $\mathbb{E}[1-\kappa \mid y]$ (\ref{eq:lemma2_pf1}).
	
\end{proof}
\begin{lemma}\label{lem:2.3}
Suppose the condition in Lemma~\ref{lem:2.1}. Let $\mbox{Var}[\beta \mid y]$ be the posterior variance. Then, the following relations hold		
\begin{itemize}
	\item[(a)] $\mbox{Var}[\beta \mid y] \leq 1+y^2$ for all $y \in \mathbb{R}$ ;
	\item[(b)] $\mbox{Var}[\beta \mid y]$ can be bounded above by a real-valued function $\bar h(y)$ such that for any $\rho>1$,  $\bar h(y)$ satisfies $		\lim_{\tau \rightarrow 0} \sup_{|y|> \rho r_{\tau,\xi}} \bar h(y) =1.$
\end{itemize}
\end{lemma}

\begin{proof}
Conditional posterior of the $\beta$ given $\kappa,y$ is $\beta \mid \kappa,y \sim \mathcal{N}_{1}((1-\kappa)y,1-\kappa)$. Thus, the following relation holds for all values $y\in \mathbb{R}$
\begin{align}
	\nonumber	\mbox{Var}[\beta \mid y] &= \mathbb{E} [\mbox{Var}[\beta \mid \kappa ,y ]\mid y] + \mbox{Var} [\mathbb{E}[\beta \mid \kappa,y]\mid y] \\
	\nonumber	&=\mathbb{E}[1-\kappa \mid y] + \mbox{Var}[(1-\kappa) y \mid y] \\
	\label{eq:thm2_later_use}&= \mathbb{E}[1-\kappa \mid y] + y^2 \mathbb{E}[\kappa^2 \mid y] - y^2 \left( \mathbb{E}[\kappa \mid y]\right)^{2} \\
	& \label{eq:lemma3_pf1}
	\leq 1 + y^2 \mathbb{E}[\kappa^2 \mid y],
\end{align}
where the inequality (\ref{eq:lemma3_pf1}) holds since $\mathbb{E}[1-\kappa\mid y] = T(y)/y  \leq 1$ for all $y\in \mathbb{R}$.

The statement (a) is then immediately followed by using an upper bound of the second term in (\ref{eq:lemma3_pf1}): $y^2 \mathbb{E}[\kappa^2 \mid y] \leq y^{2}$.

We now prove the statement (b). This follows the same reasoning adopted in Lemma~\ref{lem:2.1}. That is, we focus on finding some upper bounding function of the term $y^2\mathbb{E}[\kappa^2  |y]$ in (\ref{eq:lemma3_pf1}), particularly, letting $\tau$ converging to zero. Given any $\rho>1$, choose arbitrary values $\eta,\delta\in (0,1)$ satisfying $\rho^2 >1/\{\eta(1-\delta)\}$. Then we can separate the term $y^2\mathbb{E}[\kappa^2  |y]$ as
$$y^2 \mathbb{E}[\kappa^2 \mid y] = y^2 \mathbb{E}[\kappa^2 \ind_{\{\kappa\leq \eta\}} \mid y] + y^2 \mathbb{E}[\kappa^2 \ind_{\{\kappa>\eta\}} \mid y],$$ where the first and second terms of the right hand side are further bounded by
\begin{align*}
	y^2 \mathbb{E}[\kappa^2\ind_{\{\kappa\leq \eta\}} \mid y] \lesssim y^{-2}\left( {\int_0^{y^2/2} t^{1/(2\xi)-1}e^{-t}dt} \right)^{-1} ,\\
	y^2 \mathbb{E}[\kappa^2 \ind_{\{\kappa>\eta\}} \mid y] \lesssim y^2 \tau^{-1/\xi} \exp \left(-\frac{\eta(1-\delta)}{2} y^2\right).
\end{align*}

By combining the above results, we can show that there exists a function $\tilde{h}(y)$ satisfying the inequality $y^2\mathbb{E}[\kappa^2 \mid y] \lesssim \tilde{h}(y)$ for all $y\in \mathbb{R}$, and the limiting equation $\lim_{\tau \rightarrow 0} \sup_{|y|> \rho r_{\tau,\xi}} \tilde{h}(y) =0.$ Denote $\bar{h}(y) = 1+ \tilde{h}(y)$: then, it holds $	\mbox{Var}[\beta \mid y] \lesssim 1 + \tilde{h}(y) = \bar{h}(y)$ for all $y\in \mathbb{R}$ and $\lim_{\tau \rightarrow 0} \sup_{|y|> \rho r_{\tau,\xi}} \bar{h}(y) =1.$.
\end{proof}

\subsection{\textbf{Proofs of Theorem \ref{thm:MSE}, \ref{thm:variance}, and \ref{thm:posterior}}\label{subsec:Proofs of theorems}}
\paragraph{Proof-- Theorem \ref{thm:MSE}}
We start with separating the total MSE with the summation of the nonzero and zero means parts
\begin{align*}
	\mathbb{E}_{\bm{\beta}_0}\|T(\textbf{y})-\bm{\beta}_0\|_{2}^2 &= \sum_{i\in S} \mathbb{E}_{\beta_{0i}} \left[ (T(y_i)-\beta_{0i})^2\right] + \sum_{i\in S^c} \mathbb{E}_{\beta_{0i}}\left[ (T(y_i)-\beta_{0i})^2\right],
\end{align*}
with the support $S=\{i: \beta_{0i}\ne 0\} \subset \{1,2,\cdots,n \}$ with the cardinality $|S|=q$.\\

\noindent  $\bullet$ \emph{Nonzero mean}

For the coefficient index $i \in S$, we have
\begin{equation*}
	\begin{aligned}
		\mathbb{E}_{\beta_{0i}} \left [ ( T(y_i)-\beta_{0i})^2\right]  &= \mathbb{E}_{\beta_{0i}} \left[ (T(y_i)-y_i+y_i-\beta_{0i})^2\right]  \\
		& \leq 2\mathbb{E}_{\beta_{0i}} \left   [(T(y_i)-y_i )^2\right] +2\mathbb{E}_{\beta_{0i}} \left[(y_i-\beta_{0i})^2\right]  \\
		&= 2\mathbb{E}_{\beta_{0i}}\left [ (T(y_i)-y_i)^2\right] +2,
	\end{aligned}
\end{equation*}
where the last equation holds due to the assumption of the unit standard deviation for the normal mean model (that is, $\mathbb{E}_{\beta_{0i}} \left[(y_i-\beta_{0i})^2\right] =1$). 

Now, we shall split the first term of the last equation into two terms by using the values $\pm \rho r_{\tau,\xi} = \pm \rho \sqrt{ (2/\xi) \log\left(1/\tau\right)}$ as the cut-off values on the real line
\begin{footnotesize}
\begin{align}
	\nonumber
	\mathbb{E}_{\beta_{0i}}   \left[ (T(y_i)-y_i  )^2 \right] &=  \mathbb{E}_{\beta_{0i}} \left[ (T(y_i)-y_i )^2 \ind_{\{|y_i|\leq \rho r_{\tau,\xi}\}} \right] + \mathbb{E}_{\beta_{0i}}  \left[ (T(y_i)-y_i )^2 \ind_{\{|y_i|>\rho r_{\tau,\xi}\}} \right]\\
	\label{eq:thm_1_pf1}
	& \lesssim 
	\mathbb{E}_{\beta_{0i}} \left[ y_i^2 \ind_{\{|y_i|\leq \rho r_{\tau,\xi}\}} \right] + \mathbb{E}_{\beta_{0i}}  \left[ h(y_{i})^2 \ind_{\{|y_i|>\rho r_{\tau,\xi}\}} \right]\\
	\nonumber
	&\lesssim \rho^2 r_{\tau,\xi}^2 + \left(\sup_{|y|> \rho r_{\tau,\xi}} h(y) \right)^2\\
	\nonumber
	& = \rho^2 r_{\tau,\xi}^2 +o(1)\\
	& \lesssim  r_{\tau,\xi}^2.
\end{align}
\end{footnotesize}
Here, we applied the inequality $|T(y_i)-y_i| \leq |y_i|$ and Lemma \ref{lem:2.1} to derive the inequality (\ref{eq:thm_1_pf1}). \\



\noindent $\bullet$ \emph{ Zero means}	

For the coefficient index $i \in S^{c}$, we have:
\begin{align}
	\label{eq:thm_1_pf2}
	\mathbb{E}_{\beta_{0i}} \left[T(y_i)^2\right] = 	\mathbb{E}_{\beta_{0i}} \left[T(y_i)^2 \ind_{ \{|y_i| \leq r_{\tau,\xi}\}}\right]+\mathbb{E}_{\beta_{0i}} \left[T(y_i)^2 \ind_{ \{|y_i| > r_{\tau,\xi}\}}\right].
\end{align}
As for the first term on the right hand side of (\ref{eq:thm_1_pf2}), we have the following upper bound:
\begin{align}
	\nonumber
	\mathbb{E}_{\beta_{0i}} \left[T(y_i)^2 \ind_{\{|y_i| \leq r_{\tau,\xi}\}}\right] &  = \mathbb{E}_{\beta_{0i}} \left[ \left( \mathbb{E}[1 - \kappa \mid y_{i}] y_{i}\right) ^{2}
	\ind_{\{|y_i| \leq r_{\tau,\xi}\}}\right]\\
	\label{eq:thm_1_pf3}
	&\lesssim   \tau^{2/\xi} \int_0^{r_{\tau,\xi}} y^2\exp(y^2)\cdot  \frac{1}{\sqrt{2\pi}}\exp(-y^2/2) dy  \\
	\nonumber
	& \lesssim  \tau^{2/\xi} \cdot	 \int_0^{r_{\tau,\xi}} y^2 e^{y^{2}/2}dy  \\
	\label{pf:zero_means_1}
	& = \tau^{2/\xi}
	\cdot	\left(
	r_{\tau,\xi} \exp(r_{\tau,\xi}^2/2) 
	-\int_{0}^{r_{\tau,\xi}} e^{y^{2}/2}dy \right)
	\\
	\label{eq:thm_1_pf4}
	&\lesssim \tau^{1/\xi}  \sqrt{\log \frac{1}{\tau}}.
\end{align}	

The inequality (\ref{eq:thm_1_pf3}) holds due to Lemma 2.2 (b), that is, $\mathbb{E}[1-\kappa \mid y_{i}]\lesssim e^{y_{i}^2/2} \tau^{1/\xi}$ $\{ 1+ e^{-\eta y_{i}^2/2}+ $ $e^{-(1-\tau^2) y_{i}^2/2} \tau^{2-1/\xi}\}$, where the last term vanishes as $\tau$ converges to zero (that is, $e^{-(1-\tau^2) y^2/2} \tau^{2-1/\xi}=o(1)$). We used the integration by parts to obtain the equality (\ref{pf:zero_means_1}). 

In addition, let $\phi(x)$ and $\Phi(x)$ denote the probability density and cumulative distribution functions for the standard Gaussian random variable: then, we can show that the second term on the right hand side of (\ref{eq:thm_1_pf2}) is upper bounded by the same term (\ref{eq:thm_1_pf4})
\begin{align}
	\label{eq:thm_1_pf5} 
	\mathbb{E}_{\beta_{0i}} \left[T(y_i)^2 \ind_{ \{|y_i| > r_{\tau,\xi}\}}\right] &\leq  2 \int_{r_{\tau,\xi}}^{\infty} y^2 \phi(y) dy \\
	&= 2 \int_{r_{\tau,\xi}}^{\infty} \phi(y) - \frac{d}{dy} [y\phi(y)] dy \label{eq10_(i)}\\
	\nonumber	&= 2(1-\Phi(r_{\tau,\xi})) +2 r_{\tau,\xi} {\phi(r_{\tau,\xi})} \\
	\label{eq10_(ii)}
	&\lesssim \frac{\phi(r_{\tau,\xi})}{r_{\tau,\xi}} + r_{\tau,\xi} \phi (r_{\tau,\xi}) \\	
	\nonumber	&= \frac{1}{r_{\tau,\xi} } \frac{e^{-r_{\tau,\xi}^2/2}}{\sqrt{2\pi}} + r_{\tau,\xi} \frac{e^{-r_{\tau,\xi}^2/2}}{\sqrt{2\pi}}\\
	\nonumber	
	&\lesssim r_{\tau,\xi} e^{-r_{\tau,\xi}^2/2} (1+o(1))
	\\
	\label{eq10_(iii)}
	& \lesssim \tau^{1/\xi}  \sqrt{\log \frac{1}{\tau}} ,
\end{align}
where the inequality (\ref{eq:thm_1_pf5}), equality \eqref{eq10_(i)}, and inequality (\ref{eq10_(ii)}) are based on the inequality $|T(y_{i})|\leq |y_{i}|$ for $y_{i}\in \mathbb{R}$, and the identity $y^2\phi(y)=\phi(y)-d[y\phi(y)]/(dy)$, and the Mill's ratio (see Equation (29) of \citep{VanDerPas2014}), respectively.	 
\\

\noindent  $\bullet$ \emph{Conclusion}	

We demonstrated that (i) for the $q$ number of nonzero means, it holds $\mathbb{E}_{\beta_{0i}} $ $ [ (T(y_i)-{\beta}_{0i})^2] \lesssim \rho^{2} r_{\tau,\xi}^2$; and (ii) for the $n-q$ number of zero means, it holds $\mathbb{E}_{\beta_{0i}} [ T(y_i)^2] \lesssim \tau^{1/\xi}\sqrt{\log(1/\tau)}$ as $\tau$ converges to zero. Finish the proof by summing up the results.\\

\paragraph{Proof-- Theorem \ref{thm:variance}}
As similar with the proof of Theorem~\ref{thm:MSE}, we consider the decomposition	
\begin{equation*}
	\mathbb{E}_{\bm{\beta}_0} \sum_{i=1}^{n} \mbox{Var} [\beta_i \mid y_i]= 	\sum_{i \in S} \mathbb{E}_{\beta_{0i}} \left[ 
	\mbox{Var}[\beta_i \mid y_i] \right] +	\sum_{i \in S^c} \mathbb{E}_{\beta_{0i}}  \left[\mbox{Var}[\beta_i \mid y_i]\right],
\end{equation*}
with the support $S=\{i: \beta_{0i}\ne 0\} \subset \{1,2,\cdots,n \}$ with the cardinality $|S|=q$.\\

\noindent  $\bullet$ \emph{Nonzero means}		

For the coefficient index $i \in S$, we have 
\begin{footnotesize}
\begin{align}
	\nonumber	\mathbb{E}_{\beta_{0i}} \left[\mbox{Var}[\beta_i \mid y_i]\right] &= \mathbb{E}_{\beta_{0i}} \left[\mbox{Var}[\beta_i \mid y_i]\ind_{\{|y_i|\leq \rho r_{\tau,\xi}\}}\right] + \mathbb{E}_{\beta_{0i}} \left[\mbox{Var}[\beta_i \mid y_i]\ind_{\{|y_i|> \rho r_{\tau,\xi}\}}\right] \\
	\label{eq:nonzero2}
	& \lesssim \mathbb{E}_{\beta_{0i}} \left[(1 + y_{i}^{2}) \ind_{\{|y_i|\leq \rho r_{\tau,\xi}\}}\right] + \mathbb{E}_{\beta_{0i}} \left[\bar{h}(y_{i})\ind_{\{|y_i|> \rho r_{\tau,\xi}\}}\right] \\
	& \nonumber	
	\lesssim 1+ \rho^2 r^2_{\tau,\xi} + \sup_{|y|> \rho r_{\tau,\xi}} \bar h(y) \\
	\nonumber	&\lesssim 2+r_{\tau,\xi}^2,\\
	\nonumber	&\lesssim r_{\tau,\xi}^2,
\end{align}
\end{footnotesize}
where the the inequality ~\eqref{eq:nonzero2} is obtained by applying Lemma~\ref{lem:2.3} -(a) and (b) to the two separated terms. Because the term $r_{\tau,\xi}^2$ diverges to infinity as $\tau \rightarrow 0$, the last expression is valid. \\	

\noindent  $\bullet$ \emph{Zero means}		

For the coefficient index $i \in S$, we again work with the decomposition
\begin{align}
	\nonumber
	\mathbb{E}_{\beta_{0i}} &\left[\mbox{Var}[\beta_i \mid y_i]\right] \\
	\label{eq:thm2_zero_pf1}
	&= \mathbb{E}_{\beta_{0i}} \left[\mbox{Var}[\beta_i \mid y_i]\ind_{\{|y_i|\leq  \rho r_{\tau,\xi}\}}\right] + \mathbb{E}_{\beta_{0i}} \left[\mbox{Var}[\beta_i \mid y_i]\ind_{\{|y_i|>  \rho r_{\tau,\xi}\}}\right].
\end{align}

The second term on the right hand side of the (\ref{eq:thm2_zero_pf1}) is bounded above as follow:
\begin{align}
	\mathbb{E}_{\beta_{0i}} \left[\mbox{Var}[\beta_i \mid y_i]\ind_{\{|y_i|> \rho r_{\tau,\xi}\}}\right] &\leq 2 \int_{\rho r_{\tau,\xi}}^{\infty} (1+y)^2 \phi(y) dy  \nonumber \\
	& \leq 4 \int_{\rho r_{\tau,\xi}}^{\infty} (  y^{2} +1) \phi(y) dy  \label{eq42_1} \\
	&\lesssim  r_{\tau,\xi} e^{-\rho^2 r^2_{\tau,\xi}/2} \label{eq42_2}\\
	& \lesssim \tau^{\rho^2/\xi} \sqrt{\log \frac{1}{\tau}} \nonumber \\
	& \leq \tau^{1/\xi} \sqrt{\log \frac{1}{\tau}} \label{eq42_3}.
\end{align}
Here, we used the Cauchy-Schwarz inequality ($(1 +y)^{2} \leq 2(1+y^{2})$ for all $y\in \mathbb{R}$) to derive the inequality~\eqref{eq42_1}. Derivation of the inequality~\eqref{eq42_2} follows the same reasoning used in obtaining the upper bound for the zeros means in  Theorem~\ref{thm:MSE}. The last inequality~\eqref{eq42_3} holds since $\tau^{\rho^{2}/\xi} \leq \tau^{1/\xi}$ on the interval $0 < \tau < 1/2$ for any value $\rho>1$.

Now, we focus on the first term on the right hand side of (\ref{eq:thm2_zero_pf1}). To that end, we first derive some tight upper bound of the posterior variance $\mbox{Var}[\beta \mid y]$ (\ref{eq:thm2_later_use}), splitted with the following two terms:

\begin{footnotesize}
\begin{align*}
	\mbox{Var}[\beta \mid y] &=\mathbb{E}[1-\kappa \mid y] - y^{2}\left( \mathbb{E}[\kappa \mid y]\right)^{2} + y^2 \mathbb{E}[\kappa^2 \mid y] \\
	&=\mathbb{E}[1-\kappa \mid y] - y^{2} \left(\mathbb{E}[1-\kappa \mid y]-1\right)^{2} + y^2 \mathbb{E}[\kappa^2 \mid y]\\
	&= \mathbb{E}[1-\kappa \mid y]  - y^2 \left(\mathbb{E}[1-\kappa \mid y]\right)^2 -y^2 +2y^2 \mathbb{E}[1-\kappa \mid y]+ y^2 \mathbb{E}[\kappa^2 \mid y] \\
	&= \mathbb{E}[1-\kappa \mid y]  - y^2 \left(\mathbb{E}[1-\kappa \mid y] \right)^2 \\
	& \quad  + y^2 \mathbb{E}[1-\kappa \mid y]+y^2 \mathbb{E}[1-\kappa \mid y]+ y^2 \left(\mathbb{E}[\kappa^2 \mid y] - 1\right) \\
	& \leq \mathbb{E}[1-\kappa \mid y]  +  y^2 \mathbb{E}[1-\kappa \mid y] (1-\mathbb{E}[1-\kappa \mid y] ) + y^2 \mathbb{E}[1-\kappa \mid y] \\
	&\leq \mathbb{E}[1-\kappa \mid y]  +2  y^2 \mathbb{E}[1-\kappa \mid y],
\end{align*}
\end{footnotesize}
where the first and second inequalities follow from inequalities $y^2 \left(\mathbb{E}[\kappa^2 \mid y] - 1\right)$ $  \leq 0$ and $\mathbb{E}[1-\kappa \mid y] \geq 0$, respectively. Thus, the first term on the right hand side of (\ref{eq:thm2_zero_pf1}) can be re-expressed as
\begin{align}
	\label{eq:thm3_final_pf1}
	\mathbb{E}_{\beta_{0i}} &\left[\mbox{Var}[\beta_i \mid y_i]\ind_{\{|y_i|\leq  \rho r_{\tau,\xi}\}}\right] \lesssim	
	\\
	\nonumber
	&\mathbb{E}_{\beta_{0i}} \left[\mathbb{E}[1-\kappa \mid y_{i}] \ind_{\{|y_i|\leq  \rho r_{\tau,\xi}\}}\right] 
	+
	\mathbb{E}_{\beta_{0i}} \left[y_{i}^2 \mathbb{E}[1-\kappa \mid y_{i}]\ind_{\{|y_i|\leq  \rho r_{\tau,\xi}\}}\right].
\end{align}

By using Lemma~\ref{lem:2.2}-(b), we can derive an upper bound of the the first term on the right hand side of (\ref{eq:thm3_final_pf1}) as follows: 
\begin{footnotesize}
\begin{align*}
	\mathbb{E}_{\beta_{0i}} \left[\mathbb{E}[1-\kappa \mid y_{i}] \ind_{\{|y_i|\leq  \rho r_{\tau,\xi}\}}\right] 
	&\lesssim
	\int_{0}^{\rho r_{\tau,\xi}}
	e^{y^2/2} \tau^{1/\xi}
	\left\{ 1+ e^{-\eta y^2/2}+o(1)\right\}
	\frac{1}{\sqrt{2\pi}}e^{-y^{2}/2}dy\\
	&\lesssim
	\tau^{1/\xi}  \sqrt{\log \frac{1}{\tau}},
\end{align*}	
\end{footnotesize}

and by using Lemma~\ref{lem:2.2}-(a), we can obtain an upper bound of the second term on the right hand side of (\ref{eq:thm3_final_pf1}) as follows
\begin{footnotesize}
\begin{align*}
	&\mathbb{E}_{\beta_{0i}} \left[y_{i}^2 \mathbb{E}[1-\kappa \mid y_{i}]\ind_{\{|y_i|\leq  \rho r_{\tau,\xi}\}}\right]\\
	&\quad \lesssim
	\int_{0}^{\rho r_{\tau,\xi}}
	y^{2}
	e^{y^2/2} \tau^{1/\xi}
	\left\{y^{- 1-1/\xi}+ e^{-\eta y^2/2}+o(1)\right\}
	\frac{1}{\sqrt{2\pi}}e^{-y^{2}/2}dy
	\\
	&\quad\lesssim
	\tau^{1/\xi}
	\int_{0}^{\rho r_{\tau,\xi}}
	\left\{
	y^{2}
	e^{-\eta y^2/2}
	+
	y^{1-1/\xi}
	\right\}
	dy
	\\
	&\quad = 
	\tau^{1/\xi}
	\frac{1}{2 - 1/\xi}
	\left(
	\rho r_{\tau,\xi}
	\right)^{2 - 1/\xi }
	+
	\tau^{1/\xi}
	\frac{\sqrt{2}}{\eta^{3/2}}
	\gamma
	\left(
	\frac{3}{2},
	\frac{ \eta \rho^{2} r_{\tau,\xi}^{2}}{4}
	\right) \\
	&\quad\lesssim \tau^{1/\xi} \left(\log \frac{1}{\tau} \right)^{1-1/(2\xi)},
\end{align*}
\end{footnotesize}
where we wrote the negligible term $e^{-(1-\tau^2) y^2/2} \tau^{2-1/\xi}$ as $o(1)$.\\



	%
	%
	
	\noindent  $\bullet$ \emph{Conclusion}	
	
	Finalize the proof by summing up all results for the nonzero and zero means.

	\paragraph{Proof-- Theorem \ref{thm:posterior}}
	Set the $\tau$ to be a sequence of $n$ by $\tau = (q/n)^\alpha $, with the $\xi$ satisfying $1/2<\xi<\alpha$ for any constant $\alpha>1/2$. Let $\mathcal{U}_{\tau,\xi} ={q} \log (1/\tau)/\xi +(n-q) \tau^{{1/\xi}} \{\log {(1/\tau)}\}^{\max\{\frac{1}{2},1-\frac{1}{2\xi}\}} $ denote the upper bound of the inequality (\ref{eq:thm_2_results}) in the main paper. By plugging the $\tau = (q/n)^\alpha$ in the $\mathcal{U}_{\tau,\xi}$, we can bound the ratio between the $\mathcal{U}_{\tau= (q/n)^\alpha,\xi}$ and the minimax optimal rate $q \log (n/q)$:
	\begin{footnotesize}
	\begin{equation*}
		\frac{\mathcal{U}_{\tau= (q/n)^\alpha,\xi}}{q \log (n/q)} \lesssim \frac{q\log(n/q)+ q \{\log(n/q)\}^{\max\{\frac{1}{2},1-\frac{1}{2\xi}\}} (q/n)^{\alpha/\xi-1}  }{q\log(n/q)} \lesssim 1+(q/n)^{\alpha/\xi-1} \lesssim 1,
	\end{equation*}
\end{footnotesize}
	since $\{\log(n/q)\}^{\max\{\frac{1}{2},1-\frac{1}{2\xi}\}} \leq \log(n/q)$ given $\xi >1/2$ and $(q/n)^{\alpha/\xi-1}<1$ given $\alpha/\xi >1$. (Here, the notation $A \lesssim B$ is interpreted as $A(n) \leq c B(n)$ where $c$ is independent of $n$.) Thus, there exists some constant $c$ independent of $n$ such that $0 < \mathcal{U}_{\tau= (q/n)^\alpha,\xi} \leq c  \cdot q \log (n/q)$. Furthermore, as $c$ is positive, it holds $ K  \cdot \mathcal{U}_{\tau= (q/n)^\alpha,\xi} \leq   q \log (n/q)$ with $K=1/c$.
	
	By Markov's inequality, we have 
	\begin{footnotesize}
	\begin{align*}
		\Pi\left(\|\bm{\beta}-\bm{\beta}_0\|_2^2 \geq  M_n q\log(n/q)\mid \textbf{y} \right) &\leq 		\Pi\left(\|\bm{\beta}-\bm{\beta}_0\|_2^2 \geq M_n K \mathcal{U}_{\tau= (q/n)^\alpha,\xi} \mid \textbf{y} \right) \\
		&\leq \frac{\mathbb E[\|\bm{\beta}-\bm{\beta}_0\|_2^2 \mid \textbf{y}]}{ M_n  K \mathcal{U}_{\tau= (q/n)^\alpha,\xi}} \\
		&\leq \frac{2 \mathbb E[\|T(\textbf{y})-\bm{\beta}_0\|_2^2 \mid \textbf{y}]}{ M_n K \mathcal{U}_{\tau= (q/n)^\alpha,\xi}}+\frac{2 \mathbb E[\|\bm{\beta}-T(\textbf{y})\|_2^2 \mid \textbf{y}]}{ M_n K \mathcal{U}_{\tau= (q/n)^\alpha,\xi}}
		\\
		& = \frac{2 \|T(\textbf{y})-\bm{\beta}_0\|_2^2 }{ M_n K \mathcal{U}_{\tau= (q/n)^\alpha,\xi}}+\frac{2 
			\sum_{i=1}^{n} \mbox{Var}[\beta_i \mid y_i] 			
		}{ M_n K \mathcal{U}_{\tau= (q/n)^\alpha,\xi}},
	\end{align*}
\end{footnotesize}
	for any sequence $M_n$ with $M_n \rightarrow \infty$. 
	
	Now, take the expectation $\mathbb E_{\bm \beta_0}[\cdot]$ on the both sides of the above inequality to get
	\begin{footnotesize}
	\begin{align}
		\label{eq:final_final_final}
		\mathbb E_{\bm \beta_0}\left[
		\Pi\left(\|\bm{\beta}-\bm{\beta}_0\|_2^2 \geq  M_n q\log(n/q)\mid \textbf{y} \right)\right]
		&\leq 
		\frac{2 
			\mathbb E_{\bm \beta_0}\left[
			\|T(\textbf{y})-\bm{\beta}_0\|_2^2\right] }{ M_n K \mathcal{U}_{\tau= (q/n)^\alpha,\xi}}+\frac{2 
			\mathbb E_{\bm \beta_0}\left[
			\sum_{i=1}^{n} \mbox{Var}[\beta_i \mid y_i] 			\right]}{ M_n K \mathcal{U}_{\tau= (q/n)^\alpha,\xi}}\\
		\nonumber
		&\lesssim 
		\frac{4}{M_{n} K},
	\end{align} 	
	\end{footnotesize}
	where we used Theorem~\ref{thm:MSE} and ~\ref{thm:variance} on the first and second terms on the right hand side of the inequality (\ref{eq:final_final_final}).

	\subsection{\textbf{Derivation of theoretical optimal values for} $\xi$ \textbf{inducing near-minimax rate}}\label{subsec:Derivation of theoretical optimal values for xi inducing near-minimax rate}
	
	Denote the upper bound of the inequality (\ref{eq:thm_2_results}) in the main paper as $\mathcal{U}_{\tau,\xi} ={q} \log (1/\tau)/\xi +(n-q) \tau^{{1/\xi}} \{\log {(1/\tau)}\}^{\max\{\frac{1}{2},1-\frac{1}{2\xi}\}} $. With the choice $\tau=1/n$, we wish to find values $\xi$ belongs to the set defined as
	\begin{footnotesize}
	\begin{align}
		\label{eq:optimal_xi}
		\Omega= \{ \xi \in (1/2,\infty) \mid  \mathcal{U}_{\tau=1/n,\xi}
		= q/\xi  \cdot\log n  + (n-q) (1/n)^{{1/\xi}} \{\log {n}\}^{\max\{\frac{1}{2},1-\frac{1}{2\xi}\}}
		\lesssim q \log n \}.
	\end{align}
	\end{footnotesize}
	
	Elements of the set $\Omega$ can be regarded as optimal values for the $\xi$ in the sense that, with the choice $\tau=1/n$, the upper bound $\mathcal{U}_{\tau,\xi}$ is bounded above by the ``near-minimax rate'' $q \log n $ \citep{van2017adaptive}. Such values $\xi$ may be useful when the underlying sparsity level is unknown, and provide some insights about the behaviour of the $\xi$ under diverse sparse regime.
	
	Noting from the first term of the $\mathcal{U}_{\tau=1/n,\xi}$ (\ref{eq:optimal_xi}), it holds $q/\xi  \cdot\log n \lesssim q \log n$ for any $\xi>1/2$. This implies that the first term of $\mathcal{U}_{\tau=1/n,\xi}$ do not restrict on the choice of $\xi$. Now, we concern the second term of $\mathcal{U}_{\tau=1/n,\xi}$ by finding values $\xi$ satisfying
	\begin{align*}
		&(n-q) (1/n)^{{1/\xi}} \{\log {n}\}^{\max\{\frac{1}{2},1-\frac{1}{2\xi}\}} \lesssim q \log n
		\\
		&\quad \Longleftrightarrow
		n^{1-1/\xi} \{\log {n}\}^{\max\{\frac{1}{2},1-\frac{1}{2\xi}\}} \lesssim q \log n
		\\
		&\quad \Longleftrightarrow
		n^{1-1/\xi}\lesssim q \{ \log n \}^{\min\{\frac{1}{2},\frac{1}{\xi}\}} 
		\\
		&\quad\Longleftrightarrow
		(1-1/\xi) \log n \lesssim 
		\log q 
		+
		\min\{1/2,1/\xi\} \log \log n\\
		&\quad\Longleftrightarrow
		\xi
		\lesssim \frac{\log n}{\log n - \log q - \min\{1/2,1/\xi\} \log \log n}.
	\end{align*}
	In conclusion, the optimal set for $\xi$ is characterized by
	\begin{align*}
		\Omega = \left\{ \xi \in (1/2, \infty) \mid 
		\xi
		\lesssim \frac{\log n}{\log n - \log q - \min\{1/2,1/\xi\} \log \log n }
		\right\}.
	\end{align*}
Since $\min\{1/2, 1/\xi\} \log \log n > 0$ is nearly constant as $n \rightarrow \infty$, and by the definition of the sparsity level $s = q/n$ (Equation (\ref{eq:sparsity_level_gene_expression_data}) in the main paper), the values $\xi = \log n / (\log n - \log q) = - \log n / \log s$ (up to a multiplicative constant) are optimal values for $\xi$ (i.e., elements of $\Omega$) that induce the near-minimax rate.

\end{appendix}

%

\begin{supplement}
	\stitle{Program Codes}
	\sdescription{}
\end{supplement}
\begin{supplement}
	\stitle{Supplementary Material to  ``Tail-adaptive Bayesian shrinkage''}
	\sdescription{}
\end{supplement}
\begin{acks}[Acknowledgments]
The authors would like to thank the anonymous referees, an Associate Editor and the Editor for their constructive comments that significantly improved the quality of this paper.
\end{acks}

\bibliographystyle{imsart-number} 
\bibliography{paper-ref}       

\end{document}